\documentclass[11pt,twoside]{article}
\newcommand{\ifAMS}[2]{#1}   
\newcommand{\ifau}[4]{#1}  
\newcommand{\ifbook}[2]{#1}   
\newcommand{\ifunivariate}[2]{#1} 
\usepackage{ifthen}
\usepackage[ampersand]{easylist}
\usepackage{amsmath,amssymb,amsthm}
\usepackage{mathtools}
\usepackage{natbib}
\usepackage{epsfig,graphicx}
\usepackage{comment}
\usepackage{color}
\usepackage{srcltx}
\usepackage[mathscr]{eucal}
\usepackage[math]{easyeqn}
\usepackage{etoolbox}
\usepackage{hyperref}
\hypersetup{
            colorlinks,
            linkcolor=hookersgreen,
            linktoc=blue,
            citecolor=ultramarine,
            urlcolor=black,
            filecolor=black
            }
\usepackage{nameref}
\usepackage{multirow}
\usepackage{mathrsfs}
\usepackage{bbm}
\usepackage{algorithm}
\usepackage{algpseudocode}

\numberwithin{equation}{section}
\numberwithin{figure}{section}
\newcounter{example}[section]
\numberwithin{example}{section}
\newcounter{remark}[section]
\numberwithin{remark}{section}
\newtheorem{theorem}{Theorem}[section]
\newtheorem{proposition}[theorem]{Proposition}
\newtheorem{lemma}[theorem]{Lemma}
\newtheorem{corollary}[theorem]{Corollary}

\newtheorem{exmp}[example]{Example}
\newtheorem{rmrk}[remark]{Remark}
\newenvironment{example}{\begin{exmp}\rm}{\end{exmp}}
\newenvironment{remark}{\begin{rmrk}\rm}{\end{rmrk}}

\ifbook{
    \newcommand{\Chapter}[1]{\section{#1}}
    \newcommand{\Section}[1]{\subsection{#1}}
    \newcommand{\Subsection}[1]{\subsubsection{#1}}

  }
  {
    \newcommand{\Chapter}[1]{\chapter{#1}}
    \newcommand{\Section}[1]{\section{#1}}
    \newcommand{\Subsection}[1]{\subsection{#1}}

  }

\numberwithin{equation}{section}
\numberwithin{figure}{section}

\renewcommand{\(}{$\,}
\renewcommand{\)}{\,$}

\def\nquad{\hspace{-1cm}}
\def\eqdef{\stackrel{\operatorname{def}}{=}}

\DeclareMathAlphabet{\mathbbmsl}{U}{bbm}{bx}{sl}

\newcommand{\cc}[1]{\mathscr{#1}}
\newcommand{\bb}[1]{\boldsymbol{#1}}

\renewcommand{\bar}[1]{\overset{\!\_\!\_\!\_}{#1}}

\renewcommand{\tilde}[1]{\widetilde{#1}}

\newcommand{\thankstitle}[1]{\ifthenelse{\equal{#1}{}}{}{\thanks{#1}}}
\newcommand{\thanksau}[1]{\ifthenelse{\equal{#1}{}}{}{\thanks{#1}}}

\newcommand{\aua}[6]
{\def\authora{#1}
\def\runauthora{#2}
\def\addressa{#3}
\def\emaila{#4}
\def\affiliationa{#5}
\def\thanksa{#6}}

\def\theauthors{
\ifau{ 
  \author{
    \authora
    \thanksau{\thanksa}
    \\[5.pt]
    \addressa \\
    \texttt{ \emaila}
  }
}
{  
  \author{
    \authora
    \thanksau{\thanksa}
    \\[5.pt]
    \addressa \\
    \texttt{ \emaila}
    \and
    \authorb
    \thanksau{\thanksb}
    \\[5.pt]
    \addressb \\
    \texttt{ \emailb}
  }
}
{   
  \author{
    \authora
    \thanksau{\thanksa}
    \\[5.pt]
    \addressa \\
    \texttt{ \emaila}
    \and
    \authorb
    \thanksau{\thanksb}
    \\[5.pt]
    \addressb \\
    \texttt{ \emailb}
    \and
    \authorc
    \thanksau{\thanksc}
    \\[5.pt]
    \addressc \\
    \texttt{ \emailc}
  }
} {   
  \author{
    \authora
    \thanksau{\thanksa}
    \\[5.pt]
    \addressa \\
    \texttt{ \emaila}
    \and
    \authorb
    \thanksau{\thanksb}
    \\[5.pt]
    \addressb \\
    \texttt{ \emailb}
    \and
    \authorc
    \thanksau{\thanksc}
    \\[5.pt]
    \addressc \\
    \texttt{ \emailc}
    \and
    \authord
    \thanksau{\thanksd}
    \\[5.pt]
    \addressd \\
    \texttt{ \emaild}
  }
}
}

\renewcommand{\Gamma}{\varGamma}
\renewcommand{\Pi}{\varPi}
\renewcommand{\Sigma}{\varSigma}
\renewcommand{\Delta}{\varDelta}
\renewcommand{\Lambda}{\varLambda}
\renewcommand{\Psi}{\varPsi}
\renewcommand{\Phi}{\varPhi}
\renewcommand{\Theta}{\varTheta}
\renewcommand{\Omega}{\varOmega}
\renewcommand{\Xi}{\varXi}
\renewcommand{\Upsilon}{\varUpsilon}

\def\argmax{\operatornamewithlimits{argmax}}
\def\argmin{\operatornamewithlimits{argmin}}

\def\av{\bb{a}}

\def\uv{\bb{u}}

\def\wv{\bb{w}}
\def\xv{\bb{x}}

\def\zv{\bb{z}}

\def\Av{\bb{A}}

\def\Xv{\bb{X}}
\def\Yv{\bb{Y}}

\def\gammav{\bb{\gamma}}

\def\thetav{\bb{\theta}}

\def\xiv{\bb{\xi}}

\def\sumi{\sum_{i=1}^{n}}

\usepackage{color}

\definecolor{blue(pigment)}{rgb}{0.2, 0.2, 0.6}
\definecolor{ultramarine}{rgb}{0.07, 0.04, 0.56}
\definecolor{darkspringgreen}{rgb}{0.09, 0.45, 0.27}
\definecolor{hookersgreen}{rgb}{0.0, 0.44, 0.0}
\definecolor{plum(traditional)}{rgb}{0.56, 0.27, 0.52}
\definecolor{purple(html/css)}{rgb}{0.5, 0.0, 0.5}
\definecolor{magenta(dye)}{rgb}{0.79, 0.08, 0.48}

\def\HU{\HV}

\def\xvs{\xv^{*}}

\def\amax{\nu}

\def\xvb{\bar{\xv}}

\def\Rem{\mathcal{R}}

\def\fba{\bar{f}}

\def\vH{\vA}

\def\Idd{\Gamma}

\def\feta{\phi}

\def\HUDe{\Delta}
\def\Remi{\mathcal{Q}}

\def\DV{\mathsf{D}}

\def\HU{\mathsf{D}}

\def\dltw{\delta}

\def\dltwa{\alpha}
\def\dltwb{\omega}

\def\dltwm{\kappa}

\def\dltwu{\tau}

\def\II{\mathcal{I}}

\def\R{\mathbbmsl{R}}
\def\E{\mathbbmsl{E}}

\def\P{\mathbbmsl{P}}

\def\kappa{\varkappa}

\def\diag{\operatorname{diag}}

\def\Fr{\operatorname{Fr}}

\def\ND{\mathcal{N}}

\def\oper{\operatorname{op}}

\def\Var{\operatorname{Var}}

\def\T{\top}
\def\tr{\operatorname{tr}}

\def\TV{\operatorname{TV}}

\def\CONST{\mathtt{C} \hspace{0.1em}}
\def\CONSTi{\mathtt{C}}
\def\cond{\, \big| \,}

\def\nsize{{n}}

\def\sumi{\sum_{i=1}^{\nsize}}

\def\ex{\mathrm{e}}

\def\Id{I\!\!\!I}
\def\Ind{\operatorname{1}\hspace{-4.3pt}\operatorname{I}}

\def\bias{\mathsf{b}}

\def\BB{I\!\!B}     
\def\BB{B}

\def\BBB{\cc{B}}

\def\BBH{W}

\def\DP{D}

\def\DPGP{\DP_{\GP}}

\def\DPt{\tilde{\DPc}}
\def\DPGPt{\DPt_{\GP}}

\def\DPt{\tilde{\DP}}

\def\dist{d}

\def\dimH{\dimA}

\def\dimp{p}

\def\dimA{\mathtt{p}}

\def\dimG{\dimA_{\GP}}

\def\dimq{q}

\def\err{\diamondsuit}
\def\errs{\err_{\rdomega}^{*}}

\def\eps{\epsilon}			
\def\eps{\varepsilon}

\def\gaussv{\bb{\gauss}}
\def\gauss{\gamma}

\def\GP{G}

\def\IF{\mathbbmsl{F}}

\def\kullb{\cc{K}} 

\def\LL{\cc{L}}

\def\priord{\pi}

\def\QP{Q}

\def\regrf{m}

\def\regrfv{\bb{\regrf}}

\def\riskt{\cc{R}}

\def\rr{\mathtt{r}}

\def\supA{\lambda}

\def\thetav{\bb{\theta}}
\def\thetavs{\thetav^{*}}

\def\Tau{T}

\def\uvc{\uv^{c}}

\def\UV{\mathcal{U}}

\def\vA{\mathtt{v}}

\def\vtheta{\vartheta}
\def\vthetav{\bb{\vtheta}}

\def\weight{w}

\def\wv{\bb{w}}

\def\xvd{\xv^{\circ}}
\def\xx{\mathtt{x}}

\def\xvd{\xv^{\circ}}

\def\XX{\cc{X}}

\def\zq{z}

\def\zz{\mathfrak{z}}

\renewcommand{\Chapter}[1]{\section{#1}}

\def\HVB{\mathcal{V}}

\def\lgd{f}
\def\elll{\ell}
\def\lgdL{\elll}
\def\smlc{\rho}

\def\thetitle{Dimension free non-asymptotic bounds on the accuracy of high dimensional Laplace approximation}
\def\theruntitle {Laplace approximation and effective dimension}

\def\theabstract{
This note attempts to revisit the classical results on Laplace approximation in a modern non-asymptotic and dimension free form. 
Such an extension is motivated by applications to high dimensional statistical and optimization problems.
The established results provide explicit non-asymptotic bounds on the quality of a Gaussian approximation 
of the posterior distribution in total variation distance in terms of the so called \emph{effective dimension} \( \dimG \).
This value is defined as interplay between information contained in the data and in the prior distribution.
In the contrary to prominent Bernstein - von Mises results, the impact of the prior is not negligible and 
it allows to keep the effective dimension small or moderate even if the true parameter dimension is huge or infinite.
We also address the issue of using a Gaussian approximation with inexact parameters with the focus on 
replacing the Maximum a Posteriori (MAP) value by the posterior mean and design the algorithm of Bayesian optimization based on Laplace iterations.
The results are specified to the case of nonlinear inverse problem.
}

\def\kwdp{62F15}
\def\kwds{60E15; 65C05; 65K10}

\def\thekeywords{Gaussian approximation, effective and critical dimension, posterior mean}

\def\thankstitle{The author thanks Richard Nickl for very valuable remarks leading to a substantial improvement of the paper.}
\def\thethanks{  
  Financial support by the German Research Foundation (DFG) through the Collaborative Research Center 1294 ``Data assimilation'' is gratefully acknowledged.  
}

\aua
{Vladimir Spokoiny}
{Spokoiny, V. }
{
    Weierstrass Institute and HU Berlin,  
    \\
    Mohrenstr. 39, 10117 Berlin, Germany
    }
{spokoiny@wias-berlin.de}
{Weierstrass-Institute and Humboldt University Berlin}
{\thankstitle
}

\begin{document}
\thispagestyle{empty}
{
\title{\thetitle \thanks{\thethanks}}
\theauthors

\maketitle
\begin{abstract}
{\footnotesize \theabstract}
\end{abstract}

\ifAMS
    {\par\noindent\emph{MSC Subject Classification:} Primary \kwdp. Secondary \kwds}
    {\par\noindent\emph{JEL codes}: \kwdp}

} 

\tableofcontents

\Chapter{Introduction}

High dimensional Laplace approximation has recently gained an increasing attention in connection with 
Bayesian inference for complicated nonlinear parametric models such as nonlinear inverse problems and Deep Neuronal Networks.
Laplace approximation is obtained by replacing a log density with its second order Taylor approximation around the point of maximum.
This leads to a Gaussian measure centered at the maximum with a covariance corresponding to the Hessian of the negative log-density (see, e.g., \cite[Section 4.4]{bishop2006pattern}).
The asymptotic behavior of the parametric Laplace approximation in the small noise or large data limit has been studied extensively in the past (see, e.g., \cite{wong2001asymptotic}).
The asymptotic approximation of general integrals of the form \( \int e^{\lambda f(x)} \, g(x) \, d x \) by Laplace's method is presented in \cite{olver1974asymptotics,wong2001asymptotic}.
Non-asymptotic error bounds for the Laplace approximation can be found in \cite{Olv:1968} 
for the univariate case and in \cite{MW:1983,IM:2014} for the multivariate case. 
\cite{lapinski2019multivariate} studied 
the Laplace approximation error and its convergence in the limit \( \lambda \to \infty \) in the multivariate case when the function \( \lgd \) depends on \( \lambda \).
Coefficients in the asymptotic expansion of the approximated integral are given in \cite{nemes2013explicit}.

Laplace approximation is  an important step in establishing the prominent Bernstein - von Mises (BvM) Theorem
that quantifies the convergence of the scaled posterior distribution toward a Gaussian distribution in the large data or small noise limit.
Parametric BvM theory is well-understood \cite{van2000asymptotic, le2012asymptotic}. 
Modern applications with a high dimensional parameter space and limited sample size pose new questions and identify new issues in study 
of applicability and accuracy of Laplace approximation.
We refer to \cite{Lu:2017} for a study a parametric BvM theorem for nonlinear Bayesian inverse problems with an increasing number of parameters. 
A number of papers discuss the BvM phenomenon for nonlinear inverse problems; see e.g. \cite{nickl2020, monard2019efficient, giordano2020bernstein}, where the convergence is quantified in a distance that metrizes the weak convergence.
\cite{SSW:2020} showed that the Laplace approximation error in Hellinger distance converges to zero in the order of the noise level.
The recent paper \cite{HeKr2022} provides a finite sample error of Laplace approximation for the total variation distance
with an explicit dependence on the dimension and of the nonlinearity of the forward mapping. 
The Laplace approximation is also widely utilized for different purposes in computational Bayesian statistics; see e.g. 
\cite{rue2009approximate}.


\Section{Motivation: Gaussian approximation of the posterior}
As one of the main motivation for this study, consider the problem of Bayesian inference 
for the log-likelihood function \( L(\thetav) = L(\Yv,\thetav) \) 
with data \( \Yv \), a parameter \( \thetav \in \R^{\dimp} \)
and a Gaussian prior \( \priord \sim \ND(\thetav_{0}, \GP^{-2}) \).
Here \( \GP^{-2} \) is a symmetric positive definite matrix in \( \R^{\dimp} \).
Then the posterior density \( \priord_{\GP}(\cdot) \) of \( \thetav \) given \( \Yv \) is proportional to the product 
\( \ex^{L(\thetav)} \, \ex^{-\| \GP(\thetav - \thetav_{0}) \|^{2}/2} \):
\begin{EQA}
	\vthetav_{\GP} \cond \Yv
	\sim
	\priord_{\GP}(\thetav)
	& \propto &
	\exp\bigl\{ L(\thetav) - \| \GP(\thetav - \thetav_{0}) \|^{2}/2 \bigr\} ,
\label{ewfhuhw2ewygdfftdf}
\end{EQA}
where the sign \( \propto \) means equality up to a normalizing multiplicative constant.
Assume that the penalized maximum likelihood estimator (pMLE) \( \tilde{\thetav}_{\GP} \) is well defined:
\begin{EQA}
	\tilde{\thetav}_{\GP}
	&=&
	\argmax_{\thetav} \bigl\{ L(\thetav) - \| \GP(\thetav - \thetav_{0}) \|^{2}/2 \bigr\} .
\label{btgrtgthyhdwwwesdde}
\end{EQA}
Clearly \( \tilde{\thetav}_{\GP} \) is also maximizer of \( \priord_{\GP}(\thetav) \).
That is why it is often referred to as maximum a posteriori (MAP) estimator. 
Let also the log-likelihood function \( L(\thetav) \) be twice differentiable and weakly concave.
Define 
\begin{EQA}
	\IF_{\GP}(\thetav)
	&=&
	- \nabla^{2} L(\thetav) + \GP^{2} .
\label{ggtyjjfreedcvhjoiutghbn}
\end{EQA}
Assuming the latter expression to be positive definite for all considered \( \thetav \), 
define also its square root \( \DPGP(\thetav) = \sqrt{\IF_{\GP}(\thetav)} \).
We use the shortcut \( \DPGPt = \DPGP(\tilde{\thetav}_{\GP}) \).
Laplace's approximation means that the posterior distribution \( \priord_{\GP} \) is close 
to the Gaussian distribution \( \ND(\tilde{\thetav}_{\GP},\DPGPt^{-2}) \). 
A closely related Bernstein - von Mises phenomenon claims an approximation of the posterior
by \( \ND(\tilde{\thetav},\DP^{-2}) \), where \( \tilde{\thetav} \) is the MLE
and \( \DP^{2} = \IF = - \nabla^{2} L(\thetavs) \) is  
the Fisher information matrix for the true parameter value \( \thetavs \);
see e.g. \cite{SSW:2020} for a detailed discussion in context on nonlinear inverse problems.
The mentioned results provide an efficient tool for Bayesian uncertainty quantification and constructing the elliptic 
credible sets as level sets of the approximating Gaussian distribution; see \cite{HeKr2022} or
\cite{Reich2020} for applications to drift and diffusion estimation.

\Section{Motivation: Gradient free optimization}
Suppose that the point of maximum \( \xvs \) of the function \( \lgdL(\xv) \) is not known and has to be evaluated numerically.
We also assume that the function \( \lgdL(\cdot) \) is sufficiently smooth and can be efficiently computed at any point \( \xv \), 
however, its gradient is not available.
\cite{NeSp2016} offered a powerful gradient free method based on averaging the exponent \( \ex^{\lgdL(\xv)} \) w.r.t. to some Gaussian 
distribution for \( \xv \).
Here we aim at reconsidering this idea within the Bayesian optimization framework.
Namely, we intend to design a procedure based on Laplace iterations with an updated Gaussian prior at each step.
Let \( \priord_{0} \) be a starting Gaussian prior \( \priord_{0} \sim \ND(\xv_{0},\GP_{0}^{-2}) \).
The corresponding posterior is defined by normalizing the product \( \ex^{\lgdL(\xv)} \priord_{0}(\xv) \).
Due to the Laplace approximation result, this posterior is nearly Gaussian
with the mean \( \xv_{1} = \argmax \bigl\{ \lgdL(\xv) - \| \GP_{0} \xv \|^{2}/2 \bigr\} \); see e.g. Theorem~\ref{TLaplaceTV} below.
This leads to the general idea of Bayesian optimization: 
use the posterior mean/variance to update the prior and repeat the Laplace approximation step.
One may use standard Monte-Carlo or quasi Monte-Carlo methods for numerical approximation of the posterior mean;
see e.g. \cite{SSW:2020} in context of Bayesian inverse problems.

\Section{Challenges}
In spite of numerous existing results on Laplace approximation, some questions are still open.
A small list of relevant issues is given below.

\paragraph{Prior impact in high dimension}
Modern applications in statistics and optimization force to extend the scope of applicability 
of the classical results on Laplace approximation by including the cases 
of a very large dimension even exceeding the sample size.
It appears that in the case of a high-dimensional parameter, the prior is not washing out from the posterior.
In the contrary to the classical theory, it becomes extremely important.
It is of crucial importance to understand well the impact of the prior and its proper choice.

\paragraph{Smoothness conditions}
Standard smoothness conditions in terms of a uniform bound on the third derivative of \( f \) could be too restrictive and hard to check.

\paragraph{The use of posterior mean}
Another issue is possibility of using a Gaussian approximation of the posterior with inexact parameters.
Indeed, the pMLE/MAP \( \tilde{\thetav}_{\GP} \) is hard to compute, the classical MCMC based Bayesian computations
deliver an estimate of the posterior mean and of posterior covariance.
A big question is to justify the use of these quantities for uncertainty quantification or Bayesian optimization. 

\paragraph{Non-concavity}
The classical results of Laplace approximation are stated under the condition that the function \( \lgd(\cdot) \) is strongly 
globally concave.
In many applications including nonlinear inverse problems 
\cite{nickl2020, SSW:2020, HeKr2022} or Deep Neuronal Networks such an assumption appears to be unrealistic.



\Section{This paper's contributions}
This paper aims at reconsidering the classical results on Laplace approximation and to address the above issues.
Below the list of the most important achievements in the paper.
\paragraph{Effective dimension and dimension free guarantees}
We introduce the notion of \emph{effective dimension} \( \dimG \) of the problem which can be small of moderate even 
for huge parameter dimension \( \dimp \).
The value \( \dimG \) is defined by an interplay between the information delivered by the data 
and information contained in the prior; see Section~\ref{SeffdimLa} for more details. 
Further we establish
explicit \emph{non-asymptotic}  and \emph{dimension free} guarantees for the accuracy of Gaussian approximation 
of the posterior in \emph{total variation} (TV) distance in terms of effective dimension; see Theorem~\ref{TLaplaceTV}.
In the case when the non-penalized log-likelihood function grows linearly with the sample size \( n \), 
the quality of Laplace approximation is of order \( \sqrt{\dimG^{3}/n} \).
It can be improved to \( \dimG^{3}/n \) if instead of TV-distance, we limit ourselves to the class of centrally symmetric sets.
The proofs combine classical variational arguments with sharp bounds for Gaussian quadratic forms.
Conditions require that \( \lgd \) is strongly concave  
and locally smooth with a uniform bound on the third  
Gateaux derivative of \( \lgd \) in a local vicinity of the point of maximum. 

\paragraph{Critical dimension}
The result of Theorem~\ref{TLaplaceTV} helps to address the issue of \emph{critical dimension} for applicability
of Laplace approximation:
the relation \( \dimG^{3} \ll n \) between the sample size \( n \) and the effective dimension \( \dimG \) is sufficient 
for our main results.
The result on concentration of the posterior only requires \( \dimG \ll n \).

\paragraph{Posterior mean in place of MAP} 
The use of the posterior mean \( \xvb \) in place of the MAP \( \xvs \) is justified by Theorem~\ref{TpostmeanLan}
under the same condition \( \dimG^{3} \ll n \).
The obtained results justify a Bayesian optimization procedure based on iterated Laplace approximation
with a sequentially adjusted prior. 

\paragraph{Non-concavity}
Section~\ref{Slaplnonlin} indicates 
for the special case of nonlinear inverse problem
how the assumption of strong concavity of \( \lgd(\cdot) \) can be relaxed
using a ``warm start'' condition
which means that a good starting guess is available.

\medskip
The paper is organized as follows.
Section~\ref{STaylor} presents general bounds on the error of Laplace approximation and also 
discusses the use of posterior mean in place of posterior mode as well as
the Bayesian optimization procedure based on Laplace iterations.
Section~\ref{Slaplnonlin} specifies the results to the case on nonlinear inverse problem.
Proofs and technical tools are collected in the Appendix.

\def\lgd{f}
\def\elll{\ell}
\def\PfL{\P_{\lgd}}
\def\lgdL{\elll}
\def\lgdLy{\lgd_{\elll,\nuiv}}
\def\hL{h}
\def\biasL{\bias}

\def\DU{\mathbb{D}}
\def\DV{\mathsf{D}}
\def\DVL{\DV}
\def\DVn{\DV_{\elll}}
\def\DVLy{\DV_{\elll,\nuiv}}
\def\DVb{\breve{\DV}}
\def\DVLb{\DVb}
\def\DVbin{\DVb_{0}}
\def\DVLbin{\DVbin}
\def\DVF{\cc{D}}
\def\IFV{\mathsf{F}}
\def\IFVL{\IFV}

\def\nui{s}
\def\Nui{S}
\def\nuiv{\bb{\nui}}
\def\Nuiv{\bb{\Nui}}

\def\dimLG{\dimA_{\GP}}
\def\dimLy{\dimA_{\nuiv}}
\def\dimLs{\dimLG^{*}}
\def\dimx{\dimA_{\XX_{0}}}

\def\rrin{\rr_{0}}
\def\rrL{\rr_{\GP}}
\def\rrLb{\rr}
\def\rrLy{\rr_{\nuiv}}

\def\UV{\mathcal{U}}
\def\UVL{\UV_{\elll}}
\def\UVL{\UV}
\def\UVLy{\UV_{\elll,\nuiv}}
\def\UVb{\breve{\UV}}
\def\UVLb{\UVL}

\def\rhos{\rho^{*}}
\def\hmax{\mathsf{c}}
\def\Rem{\mathcal{R}}
\def\Reme{\mathcal{R}}

\def\HU{\mathsf{H}}
\def\feta{\phi}
\def\errf{\err_{\feta}}
\def\errs{\err^{*}}
\def\lamb{\bb{\lambda}}

\def\regrf{m}
\def\regrfv{\bb{\regrf}}
\def\XX{\mathcal{X}}
\def\CDG{\CONSTi_{\GP}}

\Chapter{Dimension free bounds for Laplace approximation}
\label{STaylor}
Here we present some general results on accuracy of Laplace approximation.

\Section{Setup and conditions}
Let \( \lgd(\xv) \) be a function in a high-dimensional Euclidean space \( \R^{\dimp} \) such that
\( \int \ex^{\lgd(\xv)} \, d\xv = \CONST < \infty \),
where the integral sign \( \int \) without limits means the integral over the whole space \( \R^{\dimp} \).
Then \( \lgd \) determines a distribution \( \PfL \) with the density
\( \CONST^{-1} \ex^{\lgd(\xv)} \).
Let \( \xvs \) be a point of maximum:
\begin{EQA}
	\lgd(\xvs)
	&=&
	\sup_{\uv \in \R^{\dimp}} \lgd(\xvs + \uv) .
\label{scdygw7ytd7wqqsquuqydtdtd}
\end{EQA}
We also assume that \( \lgd(\cdot) \) is smooth, more precisely, three or even four time differentiable. 
Introduce the negative Hessian \( \DU^{2} = - \nabla^{2} \lgd(\xvs) \) and assume \( \DU^{2} \) strictly positive definite.
We aim at approximating the measure \( \PfL \) by a Gaussian measure \( \ND(\xvs,\DU^{-2}) \).
Given a function \( g(\cdot) \), define its expectation w.r.t. \( \PfL \) after centering at \( \xvs \):
\begin{EQA}
	\II(g)
	& \eqdef &
	\frac{\int g(\uv) \, \ex^{\lgd(\xvs + \uv)} \, d\uv}{\int \ex^{\lgd(\xvs + \uv)} \, d\uv} \, .
\label{IIgfigefxududu}
\end{EQA}
A Gaussian approximation \( \II_{\DU}(g) \) for \( \II(g) \) reads:
\begin{EQA}
	\II_{\DU}(g)
	& \eqdef &
	\frac{\int g(\uv) \, \ex^{- \| \DU \uv \|^{2}/2} \, d\uv}{\int \ex^{- \| \DU \uv \|^{2}/2} \, d\uv} 
	=
	\E g(\gaussv_{\DU}) ,
	\qquad
	\gaussv_{\DU} \sim \ND(0,\DU^{-2}) \, .
\label{IgiLapgeHvdu}
\end{EQA}
The choice of the distance between \( \PfL \) and \( \ND(\xvs,\DU^{-2}) \) depends on the considered class 
of functions \( g \).
The most strong total variation distance can be obtained as 
the supremum of \( |\II(g) - \II_{\DU}(g)| \) over all measurable functions \( g(\cdot) \) with 
\( |g(\uv)| \leq 1 \):
\begin{EQA}
	\TV\bigl( \PfL,\ND(\xvs,\DU^{-2}) \bigr)
	&=&
	\sup_{\|g\|_{\infty} \leq 1} \bigl| \II(g) - \II_{\DU}(g) \bigr| \, .
\label{ghdvcgftdftgegdftefd3434545}
\end{EQA}
The results can be substantially improved if only centrally symmetric functions \( g(\cdot) \)
with \( g(\xv) = g(-\xv) \) are considered.
Obviously 
\begin{EQA}
	\II(g)
	& = &
	\frac{\int g(\uv) \, \ex^{\lgd(\xvs + \uv) - \lgd(\xvs)} \, d\uv}{\int \ex^{\lgd(\xvs + \uv) - \lgd(\xvs)} \, d\uv} \, .
\label{IIgfigefxuudufxudu}
\end{EQA}
Moreover, as \( \xvs = \argmax_{\xv} \lgd(\xv) \), it holds \( \nabla \lgd(\xvs) = 0 \) and 
\begin{EQA}
	\II(g)
	& = &
	\frac{\int g(\uv) \, \ex^{\lgd(\xvs;\uv)} \, d\uv}{\int \ex^{\lgd(\xvs;\uv)} \, d\uv} \, ,
\label{IIgfifxutudufxutdu}
\end{EQA}
where \( \lgd(\xv;\uv) \) is the Bregman divergence 
\begin{EQA}
	\lgd(\xv;\uv)
	&=&
	\lgd(\xv + \uv) - \lgd(\xv) - \bigl\langle \nabla \lgd(\xv), \uv \bigr\rangle .
\label{fxufxpufxfpxu}
\end{EQA}
Implicitly we assume that the negative Hessian \( \DU^{2} = - \nabla^{2} \lgd(\xvs) \) is sufficiently large
in the sense that the Gaussian measure \( \ND(0,\DU^{-2}) \) concentrates on a local set \( \UVL \).
This allows to use a local Taylor expansion for 
\( \lgd(\xvs;\uv) \approx - \| \DU \uv \|^{2}/2 \) in \( \uv \) on \( \UVL \).
If \( \lgd(\cdot) \) is also strongly concave, then the mass of \( \PfL \)
of the complement of \( \UVL \) is exponentially small yielding the desirable Laplace approximation.

Motivated by applications to statistical inference, we consider \( \lgd \) in a special form
\begin{EQA}
	\lgd(\xv)
	&=&
	\lgdL(\xv) - \| \GP (\xv - \xv_{0}) \|^{2} / 2
\label{jhdctrdfred4322edt7y}
\end{EQA}
for a symmetric \( \dimp \)-matrix \( \GP^{2} \geq 0 \).
Here \( \lgdL(\cdot) \) stands for the log-likelihood function while the quadratic penalty 
\( \| \GP (\xv - \xv_{0}) \|^{2} / 2 \) corresponds to a Gaussian prior \( \ND(\xv_{0},\GP^{-2}) \).
We also assume that \( \lgdL(\cdot) \) is concave with \( \DV^{2} \eqdef - \nabla^{2} \lgdL(\xv) \geq 0 \).
Then clearly 
\begin{EQA}
	- \nabla^{2} \lgd(\xv)
	&=&
	- \nabla^{2} \lgdL(\xv) + \GP^{2} 
	=
	\DVL^{2} + \GP^{2} .
\label{hydsf42wdsdtrdstdg}
\end{EQA} 
In typical asymptotic setups, the log-likelihood function \( \lgdL(\xv) \) grows with the sample size 
or inverse noise variance, while the prior precision matrix \( \GP^{2} \) is kept fixed.
Decomposition \eqref{jhdctrdfred4322edt7y} is of great importance for obtaining the dimension free results.
The main reason is that the quadratic penalty does not affect smoothness properties of the function \( \lgdL(\cdot) \)
but greatly improves the quadratic approximation term.
Now we state precise conditions.

\Subsection{Concavity}
Below we implicitly assume decomposition \eqref{jhdctrdfred4322edt7y} with 
a weakly concave function \( \lgdL(\cdot) \).
More specifically, we assume the following condition.

\begin{description}
    \item[\label{LLf0ref} \( \bb{(\cc{C}_{0})} \)]
      \textit{ There exists an operator \( \GP^{2} \leq - \nabla^{2} \lgd(\xvs) \) in \( \R^{\dimp} \) such that the function 
\begin{EQA}
	\lgdL(\xvs + \uv)
	& \eqdef &
	\lgd(\xvs + \uv) + \frac{1}{2} \| \GP \uv \|^{2}
\label{fTvufupumDTpDfx}
\end{EQA}
is concave. 
      }
\end{description}

If \( \lgdL(\cdot) \) in decomposition \eqref{jhdctrdfred4322edt7y} is concave then this condition is obviously fulfilled. 
More generally, if \( \lgdL(\cdot) \) in \eqref{jhdctrdfred4322edt7y} is weakly concave, so that 
\( \lgdL(\xvs + \uv) - \| \GP_{0} \uv \|^{2}/2 \) is concave in \( \uv \) with \( \GP_{0}^{2} \leq \GP^{2} \), then 
\nameref{LLf0ref} is fulfilled with \( \GP^{2} - \GP_{0}^{2} \) in place of \( \GP^{2} \).

The operator 
\begin{EQA}
	\DV^{2} 
	&=&
	- \nabla^{2} \lgd(\xvs) - \GP^{2}
	\quad
	\bigl( = - \nabla^{2} \lgdL(\xvs) \text{ under \eqref{jhdctrdfred4322edt7y}} \bigr)
\label{kc7c322dgdf43djhd}
\end{EQA}
plays an important role in our conditions and results.

\begin{remark}
The condition of strong concavity of \( \lgd \) on the whole space \( \R^{\dimp} \)
can be too restrictive; see the example of nonlinear inverse problem in Section~\ref{Slaplnonlin}.
This condition can be replaced by its local version: there exists a set \( \XX_{0} \) 
such that the Gaussian prior \( \ND(\xv_{0},\GP^{-2}) \) concentrates on \( \XX_{0} \)
with a high probability and the maximizer \( \xvs_{\GP} \) belongs to \( \XX_{0} \).
In all the results, the integral over \( \R^{\dimp} \) has to replaced by the integral over \( \XX_{0} \).
\end{remark}

\Subsection{Effective dimension}
\label{SeffdimLa}
With decomposition \eqref{kc7c322dgdf43djhd} in mind, we use another notation for \( \DU^{2} = - \nabla^{2} \lgd(\xvs) \):
\begin{EQA}
	\DV_{\GP}^{2} 
	&=& 
	- \nabla^{2} \lgd(\xvs)
	=
	\DVL^{2} + \GP^{2} .
\label{dscytf5w2edte5rw4e24gyd}
\end{EQA}
Also we write \( \II_{\GP}(g) \) instead of \( \II_{\DU}(g) \) in \eqref{IgiLapgeHvdu}.
The \emph{effective dimension} \( \dimG \) is given by
\begin{EQA}
	\dimG
	& \eqdef & 
	\tr\bigl( \DVL^{2} \, \DV_{\GP}^{-2} \bigr) .
\label{dAdetrH02Hm2}
\end{EQA}
Of course, \( \dimG \leq \dimp \) but a proper choice of the penalty \( \GP^{2} \) in \eqref{jhdctrdfred4322edt7y}
allows to avoid the ``curse of dimensionality'' issue and ensure a small or moderate effective dimension \( \dimG \) even for 
\( \dimp \) large or infinite. 
The value \( \dimG \) helps to describe a local vicinity \( \UVL \) around \( \xvs \) such that the most of mass of \( \PfL \)
concentrates on \( \UVL \); see Section~\ref{StailLa}.
Namely, let us fix some \( \amax < 1 \), e.g. \( \amax = 2/3 \), and some \( \xx > 0 \) ensuring that \( \ex^{-\xx} \)
is our significance level.
Define
\begin{EQ}[rcl]
	\rrL
	& = &
	2 \sqrt{\dimG} + \sqrt{2 \xx} ,
	\\
	\UVL 
	&=& 
	\bigl\{ \uv \colon \| \DVL \uv \| \leq \amax^{-1} \rrL \bigr\} .
\label{UvTDunm12spT}
\end{EQ}

\Subsection{Local smoothness conditions}
Let \( \dimp \leq \infty \) and
let \( \lgd(\cdot) \) be a three times continuously differentiable function on \( \R^{\dimp} \).
We fix a reference point \( \xv \) and local region around \( \xv \) given by the local set 
\( \UVL \subset \R^{\dimp} \) from \eqref{UvTDunm12spT}. 
Also consider the second order Taylor approximation 
\( \lgd(\xv + \uv) \approx 
\lgd(\xv) + \bigl\langle \nabla \lgd(\xv), \uv \bigr\rangle + \frac{1}{2} \bigl\langle \nabla^{2} \lgd(\xv) , \uv^{\otimes 2} \bigr\rangle \)
and similarly the third order expansion and introduce the remainders 
\begin{EQ}[rcl]
	\dltw_{3}(\xv,\uv)
	&=&
	\lgd(\xv;\uv) - \frac{1}{2} \bigl\langle \nabla^{2} \lgd(\xv) , \uv^{\otimes 2} \bigr\rangle ,
	\\
	\dltw_{4}(\xv,\uv)
	&=&
	\lgd(\xv;\uv) - \frac{1}{2} \bigl\langle \nabla^{2} \lgd(\xv) , \uv^{\otimes 2} \bigr\rangle 
	- \frac{1}{6} \bigl\langle \nabla^{3} \lgd(\xv), \uv^{\otimes 3} \bigr\rangle 
\label{d4fuv1216303}
\end{EQ}
with \( \lgd(\xv;\uv) \) from \eqref{fxufxpufxfpxu}.
The use of the Taylor formula allows to bound
\begin{EQ}[rcl]
	\bigl| \dltw_{k}(\xv,\uv) \bigr|
	& \leq &
	\sup_{t \in [0,1]}
	\frac{1}{k!} 
	\Bigl| \bigl\langle \nabla^{k} \lgd(\xv + t \uv), \uv^{\otimes k} \bigr\rangle \Bigr|,
	\quad
	k\geq 3. 
\label{k3t01f12n3fvtu}
\end{EQ}
It is worth noting that the quadratic penalty \( - \| \GP (\xv - \xv_{0}) \|^{2}/2 \) in \( \lgd \) does not affect 
the remainders \( \dltw_{3}(\xv,\uv) \) and \( \dltw_{4}(\xv,\uv) \).
Indeed, 
with \( \lgd(\xv) = \lgdL(\xv) - \| \GP (\xv - \xv_{0}) \|^{2}/2 \), it holds 
\begin{EQA}
	\lgd(\xv;\uv) 
	& \eqdef & 
	\lgd(\xv + \uv) - \lgd(\xv) - \bigl\langle \nabla \lgd(\xv), \uv \bigr\rangle 
	=
	\lgdL(\xv;\uv) - \| \GP \uv \|^{2}/2 
\label{D2fuvGu2fG}
\end{EQA}
and the quadratic term in definition of the values \( \dltw_{k}(\xv,\uv) \) cancels, \( k \geq 3 \). 
Local smoothness of \( \lgd(\cdot) \) or, equivalently, of \( \lgdL(\cdot) \), at \( \xv \) will be measured by
the value \( \dltwb(\xv) \):
\begin{EQA}
	\dltwb(\xv)
	& \eqdef &
	\sup_{\uv \in \UVL} \frac{1}{\| \DVL \uv \|^{2}/2} \bigl| \dltw_{3}(\xv,\uv) \bigr| .
\label{om3esuU1H02d3}
\end{EQA}
We also denote
\begin{EQA}
	\dltwb
	& \eqdef &
	\dltwb(\xvs) .
\label{dsviu4dskjbdsiqwfcerdq2}
\end{EQA}
Our results apply under the condition \( \dltwb \ll 1 \).
Local concentration of the measure \( \PfL \) requires \( \dltwb \leq 1/3 \);
see Proposition~\ref{PlocconLa}.
The main results about Gaussian approximation of \( \PfL \) are valid under a stronger condition 
\( \dltwb \, \dimLG \leq 2/3 \) with the effective dimension \( \dimLG \) from \eqref{dAdetrH02Hm2}.


\Section{Error bounds for Laplace approximation}
Our first result describes the quality of approximation of the measure \( \PfL \) by the Gaussian measure \( \ND(\xvs,\DV_{\GP}^{-2}) \)
with mean \( \xv \) and the covariance \( \DV_{\GP}^{-2} \) in total variation distance.
In all our result, the value \( \xx \) is fixed to ensure that \( \ex^{-\xx} \) is negligible.
First we present the general results which will be specified later under the self-concordance condition.

\Subsection{General bounds}

\begin{theorem}
\label{TLaplaceTV}
Suppose \nameref{LLf0ref}.
Let also \( \dimLG \) be defined by \eqref{dAdetrH02Hm2} and \( \rrL \) and \( \UVL \) by \eqref{UvTDunm12spT}.
If \( \dltwb \) from \eqref{om3esuU1H02d3} satisfies \( \dltwb \leq 1/3 \), then
\begin{EQA}
	\PfL(\Xv - \xvs \not\in \UVL)
	& \leq &
	\ex^{-\xx} .
\label{poybf3679jd532ff2}
\end{EQA}
If \( \dltwb \, \dimLG \leq 2/3 \),
then for any \( g(\cdot) \) with \( |g(\uv)| \leq 1 \), it holds for \( \II(g) \) from \eqref{IIgfigefxududu}
\begin{EQA}
	\bigl| \II(g) - \II_{\GP}(g) \bigr|
	& \leq &
	\frac{2 (\err + \ex^{-\xx})}{1 - \err - \ex^{-\xx}}
	\leq 
	4(\err + \ex^{-\xx}) 
\label{ufgdt6df5dtgededsxd23gjg}
\end{EQA}
with 
\begin{EQA}
	\err
	&=&
	\err_{2} 
	= 
	\frac{0.75 \, \dltwb \, \dimLG}{1 - \dltwb} \, .
\label{juytr90f2dzaryjhfyf}
\end{EQA}
\end{theorem}

This section presents more advanced bounds on the error of Laplace approximatio.
Introduce the following conditions.
\begin{description}
    \item[\label{LL3fref} \( \bb{(\LL_{3})} \)]
      \textit{ There exists a 3-homogeneous function \( \dltwu_{3}(\uv) \), 
\( \dltwu_{3}(t \uv) = |t|^{3} \dltwu_{3}(\uv) \), such that }
\begin{EQA}
	\bigl| \dltw_{3}(\xvs,\uv) \bigr|
	& \leq &
	\frac{1}{6} \dltwu_{3}(\uv) \, .
\label{bd3xu16f3uo3}
\end{EQA}
\end{description}
 
\begin{description}
    \item[\label{LL4fref} \( \bb{(\LL_{4})} \)]
      \textit{ There exists a 4-homogeneous function \( \dltwu_{4}(\uv) \), 
\( \dltwu_{4}(t \uv) = |t|^{4} \dltwu_{3}(\uv) \), such that }
\begin{EQA}
	\bigl| \dltw_{4}(\xvs,\uv) \bigr|
	& \leq &
	\frac{1}{24} \dltwu_{4}(\uv) \, .
\label{1mffmxum5}
\end{EQA}
\end{description}

\begin{theorem}
\label{TLaplaceTV2}
Suppose \nameref{LLf0ref} and \nameref{LL3fref} and let \( \dltwb \, \dimLG \leq 2/3 \).
Then the accuracy bound \eqref{ufgdt6df5dtgededsxd23gjg} applies with
\begin{EQA}
	\err
	&=&
	\err_{3}
	\eqdef
	\frac{\E \dltwu_{3}(\gaussv_{\GP})}{4 (1 - \dltwb)^{3/2}} \, .
\label{5qw7dyf4e4354coefw9dufih}
\end{EQA}
Moreover, under \nameref{LL4fref}, the accuracy bound \eqref{ufgdt6df5dtgededsxd23gjg} applies to any 
symmetric function \( g(\uv) = g(-\uv) \), \( |g(\uv)| \leq 1 \), with 
\begin{EQA}
	\err
	&=&
	\err_{4}
	=
	\frac{1}{16 (1 - \dltwb)^{2}} \Bigl\{ \E \bigl\langle \nabla^{3} \lgd(\xvs) , \gaussv_{\GP}^{\otimes 3} \bigr\rangle^{2} 
	+ 2 \E \dltwu_{4}(\gaussv_{\GP}) \Bigr\} \, .
\label{hg25t6mxwhydseg3hhfdr}
\end{EQA}
\end{theorem}

Let \( \BBB(\R^{\dimp}) \) be the \( \sigma \)-field of all Borel sets in \( \R^{\dimp} \),
while \( \BBB_{s}(\R^{\dimp}) \) stands for all centrally symmetric sets from \( \BBB(\R^{\dimp}) \).
By \( \Xv \) we denote a random element with the distribution \( \PfL \), while \( \gaussv_{\GP} \sim \ND(0,\DV_{\GP}^{-2}) \).

\begin{corollary}
\label{CTLaplaceTV}
Under the conditions of Theorem~\ref{TLaplaceTV2}, it holds for \( \Xv \sim \PfL \)
\begin{EQA}
	\sup_{A \in \BBB(\R^{\dimp})} \bigl| \PfL(\Xv - \xvs \in A) - \P(\gaussv_{\GP} \in A) \bigr|
	& \leq &
	4(\err_{3} + \ex^{-\xx}),
	\\
	\sup_{A \in \BBB_{s}(\R^{\dimp})} \bigl| \PfL(\Xv - \xvs \in A) - \P(\gaussv_{\GP} \in A) \bigr|
	& \leq &
	4(\err_{4} + \ex^{-\xx}) .
\label{cbc5dfedrdwewwgerg}
\end{EQA}
\end{corollary}

\Subsection{Bounds under self-concordance}
Often the smoothness properties of \( \lgd \) are described in terms of the third derivative.
Our approach does not require a bounded third or fourth full dimensional derivative, instead we consider 
directional Gateaux derivatives.
Let \( \lgdL(\cdot) \) be three times continuously differentiable and \( \nabla^{3} \lgdL(\cdot) \) 
stand for the third derivative.
Then the Taylor expansion of the third order implies \( \dltwb \leq \dltwa/3 \) with
\begin{EQA}
	\dltwa
	&=&
	\sup_{\uv \in \UVL} \, \sup_{t \in [0,1]} 
	\frac{\bigl| \langle \nabla^{3} \lgdL(\xvs + t \uv), \uv^{\otimes 3} \rangle \bigr|}{\| \DVL \uv \|^{2}} \, .
\label{fc4gegfet5df4rt3wg7ygrrftg}
\end{EQA} 
Moreover, in many applications, the function \( \lgdL(\cdot) \) is of the form \( - \lgdL(\xv) = n \hL(\xv) \) 
for a fixed smooth function \( h(\cdot) \) satisfying the following condition:
\begin{description}
    \item[\label{LLcS3ref} \( \bb{(\cc{S}_{3})} \)]
      \emph{ \( - \lgdL(\cdot) = n \hL(\cdot) \) for a strongly convex function \( h(\cdot) \) and 
\begin{EQA}
	\sup_{\uv \in \UVL} \, \sup_{t \in [0,1]} 
	\frac{\bigl| \langle \nabla^{3} \hL(\xvs + t \uv), \uv^{\otimes 3} \rangle \bigr|}{\langle \nabla^{2} \hL(\xvs), \uv^{\otimes 2} \rangle^{3/2}}
	& \leq &
	\hmax_{3} \, .
\end{EQA}
}
    \item[\label{LLcS4ref} \( \bb{(\cc{S}_{4})} \)]
      \emph{ \( - \lgdL(\cdot) = n \hL(\cdot) \) for a strongly convex function \( h(\cdot) \) and 
\begin{EQA}
	\sup_{\uv \in \UVL} \, \sup_{t \in [0,1]} 
	\frac{\bigl| \langle \nabla^{4} \hL(\xvs + t \uv), \uv^{\otimes 4} \rangle \bigr|}{\langle \nabla^{2} \hL(\xvs), \uv^{\otimes 2} \rangle^{2}}
	& \leq &
	\hmax_{4} \, .
\end{EQA}
}
\end{description}

This is a local version of the so called self-concordance condition; see \cite{Ne1988}.
Under this condition, we can easily bound the values \( \dltw_{k}(\xvs,\uv) \) for \( k=3,4 \),
\( \dltwb \) and \( \dltwa \):
\begin{EQA}
	\dltw_{k}(\xvs,\uv)
	\leq 
	\frac{\hmax_{k} \, n^{-1/2}}{k!} \| \DVL \uv \|^{k} ;
	\qquad
	\dltwb
	& \leq &
	\dltwa/3
	\leq 
	\hmax_{3} \, \rrL \, n^{-1/2} / 3 \, ;
	\qquad
\label{gtcdsftdfvtwdsefhfdvfrvsewse3}
\end{EQA}
see Lemma~\ref{LdltwLa} later.

Finally we state the bounds under \nameref{LLcS3ref} and \nameref{LLcS4ref}.
\begin{theorem}
\label{TLaplaceTV34}
Suppose \nameref{LLf0ref}, \nameref{LLcS3ref}, and let \( \hmax_{3} \, \rrL \, n^{-1/2} \leq 3/4 \)
for \( \rrL \) from \eqref{UvTDunm12spT}.
Then 
\begin{EQA}
	\sup_{A \in \BBB(\R^{\dimp})} \bigl| \PfL(\Xv - \xvs \in A) - \P(\gaussv_{\GP} \in A) \bigr|
	& \leq &
	2 \hmax_{3} \, \sqrt{\frac{(\dimLG+1)^{3}}{n}} + 4 \ex^{-\xx} \, .
\label{scdugfdwyd2wywy26e6de}
\end{EQA}
If \nameref{LLcS4ref} is also satisfied then
\begin{EQA}
	\sup_{A \in \BBB_{s}(\R^{\dimp})} \bigl| \PfL(\Xv - \xvs \in A) - \P(\gaussv_{\GP} \in A) \bigr|
	& \leq &	
	\frac{\hmax_{3}^{2} \, (\dimG + 2)^{3} + 2 \hmax_{4} (\dimLG + 1)^{2}}{2n} + 4 \ex^{-\xx} \, .
	\qquad
\label{scdugfdwyd2wywy26e6de4}
\end{EQA}
\end{theorem}

\Subsection{Critical dimension}
Here we briefly discuss the important issue of \emph{critical dimension}.
Theorem~\ref{TLaplaceTV} states concentration of \( \PfL \) under the condition
\( \dltwb \leq 1/3 \) and Gaussian approximation under the stronger condition
\( \dltwb \, \dimLG \leq 2/3 \).
Moreover, self-concordance \nameref{LLcS3ref} and Lemma~\ref{LdltwLa} enables us to bound 
\( \dltwb \lesssim \sqrt{\dimLG/n} \)
and \( \dltwb \, \dimLG \lesssim \sqrt{\dimLG^{3}/n} \).
Hence, concentration of \( \PfL \) requires the \emph{critical dimension} condition 
\( \dimLG \ll n \), while Gaussian approximation applies under \( \dimLG^{3} \ll n \).
We see that there is a gap between these conditions. 
We guess that in the region \( n^{1/3} \lesssim \dimLG \lesssim n \), 
a non-Gaussian approximation of the posterior is possible. 
\cite{bochkina2014} provides examples of non-Gaussian posterior limits
for non-regular models.

\Subsection{Kullback-Leibler divergence}
Theorem~\ref{TLaplaceTV} through \ref{TLaplaceTV34} quantify the approximation \( \PfL \approx \ND(\xvs,\DV_{\GP}^{-2}) \)
in the total variation distance.
Another useful characteristic could be the Kullback-Leibler (KL) divergence between \( \PfL \) and  
\( \ND(\xvs,\DV_{\GP}^{-2}) \).
The KL divergence \( \kullb(\P_{1},\P_{2}) = \E_{1} \log(d\P_{1}/d\P_{2}) \) is asymmetric, 
\( \kullb(\P_{1},\P_{2}) \neq \kullb(\P_{2},\P_{1}) \) with few exceptions 
like the case of Gaussian measures \( \P_{1} \) and \( \P_{2} \).
Moreover, \( \kullb(\P_{1},\P_{2}) \) can explode if \( \P_{2} \) is not absolutely continuous w.r.t. \( \P_{1} \). 
We present two bounds for each ordering.
For ease of presentation, we limit ourselves to the case when either \nameref{LL3fref} or \nameref{LLcS3ref} meets.

\begin{theorem}
\label{TLaplaceKL}
Suppose \nameref{LLf0ref} and \nameref{LL3fref} and let \( \dltwb \, \dimLG \leq 2/3 \).
With \( \P_{\GP} = \ND(\xvs,\DV_{\GP}^{-2}) \),
\begin{EQA}
	\kullb(\PfL,\P_{\GP})
	& \leq &
	4 \err_{3} + 4 \ex^{-\xx} 
	\leq 
	\frac{\E \dltwu_{3}(\gaussv_{\GP})}{(1 - \dltwb)^{3/2}} + 4 \ex^{-\xx} .
\label{vlgvi8ugu7tr4ry43et31}
\end{EQA}
Moreover, under \nameref{LLcS3ref}
\begin{EQA}
	\kullb(\PfL,\P_{\GP})
	& \leq &
	2 \hmax_{3} \, \sqrt{\frac{(\dimLG+1)^{3}}{n}} + 4 \ex^{-\xx} .
\label{cheiufheurhrtdgwb3wesrtd}
\end{EQA}
\end{theorem}

Now we briefly discuss the value \( \kullb(\P_{\GP},\PfL) \).
We already know that \( \PfL \) concentrates on \( \UVL \) and can be well approximated by \( \P_{\GP} \) on \( \UVL \).
However, this does not guarantee a small value of \( \kullb(\P_{\GP},\PfL) \).
It can even explode if e.g. \( \PfL \) has a compact support.
In fact, 
the log-density of \( \P_{\GP} \) w.r.t. \( \PfL \) reads
\begin{EQA}
	\log \frac{d\P_{\GP}}{d\PfL}(\xv)
	&=&
	- \lgd(\xv) - \frac{1}{2} \| \DV_{\GP} (\xv - \xv_{0}) \|^{2} - \CDG
\label{fduhefuifhew2tdw6ww3}
\end{EQA}
for some constant \( \CDG \),
and an upper bound on \( \kullb(\P_{\GP},\PfL) \) requires that the integral of \( \lgd(\xv) \) w.r.t. the measure \( \P_{\GP} \) is finite.

\begin{theorem}
\label{TLaplaceKLi}
Suppose \nameref{LLf0ref} and \nameref{LLcS3ref} and let \( \dltwb \, \dimLG \leq 2/3 \).
Let also \( \rho_{\GP} = 2\xx /\rrL^{2} \); see \eqref{UvTDunm12spT}.
If \( \lgdL(\xvs;\uv) = \lgdL(\xvs + \uv) - \lgdL(\xvs) - \langle \nabla \lgdL(\xvs), \uv \rangle \) fulfills
\begin{EQA}
	\int \bigl| \lgdL(\xvs;\uv) \bigr| \, \exp \bigl\{ - \| \DV_{\GP} \uv \|^{2}/2 + \rho_{\GP} \| \DVL \uv \|^{2}/2 \bigr\} \, d\uv
	& \leq &
	\CONSTi_{\lgdL} \, 
\label{dchbhwdhwwdgscsn2efty2162}
\end{EQA} 
for some fixed constant \( \CONSTi_{\lgdL} \) then
\begin{EQA}
	\kullb(\P_{\GP},\PfL)
	& \leq &
	\hmax_{3} \, \sqrt{\frac{(\dimLG+1)^{3}}{n}} + (2 + \CONSTi_{\lgdL}) \ex^{-\xx} .
\label{cheiufheurhrtdgwb3wesrtd}
\end{EQA}
\end{theorem}

\Subsection{Mean and MAP}
Here we present the bound on \( |\II(g) - \II_{\GP}(g)| \) for the case of a linear vector function \( g(\uv) = \QP \uv \) 
with \( \QP \colon \R^{\dimp} \to \R^{\dimq} \), \( \dimq \geq 1 \). 
A special case of \( \QP = \Id_{\dimp} \) corresponds to the mean value \( \xvb \) of \( \PfL \).
The next result presents an upper bound for the value \( \QP (\xvb - \xvs) \) under the conditions of Theorem~\ref{TLaplaceTV34}.

\begin{theorem}
\label{TpostmeanLa}
Assume the conditions of Theorem~\ref{TLaplaceTV34} and let \( \QP^{\T} \QP \leq \DVL^{2} \).
Then it holds with some absolute constant \( \CONST \)
\begin{EQA}
	\| \QP (\xvb - \xvs) \|
	& \leq &
	2.4 \, \hmax_{3} \, \| \QP \DV_{\GP}^{-2} \QP^{\T} \|^{1/2} \,
	\frac{(\dimLG + 1)^{3/2}}{n^{1/2}} + \CONST \ex^{-\xx} .
\label{hcdtrdtdehfdewdrfrhgyjufger}
\end{EQA}
\end{theorem}

Now we specify the result for the special choice \( \QP = \DVL \).

\begin{corollary}
\label{CTpostmeanLa}
Assume the conditions of Theorem~\ref{TpostmeanLa} and \nameref{LLcS3ref}.
Then 
\begin{EQA}
	\| \DVL (\xvb - \xvs) \|
	& \leq &
	2.4 \, \hmax_{3} \, \frac{(\dimLG + 1)^{3/2}}{n^{1/2}} + \CONST \ex^{-\xx} \, .
\label{klu8gitfdgregfkhj7yt}
\end{EQA}
\end{corollary}

\begin{remark}
\label{RTpostmeanc}
An interesting question is whether the result of Theorem~\ref{TpostmeanLa} or Corollary~\ref{CTpostmeanLa} applies 
with \( \QP = \DV_{\GP} \).
This issue is important in connection to inexact Laplace approximation; see the next section.
The answer is negative. 
The problem is related to the last term \( \CONST \ex^{-\xx} \) in the right hand-side of \eqref{hcdtrdtdehfdewdrfrhgyjufger}.
The constant \( \CONST \) involves the moments of \( \| \QP (\Xv - \xvs) \|^{2} \)
which explode for \( \QP = \DV_{\GP} \) and \( \dimp = \infty \).
\end{remark}

\Section{Inexact approximation and the use of posterior mean}
\label{SLaplinexact}

Now we change the setup.
Namely, we suppose that the true maximizer \( \xvs \) of the function \( \lgd \) is not available, but 
\( \xv \) is somehow close to the point of maximum \( \xvs \).
Similarly, the negative Hessian \( \DV_{\GP}^{2}(\xvs) = - \nabla^{2} \lgd(\xvs) \) is hard to obtain 
and we use a proxy \( \HU^{2} \).
%
We already know that \( \PfL \) can be well approximated by \( \ND(\xvs,\DV^{-2}) \) with 
\( \DV^{2} = \DV^{2}(\xvs) \).
This section addresses the question whether \( \ND(\xv,\HU^{-2}) \) can be used instead.
Here we may greatly benefit from the fact that Theorem~\ref{TLaplaceTV} provides a bound on the total variation distance 
between \( \PfL \) and \( \ND(\xvs,\DV_{\GP}^{-2}(\xvs)) \) yielding
\begin{EQA}
	\TV\bigl( \PfL, \ND(\xv,\HU^{-2}) \bigr)
	& \leq &
	\TV\bigl( \PfL, \ND(\xvs,\DV_{\GP}^{-2}) \bigr)
	+ \TV\bigl( \ND(\xv,\HU^{-2}), \ND(\xvs,\DV_{\GP}^{-2}) \bigr) .
	\qquad
\label{iuhvfuggeyfgeyjhefhgdy}
\end{EQA}
Therefore, it suffices to bound the TV-distance between two Gaussian distributions: 
\( \ND(\xvs,\DV_{\GP}^{-2}) \) naturally appears in the Laplace approximation, 
the second one is a proxy used instead.
Pinsker's inequality provides an upper bound: for any two measures \( P,Q \)
\begin{EQA}
	\TV(P,Q)
	& \leq &
	\sqrt{\kullb(P,Q)/2} ,
\label{gffwdygweufeuyfeygfdwyd}
\end{EQA}
where \( \kullb(P,Q) \) is the Kullback-Leibler divergence between \( P \) and \( Q \).
The KL-divergence between two Gaussians has a closed form:
\begin{EQA}
	\kullb\bigl( \ND(\xvs,\DV_{\GP}^{-2}), \ND(\xv,\HU^{-2}) \bigr)
	&=&
	\frac{1}{2} \bigl\{ \| \DV_{\GP} (\xv - \xvs) \|^{2} + \tr (\HU^{-2} \DV_{\GP}^{2} - \Id_{\dimp}) + \log \det (\HU^{-2} \DV_{\GP}^{2}) \bigr\} .
\label{h0bede4bweve7yjdyeghy}
\end{EQA}
Moreover, if the matrix \( \BB = \HU^{-1} \DV_{\GP}^{2} \HU^{-1} - \Id_{\dimp} \) satisfies \( \| \BB_{\GP} \| \leq 2/3 \) then
\begin{EQA}
	\TV\bigl( \ND(\xvs,\DV_{\GP}^{-2}), \ND(\xv,\HU^{-2}) \bigr)
	& \leq &
	\frac{1}{2} \Bigl( \| \DV_{\GP} (\xv - \xvs) \| + \sqrt{\tr \BB_{\GP}^{2}} \Bigr) .
\label{dferrfwbvf6nhdnghfkeiry}
\end{EQA}
However, 
Pinsker's inequality is only a general upper bound which applied to any two distributions \( P \) and \( Q \).
If \( P \) and \( Q \) are Gaussian, it might be too rough.
Particularly the use of \( \tr \BB^{2} \) is disappointing, this quantity is full dimensional even 
if each of \( \DV_{\GP}^{-2} \) and \( \HU^{-2} \) has a bounded trace.
Also dependence on \( \| \DV_{\GP} (\xv - \xvs) \| \) is very discouraging; see Remark~\ref{RTpostmeanc}.
\cite{Devroy2022} provides much sharper results, however, limited to the case of the same mean.
Even stronger results can be obtained if we restrict ourselves to the class \( \BBB_{el}(\R^{\dimp}) \) 
of elliptic sets \( A \) in \( \R^{\dimp} \) of the form
\begin{EQA}
	A
	&=&
	\bigl\{ \uv \in \R^{\dimp} \colon \| \QP (\uv - \xv) \| \leq \rr \bigr\}
\label{vhg4dfe5w3tfdf54wteg}
\end{EQA}
for some linear mapping \( \QP \colon \R^{\dimp} \to \R^{\dimq} \), \( \xv \in \R^{\dimp} \), and \( \rr > 0 \).
Given two symmetric \( \dimq \)-matrices \( \Sigma_{1},\Sigma_{2} \) and a vector \( \av \in \R^{\dimq} \), define 
\begin{EQA}
	\dist(\Sigma_{1},\Sigma_{2},\av)
	& \eqdef &
	\biggl( \frac{1}{\| \Sigma_{1} \|_{\Fr}} + \frac{1}{\| \Sigma_{2} \|_{\Fr}} \biggr)
	\biggl( \| \lamb_{1} - \lamb_{2} \|_{1} + \| \av \|^{2} \biggr),
\label{btegdhertrewfdgvfddffd}
\end{EQA}
where \( \lamb_{1} \) is the vector of eigenvalues of \( \Sigma_{1} \) arranged in the non-increasing order
and similarly for \( \lamb_{2} \).
\cite{GNSUl2017} stated the following bound for \( \gaussv_{1} \sim \ND(0,\Sigma_{1}) \) and \( \gaussv_{2} \sim \ND(0,\Sigma_{2}) \):
with an absolute constant \( \CONST \)
\begin{EQA}
	\sup_{\rr > 0}
	\Bigl| \P\bigl( \| \gaussv_{1} + \av \| \leq \rr \bigr) - \P\bigl( \| \gaussv_{2} \| \leq \rr \bigr) \Bigr|
	& \leq &
	\CONST  \, \dist(\Sigma_{1},\Sigma_{2},\av) 
\label{gttegderqwerwvfdhewfs}
\end{EQA}
provided that \( \| \Sigma_{k} \|^{2} \leq \| \Sigma_{k} \|_{\Fr}^{2}/3 \), \( k=1,2 \).
By the Weilandt--Hoffman inequality, 
\( \| \lamb_{1} - \lamb_{2} \|_{1} \leq \| \Sigma_{1} - \Sigma_{2} \|_{1} \) , see e.g. 
\cite{MarkusEng}.
Here \( \| M \|_{1} = \tr |M| = \sum_{j} |\lambda_{j}(M)| \) for a symmetric matrix \( M \) with eigenvalues \( \lambda_{j}(M) \).

\begin{theorem}
\label{TLaplaceTVin}
Assume the conditions of Theorem~\ref{TLaplaceTV} with \( \xvs \) being the maximizer of \( \lgd \) and 
\( \DV_{\GP}^{2} = - \nabla^{2} \lgd(\xvs) \).
For any \( \xv \) and \( \HU \),
it holds with \( \gaussv_{\HU} \sim \ND(0,\HU^{-2}) \)
\begin{EQA}
	&& \nquad
	\sup_{A \in \BBB(\R^{\dimp})} \bigl| \PfL(\Xv - \xvs \in A) - \P(\gaussv_{\HU} \in A) \bigr|
	\\
	& \leq &
	4(\err + \ex^{-\xx}) + \TV\bigl( \ND(\xv,\HU^{-2}), \ND(\xvs,\DV_{\GP}^{-2}) \bigr) 
\label{nhjrwsdgdehgdftregwbdctsh}
\end{EQA}
with \( \err = \err_{2} \), see \eqref{juytr90f2dzaryjhfyf}, or \( \err = \err_{3} \), see \eqref{5qw7dyf4e4354coefw9dufih}.

Furthermore, for \( \Xv \sim \PfL \) and \( \gaussv \sim \ND(0,\Id_{\dimp}) \), any linear mapping 
\( \QP \colon \R^{\dimp} \to \R^{\dimq} \),
it holds 
under \( 3 \| \QP \, \DV_{\GP}^{-2} \QP^{\T} \|^{2} \leq \| \QP \, \DV_{\GP}^{-2} \QP^{\T} \|_{\Fr}^{2} \)
\begin{EQA}
	&& \nquad
	\sup_{\rr > 0}
	\left| \PfL\bigl( \| \QP (\Xv - \xv) \| \leq \rr \bigr) 
	- \P\bigl( \| \QP \, \HU^{-1} \gaussv \| \leq \rr \bigr) 
	\right|
	\\
	& \leq &
	4(\err_{3} + \ex^{-\xx}) + \frac{\CONST}{\| \QP \, \DV_{\GP}^{-2} \QP^{\T} \|_{\Fr}}
	\left( 
		\| \QP (\DV_{\GP}^{-2} - \HU^{-2}) \QP^{\T} \|_{1} + \| \QP (\xv - \xvs) \|^{2} 
	\right) .
\label{jrwguvyr23jbviufdsfsdgf6w}
\end{EQA}
\end{theorem}

As a special case, we consider the use of the posterior mean 
\begin{EQA}
	\xvb
	& \eqdef &
	\frac{\int \xv \, \ex^{\lgd(\xv)} \, d\xv}{\int \ex^{\lgd(\xv)} \, d\xv}
\label{054rfgofe3rdgbrfty6ry}
\end{EQA}
in place of \( \xvs = \argmax_{\xv} \lgd(\xv) \).

\begin{theorem}
\label{TpostmeanLan}
Assume the conditions of Theorem~\ref{TpostmeanLa} and Theorem~\ref{TLaplaceTVin}. 
Then it holds for any linear mapping \( \QP \colon \R^{\dimp} \to \R^{\dimq} \) 
\begin{EQA}
	\sup_{\rr > 0}
	\left| \PfL\bigl( \| \QP (\Xv - \xvb) \| \leq \rr \bigr) 
	- \P\bigl( \| \QP \gaussv_{\GP} \| \leq \rr \bigr) 
	\right|
	& \leq & 
	4(\err_{3} + \ex^{-\xx}) + \frac{\CONST \| \QP (\xvb - \xvs) \|^{2}}{\| \QP \, \DV_{\GP}^{-2} \QP^{\T} \|_{\Fr}} \, ,
\label{jrwguvyr23jbviufdsfsdgf6wm}
\end{EQA}
where \( \| \QP (\xvb - \xvs) \| \) follows \eqref{hcdtrdtdehfdewdrfrhgyjufger} and \eqref{klu8gitfdgregfkhj7yt}.
\end{theorem}

The case \( \QP = \DVL \) is particularly transparent.
The obtained bound and \eqref{klu8gitfdgregfkhj7yt} of Corollary~\ref{CTpostmeanLa} yield
in view of \( \| \QP \, \DV_{\GP}^{-2} \QP^{\T} \|_{\Fr}^{2} = \tr(\DVL \, \DV_{\GP}^{-2} \DVL^{\T})^{2} \asymp \dimLG \) 
the following result.

\begin{corollary}
\label{CTpostmeanLab}
Under the conditions of Corollary~\ref{CTpostmeanLa}, it holds for \( \Xv \sim \PfL \)
\begin{EQA}
	\sup_{\rr > 0}
	\Bigl| \PfL\bigl( \| \DVL (\Xv - \xvb) \| \leq \rr \bigr) 
	- \P\bigl( \| \DVL \, \gaussv_{\GP} \| \leq \rr \bigr) 
	\Bigr|
	& \leq &
	\CONST \biggl( \sqrt{\frac{\dimLG^{3}}{n}} + \ex^{-\xx} \biggr) \, .
\label{klu8gitfdgregfkhj7yt3n}
\end{EQA}
\end{corollary}

The same bound holds with \( \QP = n^{1/2} \Id_{\dimp} \) in place of \( \DVL \) provided that 
\( \DVL^{2} \geq \CONSTi_{0} \, n \, \Id_{\dimp} \) for some fixed \( \CONSTi_{0} > 0 \).
We may conclude that the use of posterior mean \( \xvb \) in place of the posterior mode \( \xvs \) is justified 
under the same condition on critical dimension \( \dimLG^{3} \ll n \) as required for the main result about Gaussian 
approximation.

\Section{Bayesian optimization and iterated Laplace approximation}
Suppose that the point of maximum \( \xvs \) of the function \( \lgdL(\xv) \) is not known and has to be evaluated numerically.
We also assume that the function \( \lgdL(\cdot) \) is sufficiently smooth and can be efficiently computed at any point \( \xv \), 
however, its gradient is not available.
\cite{NeSp2016} offered a powerful gradient free method based on averaging the exponent \( \ex^{\lgdL(\xv)} \) w.r.t. to some Gaussian 
distribution for \( \xv \).
Inspired by the results of previous sections, we propose here a modified version of \cite{NeSp2016} based
on iterated Laplace approximations.
Starting from some \( \priord_{0} = \ND(\xv_{0},\GP_{0}^{-2}) \), we iteratively update 
the Gaussian prior \( \priord_{k} \sim \ND(\xv_{k},\GP_{k}^{-2}) \) and try to numerically assess the corresponding posterior obtained by normalization of 
\( \ex^{\lgd_{k}(\xv)} \) for \( \lgd_{k}(\xv) = \lgdL(\xv) - \| \GP_{k} (\xv - \xv_{k}) \|^{2}/2 \).
This posterior is known to be nearly Gaussian by Theorem~\ref{TLaplaceTV},  
its mean \( \xv_{k+1} \) can be efficiently estimated by Monte-Carlo or quasi Monte-Carlo sampling; see e.g. 
\cite{SSW:2020}.
We use \( \xv_{k+1} \) for building the Gaussian prior for the step \( k+1 \).
Moreover, due to our results, the value \( \xv_{k+1} \) is nearly the maximizer of the corresponding quadratic approximation of 
\( \lgd_{k}(\xv) \) in the vicinity of \( \xv_{k} \).
One can say that computing the posterior mean \( \xv_{k+1} \) mimics well a step of the second order Newton-Raphson 
optimization method but does not require computing the gradient and the inverse Hessian matrix. 
The prior precision matrix \( \GP_{k}^{2} \) can be taken in the form \( \lambda_{k}^{-1} \GP_{0}^{2} \), where 
\( \lambda_{k} \) replaces the step size, \( \lambda_{k} \to 0 \) as \( k \) increases.
The procedure can be written as follows.

\begin{algorithm}
\caption{Laplace iterations}\label{alg:Lapl}
\begin{algorithmic}
\State
1. Start with \( k=0 \); Fix \( \xv_{0} \) and \( \GP_{0}^{2} \);

\State
2. Draw a sample \( (\xv^{(m)})_{m \leq M} \) from \( \ND(\xv_{k},\GP_{k}^{-2}) \).
For each \( m \leq M \), compute \( \weight^{(m)} = \ex^{\lgdL(\xv^{(m)})} \) and update the sums \( W_{k} = \sum_{m} \weight^{(m)} \)
and \( S_{k} = \sum_{m} \weight^{(m)} \xv^{(m)} \).

\State 
3. Compute the value \( \xv_{k+1} = S_{k}/W_{k} \) as the next prior mean;
update the posterior precision matrix \( \GP_{k+1}^{2} = a \GP_{k}^{2} \).

\State
4. Increase \( k \) by one and repeat from Step 2. 

\State 
5. Iterate until convergence.
\end{algorithmic}
\end{algorithm}

The starting guess is important as well the choice of the step multiplier \( a \). 
Moreover, one can use the variable multiplier \( a_{k} \) and incorporate the posterior covariance 
for variance reduction schemes.
One can make a rigorous analysis of the convergence of this algorithm similarly to \cite{NeSp2016}.

The proposed procedure is closely related to ensemble Kalman filtering technique; 
see e.g. \cite{SchiSt2017}, \cite{Reich2022} and references therein.

\Chapter{Laplace approximation for non-linear inverse problem}
\label{Slaplnonlin}
Let \( \regrfv(\xv) = \bigl( \regrf_{i}(\xv), \, i \leq n \bigr) \in \R^{n} \) be a nonlinear 
mapping (operator) of the source signal \( \xv \in \R^{\dimp} \) to the target space \( \R^{n} \).
We consider the problem of inverting the relation 
\( \zv = \regrfv(\xv) \): given an image vector \( \zv \in \R^{n} \), recover the corresponding 
source \( \xv \in \R^{\dimp} \).
This leads to the nonlinear least square problem of maximizing 
the function \( \lgdL(\cdot) \) of the form 
\begin{EQA}
	\lgdL(\xv) 
	&=& 
	- \frac{1}{2} \| \zv - \regrfv(\xv) \|^{2} .
\label{o3kv6ejfgiuhbdfjk}
\end{EQA}
Given a Gaussian prior \( \ND(\xv_{0},\GP^{-2}) \),
we consider Laplace's approximation for the penalized function 
\begin{EQA}
	\lgd(\xv) 
	&=& 
	- \frac{1}{2} \| \zv - \regrfv(\xv) \|^{2} - \frac{1}{2} \| \GP (\xv - \xv_{0}) \|^{2} \, .
\label{dghoswpewgeuw9hvadouvhoisd}
\end{EQA}
Define 
\begin{EQA}
	\xvs_{\GP}
	& \eqdef &
	\argmax_{\xv} \lgd(\xv)
	=
	\argmin_{\xv} \bigl\{ \| \zv - \regrfv(\xv) \|^{2} + \| \GP (\xv - \xv_{0}) \|^{2} \bigr\} .
\label{gh0ehrb9hqeovnsflkvnouwhg}
\end{EQA}
%
%
Laplace's approximation requires weak concavity of \( \lgdL(\xv) \); see 
\nameref{LLf0ref}.
A sufficient condition is \( \IFVL(\xv) \geq 0 \), where
\begin{EQA}
	\IFVL(\xv)
	& \eqdef &
	- \nabla^{2} \lgdL(\xv) 
	=
	\sumi \nabla \regrf_{i}(\xv) \, \nabla \regrf_{i}(\xv)^{\T}
	+ \sumi \{ \regrf_{i}(\xv) - z_{i} \} \, \nabla^{2} \regrf_{i}(\xv) .
	\qquad
\label{bvudufgheuygfdtcsdwedeuf}
\end{EQA}
%
%
However, such a condition seems can be hard to ensure for all \( \xv \in \R^{\dimp} \)
unless \( \regrfv(\cdot) \) is linear.
%
Instead we restrict ourselves to some subset \( \XX_{0} \subset \R^{\dimp} \)
which can be viewed as a concentration set of the Gaussian measure 
\( \ND(\xv_{0},\GP^{-2}) \). 
As the prior mass of the complement of \( \XX_{0} \) is exponentially small,
when sampling from the prior \( \ND(\xv_{0},\GP^{-2}) \),
one would need exponentially many samples to hit any set \( A \) outside of \( \XX_{0} \).
Therefore, we restrict the prior to the set \( \XX_{0} \) by skipping all draws of \( \xv \) with 
\( \xv \not\in \XX_{0} \).
In what follows we also assume that \( \xv_{0} \) is a reasonable guess ensuring that the target
\( \xvs_{\GP} \) belongs to \( \XX_{0} \).
Then the result about Laplace approximation applies even after restricting the parameter space to \( \XX_{0} \).
This local set \( \XX_{0} \) is defined as elliptic concentration set for the Gaussian prior \( \ND(\xv_{0},\GP^{-2}) \) 
in the form 
\begin{EQA}
	\XX_{0}
	&=&
	\bigl\{ \xv \colon \| \QP (\xv - \xv_{0}) \| \leq \rrin \bigr\} ,
\label{ocnfhyetetgdgrwebgdgtdgtL}
\end{EQA}
where \( \QP \) is a linear operator in \( \R^{\dimp} \) and
\( \rrin \) is fixed to ensure that the most of prior mass is within \( \XX_{0} \).
Theorem~\ref{TexpbLGA} suggests 
\begin{EQA}
	\rrin 
	&=& 
	\sqrt{\tr (\QP^{2} \, \GP^{-2})} + \sqrt{2\xx \, \| \QP \, \GP^{-2} \QP \|} 
\label{987655rtfvsdfghjoiuytr2w3eL}
\end{EQA}
yielding \( \P(\GP^{-1} \gaussv \not\in \XX_{0}) \leq \ex^{-\xx} \) with \( \gaussv \) standard normal.
A simple choice \( \QP = \Id_{\dimp} \) yields \( \XX_{0} \) in form of a ball around \( \xv_{0} \)
with the radius of order \( \sqrt{\tr \GP^{-2}} \).
A proper choice of \( \GP^{2} \geq \Id_{\dimp} \) has to ensure a small value of \( \tr (\GP^{-2}) \).
Note that the operator \( \QP \) can be scaled by any factor.
Then the radius \( \rrin \) in \eqref{987655rtfvsdfghjoiuytr2w3eL} will be scaled correspondingly.

Introduce for any \( \xv \in \XX_{0} \) a \( \dimp \)-symmetric matrix 
\begin{EQA}
	\DVLb^{2}(\xv)
	&=&
	\nabla \regrfv(\xv) \, \nabla \regrfv(\xv)^{\T}
	=
	\sumi \nabla \regrf_{i}(\xv) \, \nabla \regrf_{i}(\xv)^{\T} \, . 
\label{tsfdc6frdw6red6wt7qagdudwgL}
\end{EQA}
This is the first term in expansion \eqref{bvudufgheuygfdtcsdwedeuf}
corresponding to a linear approximation of \( \regrfv(\cdot) \) at \( \xv \).
Injectivity of \( \regrfv(\cdot) \) means that \( \DVLb^{2}(\xv) \) is positive definite and well conditioned. 
If \( n^{-1} \DVLb^{2}(\xv_{0}) \leq \GP^{2} \) then
the choice \( \QP = \DVLb(\xv_{0}) \) in \eqref{ocnfhyetetgdgrwebgdgtdgtL} is our alternative to 
\( \QP = \Id_{\dimp} \).

Below we assume the following regularity condition.

\begin{description}
    \item[\label{LLS2ref} \( \bb{(\nabla \regrfv)} \)]
      \emph{ Let the local set \( \XX_{0} \) be defined by \eqref{ocnfhyetetgdgrwebgdgtdgtL} and \eqref{987655rtfvsdfghjoiuytr2w3eL} with \( \QP = \DVLbin \eqdef \DVLb(\xv_{0}) \).
      For some fixed \( \CONSTi_{2} \), \( \CONSTi_{n} \), 
      and any \( \xv \in \XX_{0} \), \( \uv \in \R^{\dimp} \), it holds
\begin{EQA}
\label{hjwetwee3ee5555t566dsL}
	\sumi \left| \bigl\langle \nabla^{2} \regrf_{i}(\xv), \uv^{\otimes 2} \bigr\rangle \right|
	& \leq &
	\CONSTi_{2} \, \sumi \bigl\langle \nabla \regrf_{i}(\xv), \uv \bigr\rangle^{2}
	=
	\CONSTi_{2} \, \| \DVLb(\xv) \, \uv \|^{2} 
	\, ,
	\\
	\max_{i \leq n} \, \nabla \regrf_{i}(\xv) \, \nabla \regrf_{i}(\xv)^{\T}
	& \leq &
	\frac{\CONSTi_{n}^{2}}{n} \, \sumi \nabla \regrf_{i}(\xv) \, \nabla \regrf_{i}(\xv)^{\T} 
	=
	\frac{\CONSTi_{n}^{2}}{n} \, \DVLb^{2}(\xv) \, ,
\label{wfegy7r5qrw35edfhgdyfysL}
\end{EQA}
and with some \( \CONSTi_{0} \geq 1 \)
\begin{EQA}
	\DVLb^{2}(\xv)
	& \leq &
	\CONSTi_{0}^{2} \, \DVLbin^{2} .
\label{fchghiu87686574e5rtyyuL}
\end{EQA}
}
\end{description}

\noindent
Condition \nameref{LLS2ref} only requires some local regularity of the \( \regrf_{i}(\xv) \)'s 
and injectivity of the gradient mapping \( \nabla \regrfv(\xv) \) within \( \XX_{0} \).
The shape of \( \XX_{0} \) is defined by the prior covariance \( \GP^{-2} \) which has to be selected 
to ensure that the set \( \XX_{0} \) is a local set.
More precisely, we connect the size of \( \XX_{0} \) and local regularity of \( \regrfv(\cdot) \) by the following 
condition.

\begin{description}
    \item[\label{rrinref} \( \bb{(\rrin)} \)]
      \emph{With \( \XX_{0} \) defined by \eqref{ocnfhyetetgdgrwebgdgtdgtL} and \eqref{987655rtfvsdfghjoiuytr2w3eL}
      and \( \CONSTi_{n} \), \( \CONSTi_{0} \), and \( \CONSTi_{2} \) from \nameref{LLS2ref}, it holds
\begin{EQA}
	\frac{2 \CONSTi_{n} \, \CONSTi_{0} \, \CONSTi_{2} \, \rrin}{\sqrt{n}} 
	& \leq & 
	\frac{1}{4} \,\, .
\label{yf2qsxysdrqwfhaqaraqfL}
\end{EQA}
      }           
\end{description}

To get a feeling of this condition, consider the case with \( \DVLb^{2} = n \Id_{\dimp} \).
Then \eqref{987655rtfvsdfghjoiuytr2w3eL} yields \( \rrin^{2} \asymp n \tr \GP^{-2} \) 
and \eqref{yf2qsxysdrqwfhaqaraqfL} requires \( \tr \GP^{-2} \ll 1 \), that is, 
the prior \( \ND(\xv_{0},\GP^{2}) \) is supported to a small vicinity of \( \xv_{0} \).


Now we introduce the key condition of ``warm start'' which informally means a good choice 
of the starting point \( \xv_{0} \) ensuring that the local vicinity 
\( \XX_{0} \) of \( \xv_{0} \) contains a point \( \xv \) with 
\( \| \regrfv(\xv) - \zv \|_{\infty} \) is small.

\begin{description}
    \item[\label{theta0ref} \( \bb{(\xv_{0})} \)]
      \emph{With \( \XX_{0} \) from \eqref{ocnfhyetetgdgrwebgdgtdgtL} 
      and \eqref{987655rtfvsdfghjoiuytr2w3eL} for \( \QP = \DVLbin \), 
      and \( \CONSTi_{2} \) from \eqref{hjwetwee3ee5555t566dsL}, it holds 
\begin{EQA}
	\inf_{\xv \in \XX_{0}} \| \regrfv(\xv) - \zv \|_{\infty}
	=
	\inf_{\xv \in \XX_{0}} \,\, \max_{i \leq n}\, | \regrf_{i}(\xv) - z_{i} |
	& \leq &
	\smlc_{0} \, ,
	\qquad
	\CONSTi_{2} \, \smlc_{0}
	\leq 
	{1}/{4} \, .
\label{yfswyfwsyfywsaysfyasedsqarL}
\end{EQA}      
}
\end{description}

Later we show that \nameref{LLS2ref}, \nameref{rrinref}, and \nameref{theta0ref} ensure  
\begin{EQA}
	\IFVL(\xv)
	& \geq &
	(1 - \delta ) \, \DVLb^{2}(\xv),  \qquad \xv \in \XX_{0} \, ,
\label{p160ch28mcrwiwi4fwkgey}
\end{EQA}
for some \( \delta \leq 3/4 \), and hence, 
\( \lgdL(\xv) \) is strongly concave for \( \xv \in \XX_{0} \).

Now we consider the point \( \xvs_{\GP} \) from \eqref{gh0ehrb9hqeovnsflkvnouwhg}
and its local vicinity.
It is important to secure that \( \xvs_{\GP} \) is within \( \XX_{0} \).
This appears to be automatically fulfilled if the true solution 
\( \xvs = \regrfv^{-1}(\zv) \) is within \( \XX_{0} \); see Lemma~\ref{Lnonlintrue}.
Now we are back to the setup with 
a strongly concave function \( \lgdL(\xv) \) for \( \xv \in \XX_{0} \).
Similarly to the general case, define 
\( \DVL^{2} = \DVL^{2}(\xvs_{\GP}) = \IFVL(\xvs_{\GP}) \) and
\( \DV_{\GP}^{2} = \DV_{\GP}^{2}(\xvs_{\GP}) = \IFV_{\GP}(\xvs_{\GP}) \); see \eqref{bvudufgheuygfdtcsdwedeuf}.
The effective dimension is given by \( \dimG = \tr( \DVL^{2} \, \DV_{\GP}^{-2} ) \);
see \eqref{dAdetrH02Hm2}.
The local vicinity \( \UVL_{\GP} \) of \( \xvs_{\GP} \) can be defined by the rule \eqref{UvTDunm12spT}.
Alternatively one can use 
\( \DVLb^{2} = \DVLb^{2}(\xvs_{\GP}) \) instead of \( \DVL^{2} \):
\begin{EQA}
	\UVLb
	&=&
	\bigl\{ \xv \colon \| \DVLb (\xvs_{\GP} - \xv) \| \leq \rrLb \bigr\} .
\label{UvTDunm12spTwfw535}
\end{EQA}
The radius \( \rrLb \) has to be adjusted to ensure that the such defined vicinity 
of \( \xvs \) is not smaller than \( \UVL \) from \eqref{UvTDunm12spT}:
\( \UVL_{\GP} \subseteq \UVLb \).
Relation \eqref{p160ch28mcrwiwi4fwkgey} suggests using \( \rrLb = 2 \rrL \).

For applying the general results of Theorems~\ref{TLaplaceTV} through \ref{TLaplaceTV34}, 
the function \( \elll(\xv) = - \| \regrfv(\xv) - \zv \|^{2}/2 \) has to be sufficiently smooth.
We now present some sufficient conditions 
in terms of the functions \( \regrf_{i}(\xv) \).
These conditions extend \nameref{LLS2ref}.

\begin{description}
    \item[\label{LLS3ref} \( \bb{(\nabla^{3} \regrf)} \)]
      \emph{ For some \( \CONSTi_{\GP} \) and \( \DVb^{2} = \DVb^{2}(\xvs_{\GP}) \)
\begin{EQA}
	\DVb^{2}(\xv)
	& \leq &
	\CONSTi_{\GP} \, \DVb^{2} ,
	\qquad
	\xv \in \UVLb \, ;
\label{fchghiu87686574e5rtyyuLL}
\end{EQA}
      cf. \eqref{fchghiu87686574e5rtyyuL}.
      Furthermore, for \( \CONSTi_{3} \geq 0 \), 
      it holds uniformly over \( \xv \in \UVLb \) and \( \uv \in \R^{\dimp} \)
\begin{EQA}
	\sumi\left| \bigl\langle \nabla^{3} \regrf_{i}(\xv), \uv^{\otimes 3} \bigr\rangle \right|
	& \leq &
	\CONSTi_{3} \, \sumi \left| \bigl\langle \nabla \regrf_{i}(\xv), \uv \bigr\rangle \right|^{3} \, ,
	\\
	\sumi\left| 
		\bigl\langle \nabla \regrf_{i}(\xv), \uv \bigr\rangle 
		\bigl\langle \nabla^{2} \regrf_{i}(\xv), \uv^{\otimes 2} \bigr\rangle 
	\right|
	& \leq &
	\CONSTi_{3} \, \sumi \left| \bigl\langle \nabla \regrf_{i}(\xv), \uv \bigr\rangle \right|^{3} \, .
\label{hjwetwee3ee5555t566dL}
\end{EQA}
}
    \item[\label{LLS4ref} \( \bb{(\nabla^{4} \regrf)} \)]
      \emph{ For some  \( \CONSTi_{4} \), it holds uniformly over 
      \( \xv \in \UVLb \) and \( \uv \in \R^{\dimp} \)
\begin{EQA}
	\sumi\left| \bigl\langle \nabla^{4} \regrf_{i}(\xv), \uv^{\otimes 4} \bigr\rangle \right|
	& \leq &
	\CONSTi_{4} \, \sumi \left| \bigl\langle \nabla \regrf_{i}(\xv), \uv \bigr\rangle \right|^{4} \, ,
	\\
	\sumi\left| 
		\bigl\langle \nabla \regrf_{i}(\xv), \uv \bigr\rangle 
		\bigl\langle \nabla^{3} \regrf_{i}(\xv), \uv^{\otimes 3} \bigr\rangle 
	\right|
	& \leq &
	\CONSTi_{4} \, \sumi \left| \bigl\langle \nabla \regrf_{i}(\xv), \uv \bigr\rangle \right|^{k} \, ,
	\\
	\sumi 
		\bigl\langle \nabla^{2} \regrf_{i}(\xv), \uv^{\otimes 2} \bigr\rangle^{2}
	& \leq &
	\CONSTi_{4} \, \sumi \left| \bigl\langle \nabla \regrf_{i}(\xv), \uv \bigr\rangle \right|^{4} \, .
\label{hjwetwee3ee5555t566d4L}
\end{EQA}
}
\end{description}

We are now prepared to state the main result about Laplace's approximation
which is a straightforward application of Theorem~\ref{TLaplaceTV} to the measure \( \PfL \) restricted to \( \XX_{0} \). 

\begin{theorem}
\label{TLaplnonlin}
Consider the function \( \lgdL(\xv) = - \| \regrfv(\xv) - \zv \|^{2}/2 \); see \eqref{dghoswpewgeuw9hvadouvhoisd}.
Given \( (\xv_{0}, \GP^{2}) \), let \( \XX_{0} \) be defined by \eqref{ocnfhyetetgdgrwebgdgtdgtL} and \eqref{987655rtfvsdfghjoiuytr2w3eL} with \( \QP = \DVLbin \).
Assume \nameref{LLS2ref}, \nameref{rrinref}, \nameref{theta0ref} with \( \smlc_{0} = 0 \), and \nameref{LLS3ref}. 
Let \( \xvs_{\GP} \) be from \eqref{gh0ehrb9hqeovnsflkvnouwhg}, \( \DVL^{2} = \IFV(\xvs_{\GP}) \); see 
\eqref{bvudufgheuygfdtcsdwedeuf}, and \( \DV_{\GP}^{2} = \DVL^{2} + \GP^{2} \).
If also \( \hmax_{3} \, \rrL \, n^{-1/2} \leq 3/4 \) with
\begin{EQA}
	\hmax_{3}
	& \eqdef &
	4 \CONSTi_{\GP}^{3/2} \, \CONSTi_{3} \, \CONSTi_{n}
\label{cuv6d54d3ws3tvfyfy}
\end{EQA}
then with \( \gaussv_{\GP} \sim \ND(0,\DV_{\GP}^{-2}) \)
\begin{EQA}
	\sup_{A \in \BBB(\R^{\dimp})} \bigl| \PfL(\Xv - \xvs_{\GP} \in A \cond \XX_{0}) - \P(\gaussv_{\GP} \in A) \bigr|
	& \leq &
	2 \hmax_{3} \, \sqrt{\frac{(\dimLG+1)^{3}}{n}} + 4 \ex^{-\xx} \, .
	\qquad
\label{cbc5dfedrdwewwgerg3NL}
\end{EQA}
If also \nameref{LLS4ref} is satisfied then 
\begin{EQA}
	\sup_{A \in \BBB_{s}(\R^{\dimp})} \bigl| \PfL(\Xv - \xvs_{\GP} \in A \cond \XX_{0}) - \P(\gaussv_{\GP} \in A) \bigr|
	& \leq &	
	\frac{\hmax_{3}^{2} \, (\dimG + 2)^{3} + 2 \hmax_{4} (\dimLG + 1)^{2}}{2n} + 4 \ex^{-\xx} \, 
	\qquad
\label{scdugfdwyd2wywy26e6de4NL}
\end{EQA}
with 
\begin{EQA}
	\hmax_{4}
	& \eqdef &
	8 \CONSTi_{\GP}^{2} \, \CONSTi_{4} \, \CONSTi_{n} \, .
\label{wed9u8wefuegfeyfgyfgey}
\end{EQA}
\end{theorem}

\appendix

\Chapter{Tools and proofs} 
\label{Slaplaceproofs}

Here we collect the proofs of the main results and some useful technical statements 
about the error of Laplace approximation.
Below we write \( \xv \) instead of \( \xvs \).
Note that after passing to representation \eqref{IIgfifxutudufxutdu}, 
many results below apply to any \( \xv \), not necessarily for \( \xv = \xvs \).
We only use \( \DV_{\GP}^{2} = - \nabla^{2} \lgd(\xv) \) and 
\( \dltwb \) instead of \( \dltwb(\xv) \).
Everywhere we assume the local set \( \UVL \) to be fixed by \eqref{UvTDunm12spT}.
We separately study the integrals over \( \UVL \) and over its complement. 
The local error of approximation is measured by
\begin{EQA}
	\err
	=
	\err(\UVL)
	& \eqdef &
	\biggl| 
	\frac{\int_{\UVL} \ex^{\lgd(\xv;\uv)} \, g(\uv) \, d\uv - \int_{\UVL} \ex^{- \| \DV_{\GP} \uv \|^{2}/2} \, g(\uv) \, d\uv} 
		 {\int \ex^{- \| \DV_{\GP} \uv \|^{2}/2} d\uv} 
	\biggr| \, .
\label{errdefdiUaHu}
\end{EQA}
As a special case with \( g(\uv) \equiv 1 \) we obtain an approximation of the denominator in \eqref{IIgfifxutudufxutdu}.
In addition, we have to bound the tail integrals
\begin{EQ}[rcl]
	\rho
	=
	\rho(\UVL)
	& \eqdef &
	\frac{\int \Ind(\uv \not\in \UVL) \, \ex^{\lgd(\xv;\uv)} \, d\uv}{\int \ex^{- \| \DV_{\GP} \uv \|^{2}/2} \, d\uv} \, ,
	\\
	\rho_{\GP}
	=
	\rho_{\GP}(\UVL)
	& \eqdef &
	\frac{\int \Ind(\uv \not\in \UVL) \, \ex^{- \| \DV_{\GP} \uv \|^{2}/2} \, d\uv}{\int \ex^{- \| \DV_{\GP} \uv \|^{2}/2} \, d\uv} \, .
\label{rhfiUaceHu22m}
\end{EQ}
The analysis will be split into several steps.

\Section{Overall error of Laplace approximation}
First we show how to seam together the error \( \err \) of local approximation and the bounds for the tail integrals
\( \rho \) and \( \rho_{\GP} \); see \eqref{rhfiUaceHu22m}.

\begin{proposition}
\label{PunbintLapl}
Suppose that for a function \( g(\uv) \) with \( g(\uv) \in [0,1] \) and some \( \err , \err_{g} \)
\begin{EQA}[rcl]
	\biggl| 
		\frac{\int_{\UVL} \ex^{\lgd(\xv;\uv)} \, d\uv - \int_{\UVL} \ex^{- \| \DV_{\GP} \uv \|^{2}/2} \, d\uv}			 
			 {\int \ex^{- \| \DV_{\GP} \uv \|^{2}/2} \, d\uv} 
	\biggr|
	& \leq &
	\err \, ,
\label{erifxudieHu22}
	\\
\label{erifxudieHu22g}
	\biggl| 
		\frac{\int_{\UVL} g(\uv) \, \ex^{\lgd(\xv;\uv)} \, d\uv - \int_{\UVL} g(\uv) \, \ex^{- \| \DV_{\GP} \uv \|^{2}/2} \, d\uv}			 {\int \ex^{- \| \DV_{\GP} \uv \|^{2}/2} \, d\uv} 
	\biggr|
	& \leq &
	\err_{g} \, .
\end{EQA}
Then with \( \rho \) and \( \rho_{\GP} \) from \eqref{rhfiUaceHu22m}
\begin{EQ}[rcl]
	\frac{\int g(\uv) \, \ex^{\lgd(\xv;\uv)} \, d\uv}{\int \ex^{\lgd(\xv;\uv)} \, d\uv} 
	& \leq &
	\frac{1}{1 - \rho_{\GP} - \err} \,\, 
	\frac{\int g(\uv) \, \ex^{- \| \DV_{\GP} \uv \|^{2}/2} \, d\uv}{\int \ex^{- \| \DV_{\GP} \uv \|^{2}/2} \, d\uv}
	+ \frac{\rho + \err_{g}}{1 - \rho_{\GP} - \err} \, ,
	\qquad
	\\
	\frac{\int g(\uv) \, \ex^{\lgd(\xv;\uv)} \, d\uv}{\int \ex^{\lgd(\xv;\uv)} \, d\uv} 
	& \geq &
	\frac{1}{1 + \rho + \err} \,\, 
	\frac{\int g(\uv) \, \ex^{- \| \DV_{\GP} \uv \|^{2}/2} \, d\uv}{\int \ex^{- \| \DV_{\GP} \uv \|^{2}/2} \, d\uv}
	- \frac{\rho_{\GP} + \err_{g}}{1 + \rho + \err} \, .
	\qquad
\label{igefxumiguexHu22}
\end{EQ}
\end{proposition}

\begin{proof}
It follows from \eqref{erifxudieHu22} 
\begin{EQA}
	\int \ex^{\lgd(\xv;\uv)} \, d\uv
	& \geq &
	\int_{\UVL} \ex^{\lgd(\xv;\uv)} \, d\uv
	\geq 
	\int_{\UVL} \ex^{- \| \DV_{\GP} \uv \|^{2}/2} \, d\uv - \err \int \ex^{- \| \DV_{\GP} \uv \|^{2}/2} \, d\uv
	\\
	& \geq &
	(1 - \err - \rho_{\GP}) \int \ex^{- \| \DV_{\GP} \uv \|^{2}/2} \, d\uv ,
\label{1erriexmH22du22}
	\\
	\int \ex^{\lgd(\xv;\uv)} \, d\uv
	& \leq &
	\int_{\UVL} \ex^{\lgd(\xv;\uv)} \, d\uv + \rho \int \ex^{- \| \DV_{\GP} \uv \|^{2}/2} \, d\uv
	\\
	& \leq &
	(1 + \err + \rho) \int \ex^{- \| \DV_{\GP} \uv \|^{2}/2} \, d\uv .
\label{1erriexmH22du22u}
\end{EQA}
Similarly for \( g(\uv) \geq 0 \)
\begin{EQA}
	\int g(\uv) \, \ex^{\lgd(\xv;\uv)} \, d\uv
	& \geq &
	\int_{\UVL} g(\uv) \, \ex^{- \| \DV_{\GP} \uv \|^{2}/2} \, d\uv - \err_{g} \int \ex^{- \| \DV_{\GP} \uv \|^{2}/2} \, d\uv
	\\
	& \geq &
	\int g(\uv) \, \ex^{- \| \DV_{\GP} \uv \|^{2}/2} \, d\uv
	- (\rho_{\GP} + \err_{g}) \int \ex^{- \| \DV_{\GP} \uv \|^{2}/2} \, d\uv ,
	\\
	\int g(\uv) \, \ex^{\lgd(\xv;\uv)} \, d\uv
	& \leq &
	\int_{\UVL} g(\uv) \, \ex^{\lgd(\xv;\uv)} \, d\uv + \rho \int \ex^{- \| \DV_{\GP} \uv \|^{2}/2} \, d\uv
	\\
	& \leq &
	\int g(\uv) \, \ex^{- \| \DV_{\GP} \uv \|^{2}/2} \, d\uv + (\rho + \err_{g}) \int \ex^{- \| \DV_{\GP} \uv \|^{2}/2} \, d\uv \, .
\label{iguexmHu22dreieH}
\end{EQA}
Putting together all these bounds yields \eqref{igefxumiguexHu22}.
\end{proof}

The next corollary is straightforward.

\begin{corollary}
\label{CPunbintLapl}
Let \( \rho_{\GP} \leq \rhos \), \( \rho \leq \rhos \); see \eqref{rhfiUaceHu22m}.
Let also for a function \( g(\uv) \) with \( |g(\uv)| \leq 1 \), \eqref{erifxudieHu22}, \eqref{erifxudieHu22g} hold with
\( \err_{g} \leq \err \).
If \( \err + \rhos \leq 1/2 \) then 
\begin{EQA}
	\left| 
		\frac{\int g(\uv) \, \ex^{\lgd(\xv;\uv)} \, d\uv}{\int \ex^{\lgd(\xv;\uv)} \, d\uv} 
		- \frac{\int g(\uv) \, \ex^{- \| \DV_{\GP} \uv \|^{2}/2} \, d\uv}{\int \ex^{- \| \DV_{\GP} \uv \|^{2}/2} \, d\uv}
	\right|
	& \leq &
	\frac{2 (\rhos + \err)}{1 - \rhos - \err} 
	\leq 
	4 (\rhos + \err)	\, . 
\label{2rherIHg2re}
\end{EQA}
\end{corollary}


Sometimes we need an extension to the case of an unbounded function \( g \).
This particularly arises when evaluating the moment of the posterior; see Theorem~\ref{TpostmeanLa}.
The next result corresponds to estimation of posterior mean with a linear function \( g \) 
and posterior variance with \( g \) quadratic.

\begin{proposition}
\label{PunbintLapg}
Given a function \( g(\uv) \), assume \eqref{erifxudieHu22}, \eqref{erifxudieHu22g}, and define
\begin{EQA}
	\rho_{g}
	& \eqdef &
	\frac{\int \Ind(\uv \not\in \UVL) \, | g(\uv) | \, \ex^{\lgd(\xv;\uv)} \, d\uv}{\int \ex^{- \| \DV_{\GP} \uv \|^{2}/2} \, d\uv} 
	 \, ,
	\\
	\rho_{\GP,g}
	& \eqdef &
	\frac{\int \Ind(\uv \not\in \UVL) \, | g(\uv) | \, \ex^{- \| \DV_{\GP} \uv \|^{2}/2} \, d\uv}{\int \ex^{- \| \DV_{\GP} \uv \|^{2}/2} \, d\uv} \, ,
\label{rhfiUaceHu22mg}
\end{EQA}
while \( \rho \) and \( \rho_{\GP} \) are given in \eqref{rhfiUaceHu22m}.
Then with \( \II_{\GP}(g) = \E g(\gaussv_{\GP}) \), \( \gaussv_{\GP} \sim \ND(0,\DV_{\GP}^{-2}) \), 
\begin{EQA}
	&& \nquad
	\left| 
		\frac{\int g(\uv) \, \ex^{\lgd(\xv;\uv)} \, d\uv}{\int \ex^{\lgd(\xv;\uv)} \, d\uv} 
		- \frac{\int g(\uv) \, \ex^{- \| \DV_{\GP} \uv \|^{2}/2} \, d\uv}{\int \ex^{- \| \DV_{\GP} \uv \|^{2}/2} \, d\uv}
	\right|
	\leq 
	\frac{\rho_{g} + \rho_{\GP,g} + \err_{g}}{1 - \rho_{\GP} - \err}
	+ \frac{| \II_{\GP}(g) | \, (\rho + \err)}{1 - \rho_{\GP} - \err} \, . 
\label{2rherIHg2red}
\end{EQA}
In particular, if \( \int g(\uv) \, \ex^{-\| \DV_{\GP} \uv \|^{2}/2} \, d\uv = 0 \) then
\begin{EQA}[rcl]
	\left| \frac{\int g(\uv) \, \ex^{\lgd(\xv;\uv)} \, d\uv}{\int \ex^{\lgd(\xv;\uv)} \, d\uv} \right| 
	& \leq &
	\frac{\rho_{g} + \rho_{\GP,g} + \err_{g}}{1 - \rho_{\GP} - \err} \, .
\label{igefxumiguexHu22g}
\end{EQA}
\end{proposition}

\begin{proof}
Suppose that \( \II_{\GP}(g) \geq 0 \).
Then 
\begin{EQA}
	&& \nquad
	\left| 
		\frac{\int g(\uv) \, \ex^{\lgd(\xv;\uv)} \, d\uv}{\int \ex^{\lgd(\xv;\uv)} \, d\uv} 
		- \frac{\int g(\uv) \, \ex^{- \| \DV_{\GP} \uv \|^{2}/2} \, d\uv}{\int \ex^{- \| \DV_{\GP} \uv \|^{2}/2} \, d\uv}
	\right|
	\\
	& \leq &
	\left| 
		\frac{\int g(\uv) \, \ex^{\lgd(\xv;\uv)} \, d\uv}{\int \ex^{\lgd(\xv;\uv)} \, d\uv} 
		- \frac{\int g(\uv) \, \ex^{- \| \DV_{\GP} \uv \|^{2}/2} \, d\uv}{\int \ex^{\lgd(\xv;\uv)} \, d\uv}
	\right|
	+ \II_{\GP}(g)
	\left| 
		\frac{\int \ex^{- \| \DV_{\GP} \uv \|^{2}/2} \, d\uv}{\int \ex^{\lgd(\xv;\uv)} \, d\uv} - 1
	\right| .
\label{igefxvv22igh}
\end{EQA}
By definitions
\begin{EQA}
	&& \nquad
	\left| 
		\int g(\uv) \, \ex^{\lgd(\xv;\uv)} \, d\uv - \int g(\uv) \, \ex^{-\| \DV_{\GP} \uv \|^{2}/2} \, d\uv
	\right|
	\\
	& \leq &
	\left| 
		\int_{\UVL} g(\uv) \, \ex^{\lgd(\xv;\uv)} \, d\uv - \int_{\UVL} g(\uv) \, \ex^{-\| \DV_{\GP} \uv \|^{2}/2} \, d\uv
	\right|
	\\
	&&
	+ \, \left| 
		\int \Ind(\uv \not\in \UVL) \, g(\uv) \, \ex^{\lgd(\xv;\uv)} \, d\uv 
	\right|
	+ \left| 
		\int \Ind(\uv \not\in \UVL) \, g(\uv) \, \ex^{-\| \DV_{\GP} \uv \|^{2}/2} \, d\uv
	\right|	
	\\
	& \leq &
	\bigl( \rho_{g} + \rho_{\GP,g} + \err_{g} \bigr) \int \ex^{-\| \DV_{\GP} \uv \|^{2}/2} \, d\uv
\label{lrguegfxumHu22}
\end{EQA}
and the assertion follows by \eqref{1erriexmH22du22} and \eqref{1erriexmH22du22u}.
\end{proof}

\ifunivariate{}{
\Section{Local approximation. Univariate case}
This section presents a bound of local approximation for the univariate case.
Later we provide an independent study of the multivariate case.

Define \( \lgd(x;u) = \lgd(x+u) - \lgd(x) - f'(x) u \), 
\begin{EQA}
	\dltw_{3}(x,u)
	& \eqdef &
	\lgd(x+u) - \lgd(x) - f'(x) u - f''(x) \frac{u^{2}}{2} \, .
\label{dl3xufxufxufppxu2}
\end{EQA}

\begin{proposition}
\label{Lintapproxdu}
Let a function \( \lgd(u) \) be two times continuously differentiable and satisfy 
\( f''(x) = - \DU^{2} < 0 \).
If for some \( \zq > 0 \) 
\begin{EQA}
	\sup_{|u| \leq \zq} \frac{1}{\DU^{2} u^{2}/2}\bigl| \dltw_{3}(x,u) \bigr|
	& \leq &
	\dltwb  
	\leq 
	1/3 ,
\label{123stzf3tH2}
\end{EQA}
then for any function \( g(u) \) with \( |g(u)| \leq 1 \) 
\begin{EQA}
	\frac{\bigl| \int_{-\zq}^{\zq}  \ex^{\lgd(x;u)} \, g(u) \, du 
			- \int_{-\zq}^{\zq} \ex^{- \DU^{2} u^{2}/2} \, g(u) \, du \bigr|}
		 {\int \ex^{- \DU^{2} u^{2}/2} \, du}
	& \leq &
	\err_{2}
	\eqdef
	\frac{\dltwb}{2 (1 - \dltwb)} \, .
\label{HimzzeftgttmH3}
\end{EQA}
Moreover, if, in addition to \eqref{123stzf3tH2}, for some \( \fba^{(3)} = \fba^{(3)}(x) \)
\begin{EQA}
	\bigl| \dltw_{3}(x,u) \bigr|
	& \leq &
	\frac{1}{6} \bigl| \fba^{(3)} u^{3} \bigr|,
	\qquad
	|u| \leq \zq,
\label{123stzf3tH3uzq}
\end{EQA}
then
\begin{EQA}
	\frac{\bigl| \int_{-\zq}^{\zq}  \ex^{\lgd(x;u)} \, g(u) \, du 
			- \int_{-\zq}^{\zq} \ex^{- \DU^{2} u^{2}/2} \, g(u) \, du \bigr|}
		 {\int \ex^{- \DU^{2} u^{2}/2} \, du}
	& \leq &
	\err_{3}
	\eqdef
	\frac{0.7 \bigl| \fba^{(3)} \bigr|}{\DU^{3}} \, .
\label{HimzzeftgttmH1}
\end{EQA}
\end{proposition}

\begin{proof}
By \eqref{123stzf3tH2} and \( \dltwb \leq 1/3 \), it holds 
\begin{EQA}
	&& \nquad
	\frac{\DU}{\sqrt{2\pi}}
	\biggl| 
		\int_{-\zq}^{\zq}  \ex^{\lgd(x;u)} \, g(u) \, du - \int_{-\zq}^{\zq} \ex^{- \DU^{2} u^{2}/2} \, g(u) \, du 
	\biggr|
	\\
	& \leq &
	\frac{\DU}{\sqrt{2\pi}}
	\biggl| 
		\int_{-\zq}^{\zq}  \ex^{- (1 - \dltwb) \DU^{2} u^{2}/2} \, du 
		- \int_{-\zq}^{\zq} \ex^{- \DU^{2} u^{2}/2} \, du 
	\biggr|
	\\
	& \leq &
	\frac{1}{\sqrt{2\pi}}
	\biggl| 
		\int  \ex^{- (1 - \dltwb) u^{2}/2} \, du 
		- \int \ex^{- u^{2}/2} \, du 
	\biggr|
	\leq 
	(1 - \dltwb)^{-1/2} - 1 .
\label{Himzzeftgt}
\end{EQA}
We now use that \( (1 - \dltwb)^{-1/2} - 1 \leq  0.5 \dltwb/(1 - \dltwb) \) for \( \dltwb < 1 \),
and the result \eqref{HimzzeftgttmH3} follows.

Now we show \eqref{HimzzeftgttmH1}. It holds
\begin{EQA}
	\int_{-\zq}^{\zq} \ex^{\lgd(x;u)} \, g(u) \, du
	& = &
	\int_{-\zq}^{\zq} \ex^{- \DU^{2} u^{2}/2 + \dltw_{3}(x,u)} \, g(u) \, du .
\label{iUefgu22guu}
\end{EQA}
Define for \( t \geq 0 \)
\begin{EQA}
	\Rem(t)
	&=&
	\int_{-\zq}^{\zq} \ex^{- \DU^{2} u^{2}/2 + t \dltw_{3}(x,u)} \, g(u) \, du .
\label{PhtiUexmHugu}
\end{EQA}
Then for \( t \in [0,1] \) by \eqref{123stzf3tH2} and \eqref{123stzf3tH3uzq}
\begin{EQA}
	\bigl| \Rem'(t) \bigr|
	&=&
	\left| \int_{-\zq}^{\zq} \dltw_{3}(x,u) \ex^{- \DU^{2} u^{2}/2 + t \dltw_{3}(x,u)} \, g(u) \, du \right|
	\leq 
	\int_{-\zq}^{\zq} \bigl| \dltw_{3}(x,u) \bigr| \ex^{- (1 - \dltwb) \DU^{2} u^{2}/2} \, du 
	\\
	& \leq &
	\frac{1}{6} \int_{-\zq}^{\zq} \bigl| \fba^{(3)} u^{3} \bigr| \ex^{- (1 - \dltwb) \DU^{2} u^{2}/2} \, du
	\leq 
	\frac{\bigl| \fba^{(3)} \bigr| \, \DU^{-4}}{6 (1 - \dltwb)^{2}} \int |u|^{3} \ex^{- u^{2}/2} \, du
	= 
	\frac{2\bigl| \fba^{(3)} \bigr| \, \DU^{-4}}{3 (1 - \dltwb)^{2}}
\label{Phit12u2H02u32}
\end{EQA}
and it holds
\begin{EQA}
	\bigl| \Rem(1) - \Rem(0) \bigr|
	& \leq &
	\sup_{t \in [0,1]} \bigl| \Rem'(t) \bigr|
	\leq 
	\frac{2\bigl| \fba^{(3)} \bigr| \, \DU^{-4}}{3 (1 - \dltwb)^{2}} \, 
\label{Pd2d3H0wem22}
\end{EQA}
yielding for \( \dltwb \leq 1/3 \)
\begin{EQA}
	\frac{\bigl| \Rem(1) - \Rem(0) \bigr|}{\int \ex^{- \DU^{2} u^{2}/2} du}
	& \leq &
	\frac{2\bigl| \fba^{(3)} \bigr| \, \DU^{-4}}{3 (1 - \dltwb)^{2} \sqrt{2\pi} \, \DU^{-1}}
	\leq 
	\frac{0.27 \bigl| \fba^{(3)} \bigr|}{(1 - \dltwb)^{2} \DU^{3}} \, 
\label{Phd3Ph0iUmHu22}
\end{EQA}
and \eqref{HimzzeftgttmH1} follows.
\end{proof}

Now we present some results based on a higher order Taylor approximation.
Define 
\begin{EQA}
	\dltw_{4}(x,u)
	& \eqdef &
	\lgd(x+u) - \lgd(x) - f'(x) u + \frac{\DU^{2} u^{2}}{2} - \frac{f^{(3)}(x) u^{3}}{6} 
	\\
	&=&
	\lgd(x,u) + \frac{\DU^{2} u^{2}}{2} - \frac{f^{(3)}(x) u^{3}}{6} \, .
\label{124stzf4tH24}
\end{EQA}
With \( \dltw_{3}(x,u) = \lgd(x,u) + \DU^{2} u^{2}/2 \), it obviously holds
\begin{EQA}
	\dltw_{3}(x,u)
	& = &
	\dltw_{4}(x,u) + {f^{(3)}(x) u^{3}}/{6} .
\label{d3d4uf30u36}
\end{EQA}

\begin{proposition}
\label{Lintapproxdu4}
Let \( \lgd(u) \) be three times continuously differentiable and satisfy \eqref{123stzf3tH2}.
Define \( \dltwm_{3} = f^{(3)}(x) / \DU^{2} \).
If \( \dltw_{4}(x,u) \) from \eqref{124stzf4tH24} satisfies
\begin{EQA}
	\sup_{|u| \leq \zq} \frac{24}{\DU^{2} u^{4}} \bigl| \dltw_{4}(x,u) \bigr|
	& \leq &
	\dltwm_{4}  
\label{124stzf4tH2}
\end{EQA}
then for any function \( g(u) \) with \( |g(u)| \leq 1 \) and \( g(u) = g(-u) \)
\begin{EQA}
	\frac{
	\bigl| \int_{-\zq}^{\zq}  \ex^{\lgd(x;u)} \, g(u) \, du - \int_{-\zq}^{\zq} \ex^{- \DU^{2} u^{2}/2} \, g(u) \, du \bigr|}
	{\int \ex^{- \DU^{2} u^{2}/2} \, du}
	& \leq &
	\err_{4} \, .
\label{HimzzeftgttmH2}
\end{EQA}
with 
\begin{EQA}
	\err_{4}
	& \eqdef &
	\left\{ \frac{5 \dltwm_{3}^{2}}{12 (1 - \dltwb)^{7/2} }
	+ \frac{\dltwm_{4}}{4 (1 - \dltwb)^{5/2}} \right\} \DU^{-2} .
\label{e4l1d372441m3}
\end{EQA}
\end{proposition}

\begin{proof}
We write \( \dltw_{3}(u) \) in place of \( \dltw_{3}(x,u) \) and \( f^{(3)} \) in place of \( f^{(3)}(x) \).
It holds 
\begin{EQA}
	\int_{-\zq}^{\zq} \ex^{\lgd(x;u)} \, g(u) \, du
	& = &
	\int_{-\zq}^{\zq} \exp\Bigl\{ - \frac{\DU^{2} u^{2}}{2} + \dltw_{3}(u) \Bigr\} \, g(u) \, du .
\label{P3tHmzzmt22t3}
\end{EQA}
Define for \( t \in [0,1] \)
\begin{EQA}
	\Rem(t)
	& \eqdef &
	\int_{-\zq}^{\zq} \exp\Bigl\{ - \frac{\DU^{2} u^{2}}{2} + t \dltw_{3}(u) \Bigr\} \, g(u) \, du .
\label{P3tHmzzmt22t}
\end{EQA}
Symmetricity \( g(u) = g(-u) \) implies that 
\begin{EQA}
	\Rem'(0)
	&=&
	\frac{1}{2} \int_{-\zq}^{\zq} 
	\exp\Bigl( - \frac{\DU^{2} u^{2}}{2} \Bigr) \bigl\{ \dltw_{3}(u) + \dltw_{3}(-u) \bigr\} \, g(u) \, du
	=
	\int_{-\zq}^{\zq} \exp\Bigl( - \frac{\DU^{2} u^{2}}{2} \Bigr) \dltw_{4}(u) \, g(u) \, du \, .
\label{Pp012zzmH2u212}
\end{EQA}
Moreover, as  \( |\dltw_{3}(u)| \leq \dltwb \, \DU^{2} u^{2}/2 \),
it holds for \( t \in [0,1] \) 
\begin{EQA}
	|\Rem''(t)|
	&=&
	\left| \int_{-\zq}^{\zq} \dltw_{3}^{2}(u) 
		\exp\Bigl\{ - \frac{\DU^{2} u^{2}}{2} + t \dltw_{3}(u) \Bigr\} \, g(u) \, du 
	\right|
	\leq 
	\int_{-\zq}^{\zq} \dltw_{3}^{2}(u) 
		\exp\Bigl( - \frac{(1 - \dltwb) \DU^{2} u^{2}}{2} \Bigr) \, du \, .
\label{F3pptHm2t6}
\end{EQA}
As \( \dltw_{3}(u) = f^{(3)} u^{3}/6 + \dltw_{4}(u) \) and \( |\dltw_{4}(u)| \leq 1 \), it holds 
\begin{EQA}
	|\Rem''(t)| 
	& \leq & 
	2 \int_{-\zq}^{\zq} \bigl\{ \dltw_{4}^{2}(u) + \bigl| f^{(3)} u^{3}/6 \bigr|^{2} \bigr\} 
		\exp\Bigl( - \frac{(1 - \dltwb)\DU^{2} u^{2}}{2} \Bigr) \, du
	\\
	& \leq & 
	2 \int_{-\zq}^{\zq} \bigl\{ |\dltw_{4}(u)| + \bigl| f^{(3)} \bigr|^{2} u^{6}/36 \bigr\} 
	\exp\Bigl( - \frac{(1 - \dltwb) \DU^{2} u^{2}}{2} \Bigr) \, du .
\label{u221md3u636}
\end{EQA}
The use of \( \bigl| f^{(3)} \bigr| = \dltwm_{3} \, \DU^{2} \) and \eqref{124stzf4tH2} yields
in view of \( \dltwb \leq 1/3 \) and \( |\dltw_{4}(u)| \leq \dltwm_{4} \, \DU^{2} u^{4} \)
\begin{EQA}
	&& \nquad
	\bigl| \Rem(1) - \Rem(0) \bigr|
	\leq 
	\bigl| \Rem'(0) \bigr| + \frac{1}{2} \sup_{t \leq 1} |\Rem''(t)|
	\\
	& \leq &
	2 \int_{-\zq}^{\zq} |\dltw_{4}(u)| \, \exp\Bigl( - \frac{(1 - \dltwb) \DU^{2} u^{2}}{2} \Bigr) \, du
	+
	\frac{\dltwm_{3}^{2} \, \DU^{4}}{36} \int_{-\zq}^{\zq} u^{6} \, 
		\exp\Bigl( - \frac{(1 - \dltwb)\DU^{2} u^{2}}{2} \Bigr) \, du 
	\\
	& \leq &
	\frac{\dltwm_{4} \, \DU^{2}}{12} 
	\int_{-\zq}^{\zq} u^{4} \, \exp\Bigl( - \frac{(1 - \dltwb) \DU^{2} u^{2}}{2} \Bigr) \, du
	+ \frac{\dltwm_{3}^{2} \, \DU^{4}}{36} 
	\int_{-\zq}^{\zq} u^{6} \, \exp\Bigl( - \frac{(1 - \dltwb)\DU^{2} u^{2}}{2} \Bigr) \, du 
	\\
	& \leq &
	\frac{\dltwm_{4} (1 - \dltwb)^{-5/2}}{12 \, \DU^{3}}	
	\int u^{4} \, \exp\Bigl( - \frac{\DU^{2} u^{2}}{2} \Bigr) \, du
	+ \frac{\dltwm_{3}^{2} (1 - \dltwb)^{-7/2}}{36 \, \DU^{3}} 
	\int u^{6} \, \exp\Bigl( - \frac{\DU^{2} u^{2}}{2} \Bigr) \, du 
	\\
	&=&
	\frac{3 \sqrt{2\pi} \, \dltwm_{4} (1 - \dltwb)^{-5/2} }{12 \, \DU^{3}}	
	+ \frac{15 \sqrt{2\pi} \, \dltwm_{3}^{2} (1 - \dltwb)^{-7/2}}{36 \, \DU^{3}}
\label{C4Hm2Fm4tH2t4}
\end{EQA}
yielding \eqref{HimzzeftgttmH2} in view of 
\( \int \ex^{- \DU^{2} u^{2}/2} \, du = \sqrt{2\pi} \, \DU^{-1} \)
and \( \dltwm_{3} \leq 1/3 \).
\end{proof}

\begin{remark}
Suppose that \( \lgd(\cdot) \) is 
three times continuously differentiable and \( |f^{(3)}(x+u)| \leq \dltwu_{3} \DU^{2} \)
for \( |u| \leq \zq \).
Then the expansion \eqref{123stzf3tH2} applies with \( \dltwb \leq \dltwu_{3} \zq/3 \).
The use of \( \zq = \xx/\DU \) yields
\begin{EQA}
	\dltwb \leq \dltwu_{3} \zq/ 3
	& \leq &
	\CONST \, \dltwu_{3} \, \xx \, \DU^{-1} .
\label{Cw5x5Hm3120}
\end{EQA}
Similarly, if \( \lgd(\cdot) \) is four times continuously differentiable and \( |f^{(4)}(u)| \leq \dltwm_{4} \DU^{2} \)
for \( |u| \leq \zq \), then \eqref{124stzf4tH2} holds as well.
One can see that the use of a higher order approximation of \( \lgd \) allows to improve 
the accuracy of approximation  
from \( \DU^{-1} \) as in \eqref{HimzzeftgttmH1} of Proposition~\ref{Lintapproxdu} to \( \DU^{-2} \) as in \eqref{HimzzeftgttmH2} of Proposition~\ref{Lintapproxdu4}.

\end{remark}
}

\Section{Lower and upper Gaussian measures}
This section introduces the lower and upper Gaussian measure which locally sandwich the measure \( \PfL \)
using the decomposition from condition \nameref{LLf0ref}. 
Denote \( - \nabla^{2} \lgd(\xv) = \DV_{\GP}^{2} \).
Definition \eqref{om3esuU1H02d3} enables us to bound with \( \dltwb = \dltwb(\xv) \)
\begin{EQA}
	\frac{1}{2} (\| \DV_{\GP} \uv \|^{2} - \dltwb \| \DVL \uv \|^{2})
	& \leq &
	\lgd(\xv;\uv)
	\leq 
	\frac{1}{2} (\| \DV_{\GP} \uv \|^{2} + \dltwb \| \DVL \uv \|^{2}) \, 
\label{ysdtydsfrtswedftedswfhd}
\end{EQA}
yielding two Gaussian measures which bounds \( \PfL \) locally from above and from below.
The next technical result provides sufficient conditions for their contiguity.

\begin{proposition}
\label{LlocintgrL}
Let \( \dltwb \) from \eqref{om3esuU1H02d3} satisfy \( \dltwb \leq 1/3 \).
Then with \( \dimLG \) from \eqref{dAdetrH02Hm2}
\begin{EQA}
\label{12Ip2t3t4Hp3}
	\det\bigl( \Id + \dltwb \DV_{\GP}^{-1} \DVL^{2} \, \DV_{\GP}^{-1} \bigr)
	& \leq &
	\exp (\dltwb \, \dimLG) \, ,
	\\
	\det\bigl( \Id - \dltwb \DV_{\GP}^{-1} \DVL^{2} \, \DV_{\GP}^{-1} \bigr)^{-1/2}
	& \leq &
	\exp \bigl\{ 3/2 \log (3/2) \, \dltwb \, \dimLG \bigr\} .
	\qquad
\label{12Im2t3t4Hm30}
\end{EQA}
\end{proposition}

\begin{proof}
Without loss of generality assume that \( \DV_{\GP}^{-1} \DVL^{2} \, \DV_{\GP}^{-1} \) is diagonal with 
eigenvalues \( \lambda_{j} \in [0,1] \).
As \( - x^{-1}\log(1 - x) \leq 3 \log(3/2) \) for \( x \in [0,1/3] \), it holds by \eqref{om3esuU1H02d3}
\begin{EQA}
	&& \nquad
	\log \det\bigl( \Id - \dltwb \, \DV_{\GP}^{-1} \DVL^{2} \, \DV_{\GP}^{-1} \bigr)^{-1}
	=
	- \sum_{j=1}^{\dimp} \log\bigl( 1 - \dltwb \lambda_{j} \bigr)
	\leq 
	3 \log (3/2) \sum_{j=1}^{\dimp} \dltwb \lambda_{j}
	\\
	&=&
	3 \log (3/2) \, \dltwb \, \tr\bigl( \DV_{\GP}^{-1} \DVL^{2} \, \DV_{\GP}^{-1} \bigr)
	=
	3 \log (3/2) \, \dltwb \, \dimLG \, 
\label{125d3p0ldm1}
\end{EQA}
yielding \eqref{12Im2t3t4Hm30}.
The proof of \eqref{12Ip2t3t4Hp3} is similar using \( \log(1 + x) \leq x \) for \( x \geq 0 \).
\end{proof}

\Section{Local approximation}
\label{SLaplapprmu}

This section presents the bounds on the error \( \err \) of local approximation \eqref{errdefdiUaHu}.
The first result only uses \( \dltwb \, \dimLG < 1 \) while the second one also assumes \nameref{LL3fref}.
These two results allow to bound the total variation distance between \( \PfL \) and \( \ND(\xv,\DV_{\GP}^{-2}) \).
The third result is finer and is based on \nameref{LL4fref}.
It allows to improve the error term of Gaussian approximation over the class of centrally symmetric sets.
We also present some extensions for the moments of \( \PfL \).

\begin{proposition}
\label{Lintfxupp3}
Let \( \dltwb = \dltwb(\xv) \) from \eqref{om3esuU1H02d3} and \( \dimLG \) from \eqref{dAdetrH02Hm2} satisfy 
\begin{EQA}
	\dltwb \, \dimLG
	& \leq &
	2/3 \, .
\label{w3p0le13fx}
\end{EQA}
Then for any function \( g(\uv) \) with \( |g(\uv)| \leq 1 \) 
\begin{EQA}
	\biggl| 
	\frac{\int_{\UVL} \ex^{\lgd(\xv;\uv)} \, g(\uv) \, d\uv - \int_{\UVL} \ex^{- \| \DV_{\GP} \uv \|^{2}/2} \, g(\uv) \, d\uv} 
		 {\int \ex^{- \| \DV_{\GP} \uv \|^{2}/2} d\uv} 
	\biggr|
	& \leq &
	\err
	=
	\err_{2} 
	= 
	\frac{0.75 \, \dltwb \, \dimLG}{1 - \dltwb} \, .
\label{ed3le2d3pGd}
\end{EQA}
Moreover, under \nameref{LL3fref}, the bound applies with \( \err = \err_{3} \); see \eqref{5qw7dyf4e4354coefw9dufih}.
\end{proposition}

\begin{proof}
Under \eqref{w3p0le13fx}, bound \eqref{12Im2t3t4Hm30} implies 
\begin{EQA}
	\det\bigl( \Id - \dltwb \DV_{\GP}^{-1} \DVL^{2} \, \DV_{\GP}^{-1} \bigr)^{-1/2}
	& \leq &
	\exp \bigl\{ 3/2 \log (3/2) \, \dltwb \, \dimLG \bigr\} 
	\leq 
	3/2 \, .
	\qquad
\label{12Im2t3t4Hm3}
\end{EQA}
Define for \( t \geq 0 \)
\begin{EQA}
	\Rem(t)
	&=&
	\int_{\UVL} \ex^{- \| \DV_{\GP} \uv \|^{2}/2 + t \dltw_{3}(\xv,\uv)} \, g(\uv) \, d\uv .
\label{PhtiUexmHu2mtH0u2g}
\end{EQA}
Then for \( t \in [0,1] \) by \eqref{om3esuU1H02d3} 
\begin{EQA}
	\bigl| \Rem'(t) \bigr|
	&=&
	\left| \int_{\UVL} \dltw_{3}(\xv,\uv) \ex^{- \| \DV_{\GP} \uv \|^{2}/2 + t \dltw_{3}(\xv,\uv)} \, g(\uv) \, d\uv \right|
	\\
	& \leq &
	\int_{\UVL} \bigl| \dltw_{3}(\xv,\uv) \bigr| \ex^{- (\| \DV_{\GP} \uv \|^{2} - \dltwb \| \DVL \uv \|^{2})/2} \, d\uv .
\label{Pht12u212H02u22}
\end{EQA}
Now we make change of variable \( \Idd \uv \) to \( \wv \) with
\( \Idd^{2} = \Id - \dltwb \DV_{\GP}^{-1} \DVL^{2} \, \DV_{\GP}^{-1} \).
By \eqref{12Im2t3t4Hm3} \( \det \Idd^{-1} \leq 3/2 \) and also 
\( \| \Idd^{-1} \| \leq (1 - \dltwb)^{-1/2} \).
By \eqref{om3esuU1H02d3}
\begin{EQA}
	&& \nquad
	\bigl| \Rem(1) - \Rem(0) \bigr|
	\leq 
	\sup_{t \in [0,1]} \bigl| \Rem'(t) \bigr|
	\leq 
	\frac{\dltwb}{2}
	\int_{\UVL} \| \DVL \uv \|^{2} \ex^{- (\| \DV_{\GP} \uv \|^{2} - \dltwb \| \DVL \uv \|^{2})/2} \, d\uv 
	\\
	& \leq &
	\frac{3 \dltwb}{4}
	\int \| \DVL \Idd^{-1} \wv \|^{2} \ex^{- \| \DV_{\GP} \wv \|^{2} /2} \, d\wv 
	\leq 
	\frac{3 \dltwb}{4(1 - \dltwb)}
	\int \| \DVL \wv \|^{2} \ex^{- \| \DV_{\GP} \wv \|^{2} /2} \, d\wv 	.
\label{ppt2iUH0w2edw}
\end{EQA}
In view of \( \E \| \DVL \, \gaussv_{\GP} \|^{2} = \tr\bigl( \DVL^{2} \, \DV_{\GP}^{-2} \bigr) \) for a standard normal
\( \gaussv \), we derive 
\begin{EQA}
	\frac{\bigl| \Rem(1) - \Rem(0) \bigr|}{\int \ex^{- \| \DV_{\GP} \uv \|^{2}/2} d\uv}
	& \leq &
	\frac{3 \dltwb }{4 (1 - \dltwb)}
	\frac{\int \| \DVL \wv \|^{2} \ex^{- \| \DV_{\GP} \wv \|^{2} /2} \, d\wv}
		{\int \ex^{- \| \DV_{\GP} \wv \|^{2} /2} \, d\wv}
	\leq 
	\frac{3 \dltwb \, \dimLG}{4 (1 - \dltwb)}
\label{Phd3Ph0iUmHu22}
\end{EQA}
and \eqref{ed3le2d3pGd} follows.
Under \nameref{LL3fref}
\begin{EQA}
	\bigl| \Rem(1) - \Rem(0) \bigr|
	& \leq &
	\frac{\det (\Idd^{-1})}{6} 
	\int \dltwu_{3}(\Idd^{-1} \wv) \, \ex^{- \| \DV_{\GP} \wv \|^{2} /2} \, d\wv 
	\\
	& \leq &
	\frac{1}{4(1 - \dltwb)^{3/2}}
	\int \dltwu_{3}(\uv) \, \ex^{- \| \DV_{\GP} \uv \|^{2} /2} \, d\uv 
\label{Pd3p2d3H0wem22g1}
\end{EQA}
yielding the second statement.
\end{proof}

The result can be extended to the case with a homogeneous function \( g(\uv) \),
e.g. \( g(u) = \| \QP \uv \|^{m} \).

\begin{proposition}
\label{Pd3t3hot3d3u}
Suppose the conditions of Proposition~\ref{Lintfxupp3} and \nameref{LL3fref}.
Then for \( m \geq 1 \) and any \( m \)-homogeneous function \( g(\cdot) \) with \( g(t\uv) = t^{m} g(\uv) \)
\begin{EQA}
	\biggl| 
	\frac{\int_{\UVL} \ex^{\lgd(\xv;\uv)} \, g(\uv) \, d\uv 
	- \int_{\UVL} \ex^{- \| \DV_{\GP} \uv \|^{2}/2} \, g(\uv) \, d\uv} 
		 {\int \ex^{- \| \DV_{\GP} \uv \|^{2}/2} d\uv} 
	\biggr|
	& \leq &
	\frac{\E \bigl\{ |g(\gaussv_{\GP})| \, \dltwu_{3}(\gaussv_{\GP}) \bigr\}}{4 (1 - \dltwb)^{(m+3)/2}} \, ,
	\qquad
\label{e3dEf3xGHo3}
\end{EQA}
where \( \gaussv_{\GP} \sim \ND(0,\DV_{\GP}^{-2}) \) is a Gaussian element in \( \R^{\dimp} \).
\end{proposition}

\begin{proof}
Similarly to the proof of Proposition~\ref{Lintfxupp3}, under \eqref{bd3xu16f3uo3} for \( \Rem(t) \) from
\eqref{PhtiUexmHu2mtH0u2g}
\begin{EQA}
	\bigl| \Rem(1) - \Rem(0) \bigr|
	& \leq &
	\frac{\det (\Idd^{-1})}{6} 
	\int \dltwu_{3}(\Idd^{-1} \wv) \, \bigl| g(\Idd^{-1} \wv) \bigr| \, \ex^{- \| \DV_{\GP} \wv \|^{2} /2} \, d\wv 
	\\
	& \leq &
	\frac{1}{4(1 - \dltwb)^{(m+3)/2}}
	\int \dltwu_{3}(\wv) \, \bigl| g(\wv) \bigr| \, \ex^{- \| \DV_{\GP} \wv \|^{2} /2} \, d\wv 
\label{Pd3p2d3H0wem22}
\end{EQA}
yielding 
\begin{EQA}
	\frac{\bigl| \Rem(1) - \Rem(0) \bigr|}{\int \ex^{- \| \DV_{\GP} \uv \|^{2}/2} d\uv}
	& \leq &
	\frac{1}{4(1 - \dltwb)^{(m+3)/2}}
	\frac{\int 
			\dltwu_{3}(\wv) \, |g(\wv)| \, \ex^{- \| \DV_{\GP} \wv \|^{2} /2} \, d\wv}
		{\int \ex^{- \| \DV_{\GP} \wv \|^{2} /2} \, d\wv} \, .
\label{Phd3Ph0iUmHu22}
\end{EQA}
This yields \eqref{e3dEf3xGHo3}.
\end{proof}

Important special cases correspond to \( m=1 \).

\begin{proposition}
\label{Perrbound3}
Suppose the conditions of Proposition~\ref{Pd3t3hot3d3u}.
Then with \( \gaussv_{\GP} \sim \ND(0,\DV_{\GP}^{-2}) \), it holds 
for any linear mapping \( \QP \colon \R^{\dimp} \to \R^{\dimq} \) and any vector \( \av \in \R^{\dimq} \)
\begin{EQA}
	\frac{\bigl| 
		\int_{\UVL} \ex^{\lgd(\xv;\uv)} \, \langle \QP \uv,\av \rangle \, d\uv 
		- \int_{\UVL} \ex^{- \| \DV_{\GP} \uv \|^{2}/2} \, \langle \QP \uv,\av \rangle \, d\uv 	 
		\bigr|} 
		 {\int \ex^{- \| \DV_{\GP} \uv \|^{2}/2} d\uv} 
	& \leq &
	\frac{\E \bigl\{ \bigl| \langle \QP \gaussv_{\GP},\av \rangle \bigr| \, \dltwu_{3}(\gaussv_{\GP}) \bigr\}}{4 (1 - \dltwb)^{2}} \, .
	\qquad
\label{iUafxvdvH2vk3}
\end{EQA}
\end{proposition}

Now we state a sharper result based on \nameref{LL4fref}.

\begin{proposition}
\label{Lintfxupp2}
Suppose the conditions of Proposition~\ref{Lintfxupp3} and \nameref{LL4fref}.
Then for any function \( g(\uv) \) with \( |g(\uv)| \leq 1 \) and 
\( g(\uv) = g(-\uv) \)
\begin{EQA}
	\biggl| 
	\frac{\int_{\UVL} \ex^{\lgd(\xv;\uv)} \, g(\uv) \, d\uv - \int_{\UVL} \ex^{- \| \DV_{\GP} \uv \|^{2}/2} \, g(\uv) \, d\uv} 
		 {\int \ex^{- \| \DV_{\GP} \uv \|^{2}/2} d\uv} 
	\biggr|
	& \leq &
	\err_{4} \, ,
\label{efuiguduUem}
\end{EQA}
where for a Gaussian element \( \gaussv_{\GP} \sim \ND(0,\DV_{\GP}^{-2}) \) in \( \R^{\dimp} \)
\begin{EQA}
	\err_{4}
	& \eqdef &
	\frac{1}{16 (1 - \dltwb)^{2}} \Bigl\{ \E \bigl\langle \nabla^{3} \lgd(\xv) , \gaussv_{\GP}^{\otimes 3} \bigr\rangle^{2} 
	+ 2 \E \dltwu_{4}(\gaussv_{\GP}) \Bigr\}
	\, .
\label{errdef3322Hm2}
\end{EQA}
If the function \( g(\cdot) \) is not bounded by one but it is symmetric and \( 2m \)-homogeneous,
i.e. \( g(t \uv) = t^{2m} g(\uv) \), then
\eqref{efuiguduUem} still applies with 
\begin{EQA}
	\err_{4}
	& \eqdef &
	\frac{1}{16 (1 - \dltwb)^{2+m}} \E \Bigl\{ 
		\bigl\langle \nabla^{3} \lgd(\xv) , \gaussv_{\GP}^{\otimes 3} \bigr\rangle^{2} \, g(\gaussv_{\GP})
	+ 2 \dltwu_{4}(\gaussv_{\GP}) \, g(\gaussv_{\GP}) \Bigr\}
	\, .
\label{errdef3322Hm2m}
\end{EQA}
\end{proposition}

\begin{proof}
Below we write \( f^{(3)} \) instead of \( \nabla^{3} \lgd(\xv) \) and \( \dltw_{k}(\uv) \) in place of 
\( \dltw_{k}(\xv,\uv) \), \( k=3,4 \).
It holds 
\begin{EQA}
	\int_{\UVL} \ex^{\lgd(\xv;\uv)} \, g(\uv) \, d\uv
	& = &
	\int_{\UVL} \exp\Bigl\{ - \frac{\| \DV_{\GP} \uv \|^{2}}{2} + \dltw_{3}(\uv) \Bigr\} \, g(\uv) \, d\uv .
\label{P3tHmzzmt22t3v}
\end{EQA}
Define for \( t \in [0,1] \)
\begin{EQA}
	\Rem(t)
	& \eqdef &
	\int_{\UVL} \exp\Bigl\{ - \frac{\| \DV_{\GP} \uv \|^{2}}{2} + t \dltw_{3}(\uv) \Bigr\} \, g(\uv) \, d\uv .
\label{P3tHmzzmt22t}
\end{EQA}
Symmetricity of \( \UVL \) and \( g(\uv) = g(-\uv) \) implies that 
\begin{EQA}
	\Rem'(0)
	&=&
	\frac{1}{2} \int_{\UVL} 
	\exp\Bigl( - \frac{\| \DV_{\GP} \uv \|^{2}}{2} \Bigr) \bigl\{ \dltw_{3}(\uv) + \dltw_{3}(-\uv) \bigr\} \, g(\uv) \, d\uv
	\\
	&=&
	\int_{\UVL} \exp\Bigl( - \frac{\| \DV_{\GP} \uv \|^{2}}{2} \Bigr) \bar{\dltw}_{4}(\uv) \, g(\uv) \, d\uv \, 
\label{Pp012zzmH2u212v}
\end{EQA}
with \( \bar{\dltw}_{4}(\uv) = \bigl\{ \dltw_{4}(\uv) + \dltw_{4}(-\uv) \bigr\}/2 \).
Moreover, as  \( |\dltw_{3}(\uv)| \leq \dltwb \| \DVL \uv \|^{2}/2 \),
it holds for \( t \in [0,1] \) 
\begin{EQA}
	|\Rem''(t)|
	& \leq &
	\int_{\UVL} \dltw_{3}^{2}(\uv) 
		\exp\Bigl\{ - \frac{\| \DV_{\GP} \uv \|^{2}}{2} + t \dltw_{3}(\uv) \Bigr\} \, |g(\uv)| \, d\uv
	\\
	& \leq &
	\int_{\UVL} \dltw_{3}^{2}(\uv) 
		\exp\Bigl( - \frac{\| \DV_{\GP} \uv \|^{2} - \dltwb \| \DVL \uv \|^{2}}{2} \Bigr) \, d\uv \, .
\label{F3pptHm2t6v}
\end{EQA}
As \( \dltw_{3}(\uv) = \bigl\langle f^{(3)}, \uv^{\otimes 3} \bigr\rangle/6 + \dltw_{4}(\uv) \) and \( |\dltw_{4}(\uv)| \leq 1 \), one can bound for \( t \in [0,1] \)
\begin{EQA}
	|\Rem''(t)| 
	& \leq & 
	2 \int_{\UVL} \bigl\{ \bar{\dltw}_{4}^{2}(\uv) + \bigl| \bigl\langle f^{(3)}, \uv^{\otimes 3} \bigr\rangle/6 \bigr|^{2} \bigr\} 
		\exp\Bigl( - \frac{\| \DV_{\GP} \uv \|^{2} - \dltwb \| \DVL \uv \|^{2}}{2} \Bigr) \, d\uv
	\\
	& \leq & 
	2 \int_{\UVL} \bigl\{ |\bar{\dltw}_{4}(\uv)| + \bigl\langle f^{(3)}, \uv^{\otimes 3} \bigr\rangle^{2} /36 \bigr\} 
	\exp\Bigl( - \frac{\| \DV_{\GP} \uv \|^{2} - \dltwb \| \DVL \uv \|^{2}}{2} \Bigr) \, d\uv .
\label{u221md3u636v}
\end{EQA}
This and \eqref{Pp012zzmH2u212v} yield
\begin{EQA}
	&& \nquad
	\bigl| \Rem(1) - \Rem(0) \bigr|
	\leq 
	\bigl| \Rem'(0) \bigr| + \frac{1}{2} \sup_{t \in [0,1]} |\Rem''(t)|
	\leq 
	2 \int_{\UVL} |\bar{\dltw}_{4}(\uv)| \, \exp\Bigl( - \frac{\| \DV_{\GP} \uv \|^{2} - \dltwb \| \DVL \uv \|^{2}}{2} \Bigr) \, d\uv
	\\
	&& \quad
	+ \, \frac{1}{36} \int_{\UVL} \bigl\langle f^{(3)}, \uv^{\otimes 3} \bigr\rangle^{2} \, 
		\exp\Bigl( - \frac{\| \DV_{\GP} \uv \|^{2} - \dltwb \| \DVL \uv \|^{2}}{2} \Bigr) \, d\uv .
\label{C4Hm2Fm4tH2tv}
\end{EQA}
Change of variable \( \bigl( \Id - \dltwb \, \DV_{\GP}^{-1} \DVL^{2} \, \DV_{\GP}^{-1} \bigr)^{1/2} \uv \) to \( \wv \) yields by 
\eqref{12Im2t3t4Hm3} in view of \( \dltwb \leq 1/3 \)
\begin{EQA}
	&& \nquad
	\frac{1}{36} \int_{\UVL} \bigl\langle f^{(3)}, \uv^{\otimes 3} \bigr\rangle^{2} \, 
		\exp\Bigl( - \frac{\| \DV_{\GP} \uv \|^{2} - \dltwb \| \DVL \uv \|^{2}}{2} \Bigr) \, d\uv 
	\\
	& \leq &
	\frac{3/2}{36 (1 - \dltwb)^{3}} 
	\int \bigl\langle f^{(3)}, \wv^{\otimes 3} \bigr\rangle^{2} \, \exp\Bigl( - \frac{\| \DV_{\GP} \wv \|^{2}}{2} \Bigr) \, d\wv .
\label{C4Hm2Fm4tH2t3v}
\end{EQA}
Similarly by \eqref{1mffmxum5}
\begin{EQA}
	&& \nquad
	\int_{\UVL} |\bar{\dltw}_{4}(\uv)| \,\exp\Bigl( - \frac{\| \DV_{\GP} \uv \|^{2} - \dltwb \| \DVL \uv \|^{2}}{2} \Bigr) \, d\uv
	\\
	& \leq &
	\frac{3/2}{24(1 - \dltwb)^{2}} 
	\int  
	\dltwu_{4}(\wv) \, \exp\Bigl( - \frac{\| \DV_{\GP} \wv \|^{2}}{2} \Bigr) \, d\wv .
\label{C4Hm2Fm4tH2t4v}
\end{EQA}
The use of \( \dltwb \leq 1/3 \) implies that
\begin{EQA}
	\frac{\bigl| \Rem(1) - \Rem(0) \bigr|}{\int_{\UVL} \ex^{- \| \DV_{\GP} \uv \|^{2}/2} d\uv}
	& \leq &
	\frac{3/2}{24 (1 - \dltwb)^{2}}
	\Bigl\{  
	\E \bigl\langle f^{(3)} , \gaussv_{\GP}^{\otimes 3} \bigr\rangle^{2}
		+ 2 \E \dltwu_{4}(\gaussv_{\GP})
	\Bigr\} 
	\leq 
	\err_{4} \, 
\label{dwHw22n3n432}
\end{EQA}
and \eqref{efuiguduUem} follows.
The proof of \eqref{errdef3322Hm2m} is similar.
\end{proof}

\Section{Tail integrals}
\label{StailLa}
In this section we also write \( \xv \) in place of \( \xvs \).
Below we evaluate \( \rho \) from \eqref{rhfiUaceHu22m} which bounds the integral of \( \ex^{\lgd(\xv;\uv)} \) 
over the complement of the local set \( \UVL \) of a special form 
\( \UVL = \bigl\{ \uv \colon \| \DVL \uv \| \leq \amax^{-1} \rrL \bigr\} \) for \( \DVL \) from \nameref{LLf0ref}.
Our results help to understand how the radius \( \rrL \) should be fixed to ensure \( \rho \) sufficiently small.
The main tools for the analysis are deviation probability bounds for Gaussian quadratic forms; see Section~\ref{Sdevboundgen}.

\begin{proposition}
\label{PUVstarLLfx}
Suppose \nameref{LLf0ref}.
Given \( \amax < 1 \) and \( \xx > 0 \), 
let \( \UVL \) and \( \rrL \) be defined by \eqref{UvTDunm12spT}.
Let also \( \dltwb \) from \eqref{om3esuU1H02d3} satisfy
\begin{EQA}
	\dltwb
	& \leq &
	1 - \amax. 
\label{1w3rT2spTs2xfx}
\end{EQA}
Then 
\begin{EQA}
\label{fiinIHUniUef0}
	\frac{\int \Ind\bigl( \uv \not\in \UVL \bigr) \, \ex^{\lgd(\xv;\uv)} \, d\uv}
		 {\int \ex^{ - \| \DV_{\GP} \uv \|^{2}/2} \, d \uv}
	& \leq &
	4 \ex^{-\xx - \dimLG/2} \, ,
	\\
	\frac{\int \Ind\bigl( \uv \not\in \UVL \bigr) \, \ex^{ - \| \DV_{\GP} \uv \|^{2}/2} \, d\uv}
		 {\int \ex^{ - \| \DV_{\GP} \uv \|^{2}/2} \, d \uv}
	& \leq &
	\ex^{-\xx - \dimLG/2} \, .
\label{fiinIHUniUef}
\end{EQA}
\end{proposition}

\begin{proof}
Let \( \uv \not \in \UVL \), i.e. \( \| \DVL \uv \| > \rr  \) with \( \rr = \amax^{-1} \rrL \).
Define \( \uvc = \rr \| \DVL \uv \|^{-1} \uv \) yielding \( \| \DVL \uvc \| = \rr \).
We also write \( \uv = (1 + \tau) \uvc \) for \( \tau > 0 \).
By \eqref{om3esuU1H02d3} 
and \( \nabla^{2} \lgdL(0) = - \DVL^{2} \)
\begin{EQ}[rcl]
	\lgdL(\uvc) - \lgdL(0) - \bigl\langle \nabla \lgdL(0), \uvc \bigr\rangle 
	& \leq &
	- (1 - \dltwb) \| \DVL \uvc \|^{2}/2 ,
	\\
	\bigl\langle \nabla \lgdL(\uvc) - \nabla \lgdL(0), \uv - \uvc \bigr\rangle
	& \leq &
	- (1 - \dltwb) \bigl\langle \DVL^{2} \uvc, \uv - \uvc \bigr\rangle .
\label{2ud1mud3Hv1dfx}
\end{EQ}
Concavity of \( \lgdL(\uv) \) implies for \( \uv = (1 + \tau) \uvc \),
\begin{EQA}
	\lgdL(\uv) 
	& \leq &
	\lgdL(\uvc) + \bigl\langle \nabla \lgdL(\uvc), \uv - \uvc \bigr\rangle 
\label{fudfpudumudwuld}
\end{EQA}
yielding by \eqref{2ud1mud3Hv1dfx} in view of 
\( \bigl\langle \DVL \uvc, \DVL \uv \bigr\rangle = \| \DVL \uvc \| \, \| \DVL \uv \| \)
\begin{EQA}
	&& \nquad
	\lgdL(\uv) - \lgdL(0) - \bigl\langle \nabla \lgdL(0), \uv \bigr\rangle
	\\
	& = &
	\lgdL(\uv) - \lgdL(\uvc) - \bigl\langle \nabla \lgdL(\uvc), \uv - \uvc \bigr\rangle
	\\
	&&
	+ \, \lgdL(\uvc) - \lgdL(0) - \bigl\langle \nabla \lgdL(0), \uvc \bigr\rangle
	+ \bigl\langle \nabla \lgdL(\uvc) - \nabla \lgdL(0), \uv - \uvc \bigr\rangle
	\\
	& \leq &
	(1 - \dltwb) \| \DVL \uvc \|^{2}/2
	- (1 - \dltwb) \bigl\langle \DVL \uvc, \DVL \uv \bigr\rangle 
	\leq 
	- (1 - \dltwb) \| \DVL \uvc \| \, \| \DVL \uv \| /2.
\label{tnfudnfHm1Hm1d}
\end{EQA}
We now use that 
\( \| \DVL \uvc \| = \rr \),  
\( \uvc = \uv/(1 + \tau) \), and thus,
\begin{EQA}
	&& \nquad
	\lgd(\xv + \uv) - \lgd(0) - \bigl\langle \nabla \lgd(\xv), \uv \bigr\rangle
	\\
	&=&
	\lgdL(\uv) - \lgdL(0) - \bigl\langle \nabla \lgdL(0), \uv \bigr\rangle
	- \| \DV_{\GP} \uv \|^{2}/2 + \| \DVL \uv \|^{2}/2 
	\\
	& \leq &
	- (1 - \dltwb) \rr \| \DVL \uv \|/2 
	- \| \DV_{\GP} \uv \|^{2}/2 + \| \DVL \uv \|^{2}/2.
\label{DGa22DTa221mw3rT22}
\end{EQA}
This yields by \( \rrL = \amax \, \rr \leq (1 - \dltwb) \rr \) with \( \Tau = \DVL \DV_{\GP}^{-1} \)
\begin{EQA}
	&& \nquad
	\frac{\int \Ind\bigl( \uv \not\in \UVL \bigr) 
		\exp\bigl\{  \lgd(\xv + \uv) - \lgd(\xv) - \langle \nabla \lgd(\xv), \uv \rangle \bigr\} \, d\uv}
		{\int \exp \bigl( - \| \DV_{\GP} \uv \|^{2}/{2} \bigr) \, d \uv}
	\\
	& \leq &
	\frac{\int \Ind\bigl( \| \DVL \uv \| > \rr \bigr)
		\exp \bigl\{ - {(1 - \dltwb) \rr} \| \DVL \uv \|/2 - \| \DV_{\GP} \uv \|^{2}/2 + \| \DVL \uv \|^{2}/2 \bigr\} \, d\uv}
		{\int \exp \bigl( - \| \DV_{\GP} \uv \|^{2}/{2} \bigr) \, d \uv}
	\\
	& \leq &
	\E \exp \Bigl\{ 
		- \frac{\rrL}{2} \| \Tau \gaussv \| 
		+ \frac{1}{2} \| \Tau \gaussv \|^{2}  
		\Bigr\} \Ind\bigl( \| \Tau \gaussv \| > \rrL \bigr) 
\label{2C2emx1212IganiU}
\end{EQA}
with \( \gaussv \) standard normal in \( \R^{\dimp} \).
Next, define
\begin{EQA}
	\riskt_{0}(\rr)
	& \eqdef &
	\E \exp\Bigl( - \frac{\rr}{2} \| \DVL \gaussv \| + \frac{1}{2} \| \DVL \gaussv \|^{2} \Bigr) 
		\Ind(\| \DVL \gaussv \| > \rr) .
\label{mp2mxr2g22T}
\end{EQA}
Integration by parts allows to represent the last integral as
\begin{EQA}
	\riskt_{0}(\rr)
	&=&
	- \int_{\rr}^{\infty} \exp\Bigl( - \frac{ \rr \, \zq}{2} + \frac{\zq^{2}}{2} \Bigr) \, 
		d\P\bigl( \| \Tau \gaussv \| > \zq \bigr)
	\\
	&=&
	\P\bigl( \| \Tau \gaussv \| > \rr \bigr) 
	+ \int_{\rr}^{\infty} (\zq - \rr/2) \exp\Bigl( - \frac{\rr \zq}{2} + \frac{\zq^{2}}{2} \Bigr) \, 
		\P\bigl( \| \Tau \gaussv \| > \zq \bigr) \, d\zq \, .
\label{dzzTPzr2rT}
\end{EQA}
By Theorem~\ref{TexpbLGA}, for any \( \zq \geq \sqrt{\dimLG} \) for \( \dimLG = \tr(\Tau \, \Tau^{\T}) = \tr (\DVL^{2} \, \DV_{\GP}^{-2}) \)
\begin{EQA}
	\P\bigl( \| \Tau \gaussv \| > \zq \bigr)
	& \leq & 
	\exp\bigl\{ - (\zq - \sqrt{\dimLG})^{2}/2 \bigr\}
\label{2emxPTgasps2d}
\end{EQA}
yielding for \( \zq \geq \rrL = 2 \sqrt{\dimLG} + \sqrt{2\xx} \)
\begin{EQA}
	\P\bigl( \| \Tau \gaussv \| > \zq \bigr)
	& \leq &
	\exp\bigl\{ - (\zq - \sqrt{\dimLG})^{2}/2 \bigr\}
	\leq 
	\ex^{-\xx - \dimLG/2}
\label{2exmp22ezdp2}
\end{EQA}
and for \( \rr \geq 2 \sqrt{\dimLG} + \sqrt{2\xx} \) and \( \xx \geq 2 \)
\begin{EQA}
	\riskt_{0}(\rr)
	& \leq &
	\ex^{-\xx - \dimLG/2}
	+ \int_{\rr}^{\infty} (\zq - \rr/2) 
	\exp\Bigl\{ - \frac{\rr \zq}{2} + \frac{\zq^{2}}{2} - \frac{(\zq - \sqrt{\dimLG})^{2}}{2} \Bigr\} \, d\zq
	\\
	& \leq &
	\ex^{-\xx - \dimLG/2}
	+ \exp\Bigl( - \frac{(\rr - \sqrt{\dimLG})^{2}}{2} \Bigr) \int_{0}^{\infty} \Bigl(\zq + \frac{\rr}{2} \Bigr) 
	\exp\Bigl\{ - \frac{(\rr - 2 \sqrt{\dimLG}) \zq}{2} \Bigr\} \, d\zq
	\\
	& \leq &
	2 \ex^{-\xx - \dimLG/2} .
\label{Cexmxmp222Tr22fx}
\end{EQA}
This completes the proof of the result \eqref{fiinIHUniUef0}.
The second statement \eqref{fiinIHUniUef} is about Gaussian probability 
\( \P\bigl( \| \Tau \gaussv \| \geq \rrL \bigr) \) for a standard normal element \( \gaussv \), and we derive
\begin{EQA}
	\P\bigl( \| \Tau \gaussv \| \geq 2 \sqrt{\dimLG} + \sqrt{2\xx} \bigr)
	& \leq &
	\exp\bigl\{ - \bigl( \sqrt{\dimLG} + \sqrt{2\xx} \bigr)^{2}/2 \bigr\}
	\leq 
	\exp\bigl( -\xx - \dimLG/2 \bigr) 
\label{fiinIHUniUefi}
\end{EQA}
and \eqref{fiinIHUniUef} follows.
\end{proof}

\medskip

The next result extends \eqref{fiinIHUniUef0}.
\begin{proposition}
\label{PUVstarLLgenlifx}
Assume the conditions of Proposition~\ref{PUVstarLLfx}
with
\begin{EQA}
	\rrL
	& \geq &
	2 \sqrt{\dimLG} + \sqrt{2 \xx} + m
\label{1w3rT2spTs2xlifx}
\end{EQA}
for some \( m \geq 0 \).
Then 
\eqref{fiinIHUniUef0} can be extended to
\begin{EQA}
\label{fiinIHUniUef0m}
	\frac{\int \Ind\bigl( \uv \not\in \UVL \bigr) \, \| \DVL \uv \|^{m} \, \ex^{\lgd(\xv;\uv)} \, d\uv}
		 {\int \ex^{ - \| \DV_{\GP} \uv \|^{2}/2} \, d \uv}
	& \leq &
	4 \ex^{-\xx - \dimLG/2} \, ,
	\\
	\frac{\int \Ind\bigl( \uv \not\in \UVL \bigr) \, \| \DVL \uv \|^{m} \, \ex^{ - \| \DV_{\GP} \uv \|^{2}/2} \, d\uv}
		 {\int \ex^{ - \| \DV_{\GP} \uv \|^{2}/2} \, d \uv}
	& \leq &
	\ex^{-\xx - \dimLG/2} \, .
\label{fiinIHUniUefm}
\end{EQA}
\end{proposition}

\begin{proof}
The case \( m > 0 \) can be proved similarly to \( m=0 \) using 
\( m \log z \leq m z \).
\end{proof}

\Section{Local concentration}
Here we show that the measure \( \PfL \) well concentrates on the local set \( \UVL \) from \eqref{UvTDunm12spT}.
Again we fix \( \xv = \xvs \).

\begin{proposition}
\label{PlocconLa}
Assume \( \dltwb \leq 1/3 \).
Then
\begin{EQA}
	\int_{\UVL} \ex^{\lgd(\xv;\uv)} \, d\uv  
	& \geq &
	\ex^{ - \dltwb \, \dimLG/2} \, \int_{\UVL} \ex^{- \| \DV_{\GP} \uv \|^{2}/2} d\uv \, .
\label{emdw3dp0iUfxu}
\end{EQA}
Moreover,
\begin{EQA}
	\frac{\int_{\UVL^{c}} \ex^{\lgd(\xv;\uv)} \, d\uv}{\int \ex^{\lgd(\xv;\uv)} \, d\uv}  
	& \leq &
	4 \ex^{ - \xx - (1 - \dltwb) \, \dimLG/2} 
	\leq 
	\ex^{-\xx} \, .
\label{emdw3dp0iUfxuL}
\end{EQA}
\end{proposition}

\begin{proof}
By \eqref{om3esuU1H02d3}
\begin{EQA}
	\int_{\UVL} \ex^{\lgd(\xv;\uv)} \, d\uv
	& = &
	\int_{\UVL} \ex^{- \| \DV_{\GP} \uv \|^{2}/2 + \dltw_{3}(\xv,\uv)} \, d\uv 
	\geq 
	\int_{\UVL} \ex^{- \| \DV_{\GP} \uv \|^{2}/2 - \dltwb \| \DVL \uv \|^{2}/2} \, d\uv \, .
\label{PhtiUexmHu222}
\end{EQA}
Change of variable \( \bigl( \Id + \dltwb \, \DV_{\GP}^{-1} \DVL^{2} \, \DV_{\GP}^{-1} \bigr)^{1/2} \uv \) to \( \wv \) yields by \eqref{12Ip2t3t4Hp3} 
\begin{EQA}
	\int_{\UVL} \ex^{\lgd(\xv;\uv)} \, d\uv
	& \geq &
	\det \bigl( \Id + \dltwb \, \DV_{\GP}^{-1} \DVL^{2} \, \DV_{\GP}^{-1} \bigr)^{-1/2} 
	\int_{\UVL} \ex^{- \| \DV_{\GP} \wv \|^{2}/2 } \, d\wv
	\\
	& \geq &
	\ex^{- \dltwb \, \dimLG / 2} \int_{\UVL} \ex^{- \| \DV_{\GP} \wv \|^{2}/2 } \, d\wv ,
\label{iUefxudetId3H121}
\end{EQA}
and \eqref{emdw3dp0iUfxu} follows.
This and \eqref{fiinIHUniUef0}, \eqref{fiinIHUniUef} of Proposition~\ref{PUVstarLLfx} imply
\begin{EQA}
	&& \nquad
	\frac{\int_{\UVL^{c}} \ex^{\lgd(\xv;\uv)} \, d\uv}{\int \ex^{\lgd(\xv;\uv)} \, d\uv}
	= 
	\frac{\int_{\UVL^{c}} \ex^{\lgd(\xv;\uv)} \, d\uv}{\int_{\UVL} \ex^{\lgd(\xv;\uv)} \, d\uv + \int_{\UVL^{c}} \ex^{\lgd(\xv;\uv)} \, d\uv}
	\\
	& \leq &
	\frac{4 \ex^{-\xx - \dimLG/2} \int \ex^{- \| \DV_{\GP} \uv \|^{2}/2} d\uv}
		 {\ex^{ - \dltwb \, \dimLG/2} \, \int_{\UVL} \ex^{- \| \DV_{\GP} \uv \|^{2}/2} d\uv 
			+ 4 \ex^{-\xx - \dimLG/2} \int \ex^{- \| \DV_{\GP} \uv \|^{2}/2} d\uv}
	\leq 
	4 \ex^{-\xx - (1 - \dltwb) \dimLG/2} 
\label{jyfftc434rdssdsaaa}
\end{EQA}
as required in \eqref{emdw3dp0iUfxuL}.
\end{proof}

\Section{Finalizing the proof of Theorem~\ref{TLaplaceTV}}
Theorem~\ref{TLaplaceTV} is proved by compiling the previous technical statements.
Proposition~\ref{PUVstarLLfx} provides some upper bounds for the quantities \( \rho \) and \( \rho_{\GP} \),
while Proposition~\ref{Lintfxupp3} and Proposition~\ref{Perrbound3} bound the local errors \( \err \) and \( \err_{g} \).
The final bound \eqref{ufgdt6df5dtgededsxd23gjg} follows from Corollary~\ref{CPunbintLapl}.

\Section{Proof of Theorem~\ref{TLaplaceTV34}}
We make use of the following lemma.
\begin{lemma}
\label{LdltwLa}
Assume \nameref{LLcS3ref}.
Then
\begin{EQA}
	\dltwb
	& \leq &
	\dltwa/3
	\leq 
	\hmax_{3} \, \rrL \, n^{-1/2} / 3 \, 
\label{gtcdsftdfvtwdsefhfdvfrvsewse}
\end{EQA}
and \nameref{LL3fref} holds with 
\begin{EQA}
	\dltwu_{3}(\uv) 
	&=& 
	\hmax_{3} \, n^{-1/2} \| \DVL \uv \|^{3} .
\label{vcdhfgv6e3if5fgfwedgdf}
\end{EQA}
Moreover, \nameref{LLcS4ref} implies 
\nameref{LL4fref} with \( \dltwu_{4}(\uv) = \hmax_{4} \, n^{-1} \| \DVL \uv \|^{4} \), and
\begin{EQA}
	\E \dltwu_{3}(\gaussv_{\GP})
	& \leq &
	\hmax_{3} \, n^{-1/2} (\dimLG + 1)^{3/2} \, ,
	\\
	\E \dltwu_{4}(\gaussv_{\GP})
	& \leq &
	\hmax_{4} \, n^{-1} (\dimLG + 1)^{2} \, .
\label{eciu8hef8h8we3hy87u3y7y3fe7}
\end{EQA}
\end{lemma}

\begin{proof}
By definition, for any \( \uv \in \UVL \) and \( t \in [0,1] \), it holds 
\begin{EQA}
	\frac{\bigl| \langle \nabla^{3} \lgdL(\xvs + t \uv), \uv^{\otimes 3} \rangle \bigr|}{\| \DVL \uv \|^{2}}
	&=&
	\frac{n \bigl| \langle \nabla^{3} \hL(\xvs + t \uv), \uv^{\otimes 3} \rangle \bigr|}{n \langle \nabla^{2} \hL(\xvs), \uv^{\otimes 2} \rangle}
	\leq 
	\hmax_{3} \langle \nabla^{2} \hL(\xvs), \uv^{\otimes 2} \rangle^{1/2}
	\\
	&=& 
	\hmax_{3} \, n^{-1/2} \| \DVL \uv \|
	\leq 
	\hmax_{3} \, \rrL \, n^{-1/2} 
\label{jrgeteteer2234587654}
\end{EQA}
and the first assertion follows.
The second and third ones are proved similarly.
Further, by Lemma~\ref{Gaussmoments} with \( \BB_{\GP} = \DVL \, \DV_{\GP}^{-2} \DVL \leq \Id_{\dimp} \)
\begin{EQA}
	\E \| \DVL \, \gaussv_{\GP} \|^{3}
	& \leq &
	\E^{3/4} \| \DVL \, \gaussv_{\GP} \|^{4}
	\leq 
	\bigl\{ \tr^{2} (\BB_{\GP}) + 2 \tr (\BB_{\GP}^{2}) \bigr\}^{3/4}
	\leq 
	\bigl\{ \tr^{2} (\BB_{\GP}) + 2 \tr (\BB_{\GP}) \bigr\}^{3/4}
	\\
	& < &
	\bigl\{ \tr (\BB_{\GP}) + 1 \bigr\}^{3/2}
	=
	(\dimLG + 1)^{3/2} 
\label{gvdftrfgterfyerhtrfertlfjugfyt}
\end{EQA}
and similarly 
\begin{EQA}
	\E \dltwu_{4}(\gaussv_{\GP})
	& = &
	\hmax_{4} \, n^{-1} \E \| \DVL \, \gaussv_{\GP} \|^{4} 
	\leq 
	\hmax_{4} \, n^{-1} (\dimLG + 1)^{2} \, 
\label{hsdufdkfnhdygfwydfwydgf}
\end{EQA}
yielding \eqref{eciu8hef8h8we3hy87u3y7y3fe7}.
\end{proof}

Lemma~\ref{LdltwLa} implies \( \dltwb \leq \hmax_{3} \, \rrL \, n^{-1/2}/3 \leq 1/4 \) and
\begin{EQA}
	4 \err_{3}
	& \leq &
	\frac{\hmax_{3} }{(1 - \dltwb)^{3/2} n^{1/2}} \,\, \E \| \DVL \, \gaussv_{\GP} \|^{3} 
\label{pcdkjdyt3335r6drgew2eff}
\end{EQA}
yielding \( 4 \err_{3} \leq 2 \hmax_{3} \, \sqrt{(\dimLG+1)^{3}/n} \). 
Similarly under \nameref{LLcS4ref}
\begin{EQA}
	\E \bigl\langle \nabla^{3} \lgd(\xvs) , \gaussv_{\GP}^{\otimes 3} \bigr\rangle^{2}
	& \leq &
	\frac{\hmax_{3}^{2}}{n} \, \E \| \DVL \, \gaussv_{\GP} \|^{6}
	\\
	&=&
	\frac{\hmax_{3}^{2}}{n} \, 
	\bigl\{ (\tr \BB_{\GP})^{3} + 6 \tr \BB_{\GP} \,\, \tr \BB_{\GP}^{2} + 8 \tr \BB_{\GP}^{3} \bigr\}
	\leq 
	\frac{\hmax_{3}^{2} (\dimG + 2)^{3}}{n} \, 
\label{9cuc43ds2s23sgtwsywjnw}
\end{EQA}
and statement \eqref{scdugfdwyd2wywy26e6de} follows from \eqref{ufgdt6df5dtgededsxd23gjg}
and \eqref{hg25t6mxwhydseg3hhfdr} with \( \dltwb \leq 1/4 \) and Lemma~\ref{LdltwLa}.

\Section{Proof of Theorem~\ref{TLaplaceKL} and Theorem~\ref{TLaplaceKLi}}
Define
\begin{EQA}
	\CDG 
	& \eqdef &
	\log \int \ex^{- \| \DV_{\GP} \uv \|^{2}/2} \, d\uv 
	- \log \int \ex^{\lgd(\xvs;\uv)} \, d\uv 
	\, . 
\label{sihedyfw3ytw3ydqwdbhwd}
\end{EQA}
For \( \uv = \xv - \xvs \), it holds as in \eqref{IIgfifxutudufxutdu}
\begin{EQA}
	\PfL
	& \sim &
	\frac{\ex^{\lgd(\xv)}}{\int \ex^{\lgd(\xv)} \, d\uv}
	=
	\frac{\ex^{\lgd(\xv) - \lgd(\xvs)}}{\int \ex^{\lgd(\xv) - \lgd(\xvs)} \, d\uv}
	=
	\frac{\ex^{\lgd(\xvs;\uv)}}{\int \ex^{\lgd(\xvs;\uv)} \, d\uv} 
	=
	\frac{\ex^{\lgd(\xvs;\uv) + \CDG}}{\int \ex^{- \| \DV_{\GP} \uv \|^{2}/2} \, d\uv} \, .
\label{dfhuwffwt5wtcsbydyqyqh}
\end{EQA}
Further, with \( \P_{\GP} = \ND(\xvs,\DV_{\GP}^{-2}) \)
\begin{EQA}
	\log \frac{d\PfL}{d\P_{\GP}}(\xv) 
	&=&
	\lgd(\xvs;\uv) - \| \DV_{\GP} \uv \|^{2}/2 + \CDG
	=
	\dltw_{3}(\uv) + \CDG 
\label{sufuwdeuquqeygwyg22e24}
\end{EQA}
and
\begin{EQA}
	\kullb(\PfL,\P_{\GP})
	&=&
	\ex^{\CDG} \, \frac{\int \dltw_{3}(\uv) \, \ex^{\lgd(\xvs;\uv)} \, d\uv}
		 {\int \ex^{- \| \DV_{\GP} \uv \|^{2}/2} \, d\uv}
	+ \CDG \, .
\label{dvjnue37e37e7weedhe}
\end{EQA}
Similarly to \eqref{Pht12u212H02u22} and \eqref{Pd3p2d3H0wem22g1}, we can bound
\begin{EQA}
	\left| \int_{\UVL} \dltw_{3}(\uv) \, \ex^{\lgd(\xvs;\uv)} \, d\uv \right|
	& \leq &
	\int_{\UVL} |\dltw_{3}(\uv)| \, \ex^{- \| \DV_{\GP} \uv \|^{2}/2 - \dltwb \| \DVL \uv \|^{2}/2} \, d\uv 
	\\
	& \leq &
	\frac{1}{4(1 - \dltwb)^{3/2}}
	\int \dltwu_{3}(\uv) \, \ex^{- \| \DV_{\GP} \uv \|^{2} /2} \, d\uv \, .
\label{ucuwudhuwhi2wefhw78w}
\end{EQA}
On the complement \( \uv \in \UVL^{c} \), we use that
\begin{EQA}
	\dltw_{3}(\uv)
	=
	\lgd(\xvs;\uv) + \frac{1}{2} \| \DV_{\GP} \uv \|^{2}
	&=&
	\lgdL(\xvs;\uv) + \frac{1}{2} \| \DVL \uv \|^{2}
	\leq 
	\frac{1}{2} \| \DVL \uv \|^{2} .
\label{yuqt2sstwwtwbstdywhsyw}
\end{EQA}
The last inequality is based on concavity of \( \lgdL(\cdot) \) and local approximation 
\( \lgdL(\xvs;\uv) \approx - \| \DVL \uv \|^{2}/2 \) within \( \UVL \) yielding \( \lgdL(\xvs;\uv) < 0 \) for \( \uv \in \UVL^{c} \).
By Proposition~\ref{PUVstarLLgenlifx} with \( m = 2 \),
\begin{EQA}
	\int_{\UVL^{c}} \dltw_{3}(\uv) \, \ex^{\lgd(\xvs;\uv)} \, d\uv
	& \leq &
	\frac{1}{2} \int_{\UVL^{c}} \| \DVL \uv \|^{2} \, \ex^{\lgd(\xvs;\uv)} \, d\uv
	\leq 
	\ex^{-\xx} \int \ex^{- \| \DV_{\GP} \uv \|^{2}/2} \, d\uv \, .
\label{sywdywydwbwtwgqb1qt213ba}
\end{EQA}
We conclude that
\begin{EQA}
	\frac{\int \dltw_{3}(\uv) \, \ex^{\lgd(\xvs;\uv)} \, d\uv}
		 {\int \ex^{- \| \DV_{\GP} \uv \|^{2}/2} \, d\uv}
	& \leq &
	\frac{\E \dltwu_{3}(\gaussv_{\GP})}{4(1 - \dltwb)^{3/2}} + \ex^{- \xx} 
	\leq 
	\err_{3} + \ex^{-\xx} .
\label{0fmdfeyjdwhqkshyxiwck}
\end{EQA}
Similarly
\begin{EQA}
	&& \nquad
	\bigl| \ex^{\CDG} - 1 \bigr|
	=
	\left| \frac{\int \ex^{\lgd(\xvs;\uv)} \, d\uv - \int \ex^{- \| \DV_{\GP} \uv \|^{2}/2} \, d\uv}
		 {\int \ex^{- \| \DV_{\GP} \uv \|^{2}/2} \, d\uv} 
	\right|
	\\
	& \leq &
	\left| \frac{\int_{\UVL} \ex^{\lgd(\xvs;\uv)} \, d\uv - \int_{\UVL} \ex^{- \| \DV_{\GP} \uv \|^{2}/2} \, d\uv}
		 {\int \ex^{- \| \DV_{\GP} \uv \|^{2}/2} \, d\uv} 
	\right|
	+ \frac{\int_{\UVL^{c}} \ex^{\lgd(\xvs;\uv)} \, d\uv}{\int \ex^{- \| \DV_{\GP} \uv \|^{2}/2} \, d\uv}
	+ \frac{\int_{\UVL^{c}} \ex^{- \| \DV_{\GP} \uv \|^{2}/2} \, d\uv}
		 {\int \ex^{- \| \DV_{\GP} \uv \|^{2}/2} \, d\uv}
	\\
	& \leq &
	\frac{\E \dltwu_{3}(\gaussv_{\GP})}{4(1 - \dltwb)^{3/2}} + 2 \ex^{- \xx} 
	\leq 
	\err_{3} + 2 \ex^{-\xx} .
\label{hsyc4w2tde6dyfdgyww23}
\end{EQA}
This yields \( \ex^{\CDG} \leq 1 + \err_{3} + 2 \ex^{-\xx} \) and \( \CDG \leq \err_{3} + 2 \ex^{-\xx} \).
Putting all together results in
\begin{EQA}
	\kullb(\PfL,\P_{\GP})
	& \leq &
	\bigl( 1 + \err_{3} + 2 \ex^{-\xx} \bigr) \bigl( \err_{3} + \ex^{-\xx} \bigr) + \err_{3} + 2 \ex^{-\xx}
	<  
	4 \err_{3} + 4 \ex^{-\xx} 
\label{gcusgu255eetfghbhjjkkook}
\end{EQA}
provided that \( 4 \err_{3} + 4 \ex^{-\xx} \leq 1 \),
and \eqref{vlgvi8ugu7tr4ry43et31} follows.

The proof of Theorem~\ref{TLaplaceKLi} is similar and even simpler except one special part, namely, 
the bound for the tail integral of \( - \dltw_{3}(\uv) \).
By definition \( \| \DVL \uv \| \geq \rr_{\GP} \) for \( \uv \in \UVL^{c} \).
This implies by \eqref{dchbhwdhwwdgscsn2efty2162} similarly to \eqref{yuqt2sstwwtwbstdywhsyw} 
\begin{EQA}
	&& \nquad
	- \int_{\UVL^{c}} \dltw_{3}(\uv) \, \ex^{- \| \DV_{\GP} \uv \|^{2}/2} \, d\uv
	\leq
	\int_{\UVL^{c}} \bigl| \lgdL(\xvs;\uv) \bigr| \, \ex^{- \| \DV_{\GP} \uv \|^{2}/2 + \rho_{\GP} \| \DVL \uv \|^{2}/2} \,
		\ex^{- \rho_{\GP} \| \DVL \uv \|^{2}/2} \, d\uv
	\\
	& \leq &
	\ex^{- \rho_{\GP} \rr_{\GP}^{2}/2} \, 
	\int \bigl| \lgdL(\xvs;\uv) \bigr| \, \ex^{- \| \DV_{\GP} \uv \|^{2}/2 + \rho_{\GP} \| \DVL \uv \|^{2}/2} \,
		d\uv
	\leq \CONSTi_{\lgdL} \, \ex^{-\xx} ,
\label{sguswdw2e2j2wuhq2wdw}
\end{EQA}
and the result follows.

\Section{Finalizing the proof of Theorem~\ref{TpostmeanLa} and Corollary~\ref{CTpostmeanLa}}
For Theorem~\ref{TpostmeanLa}, we follow the same line as for Theorem~\ref{TLaplaceTV34}.
Note first that
\begin{EQA}
	\QP (\xvb - \xvs)
	=
	\frac{\int \QP(\xvs + \uv) \, \ex^{\lgd(\xvs + \uv)} \, d\uv}{\int \ex^{\lgd(\xvs + \uv)} \, d\uv} - \QP \xvs
	&=&
	\frac{\int \QP \uv \, \ex^{\lgd(\xvs;\uv)} \, d\uv}{\int \ex^{\lgd(\xvs;\uv)} \, d\uv} \, 
\label{dvhjt6efedfchsijcfte4ws}
\end{EQA}
and
\begin{EQA}
	\| \QP (\xvb - \xvs) \|
	&=&
	\sup_{\av \in \R^{\dimq} \colon \| \av \| = 1} \bigl| \langle \QP (\xvb - \xvs), \av \rangle \bigr|
	= 
	\sup_{\av \in \R^{\dimq} \colon \| \av \| = 1} 
		\left| \frac{\int \langle \QP \uv, \av \rangle \, \ex^{\lgd(\xvs;\uv)} \, d\uv}{\int \ex^{\lgd(\xvs;\uv)} \, d\uv} \right| \, .
\label{p3hb893dfvfdsiuw5whlh}
\end{EQA}
Now fix \( \av \in \R^{\dimq} \) with \( \| \av \| = 1 \) and \( g(\uv) = \langle \QP \uv, \av \rangle \).
\eqref{igefxumiguexHu22g} implies
\begin{EQA}[rcl]
	\left| \frac{\int g(\uv) \, \ex^{\lgd(\xvs;\uv)} \, d\uv}{\int \ex^{\lgd(\xvs;\uv)} \, d\uv} \right| 
	& \leq &
	\frac{\rho_{g} + \rho_{\GP,g} + \err_{3,g}}{1 - \rho_{\GP} - \err_{3}} \, .
\label{igefxumiguexHu22gp}
\end{EQA}
The bound \( 1 - \err_{3} - \rho_{\GP} \geq 1/2 \) has been already checked.
Proposition~\ref{PUVstarLLgenlifx} for \( m=1 \) helps to bound the values \( \rho_{g} \) and \( \rho_{\GP,g} \) 
by \( \CONST \ex^{-\xx} \).
Next we bound \( \err_{3,g} \).
Under \nameref{LLcS3ref},
\eqref{iUafxvdvH2vk3} of Proposition~\ref{Perrbound3} combined with Lemma~\ref{LdltwLa} and Lemma~\ref{Gaussmoments} yield
\begin{EQA}
	&& \nquad
	4 \err_{3,g}
	=
	\frac{1}{(1 - \dltwb)^{2}} \,
	\E \bigl\{  | \langle \QP \gaussv_{\GP},\av \rangle | \,\, \dltwu_{3}(\gaussv_{\GP}) \bigr\}
	= 
	\frac{\hmax_{3} \, n^{-1/2}}{(1 - \dltwb)^{2}} \, 
	\E \bigl\{ | \langle \QP \, \gaussv_{\GP},\av \rangle | \,\, \| \DVL \, \gaussv_{\GP} \|^{3} \, \bigr\} 
	\\
	& \leq &
	\frac{\hmax_{3} \, n^{-1/2}}{(1 - \dltwb)^{2}} \, 
	\E^{3/4} \| \DVL \, \gaussv_{\GP} \|^{4} \,\, \E^{1/4} \langle \QP \gaussv_{\GP},\av \rangle^{4} 
	\\
	& \leq &
	\frac{3^{1/4} \, \hmax_{3} \, (\dimLG + 1)^{3/2} \, n^{-1/2}}{(1 - \dltwb)^{2}} \, 
	\bigl( \av^{\T} \QP \, \DV_{\GP}^{-2} \QP^{\T} \av \bigr)^{1/2} .
\label{sfdhysdfyf11ewds3wsded}
\end{EQA}
Here we used that \( \langle \QP \gaussv_{\GP},\av \rangle \sim \ND(0,\av^{\T} \QP \, \DV_{\GP}^{-2} \QP^{\T} \av) \)
and \( \E \langle \QP \gaussv_{\GP},\av \rangle^{4} = 3 (\av^{\T} \QP \, \DV_{\GP}^{-2} \QP^{\T} \av)^{2} \).
Now \eqref{hcdtrdtdehfdewdrfrhgyjufger} follows from \( 3^{1/4} (1 - \dltwb)^{-2} \leq 2.4 \) and
\begin{EQA}
	\sup_{\av \in \R^{\dimq} \colon \| \av \| = 1} \av^{\T} \QP \, \DV_{\GP}^{-2} \QP^{\T} \av
	&=&
	\| \QP \, \DV_{\GP}^{-2} \QP^{\T} \| \, .
\label{hdrtesw5sdghww3r4d5tdy6gh}
\end{EQA}
With \( \QP = \DVL \), this implies \eqref{klu8gitfdgregfkhj7yt}.
\Section{Proof of Theorem~\ref{TLaplnonlin}}
The basic idea of the proof is to reduce the statements of the theorem to the case of Theorem~\ref{TLaplaceTV}.
Conceptually, the most important step of the proof is to show that the function \( \lgd(\xv) \) is strongly concave
on the subset \( \XX_{0} \) which is a concentration set of \( \ND(\xv_{0},\GP^{-2}) \).

\begin{lemma}
\label{SconccalmoL}
Suppose \nameref{LLS2ref}, \nameref{rrinref}, and \nameref{theta0ref}.
Then for all \( \xv \in \XX_{0} \)
\begin{EQA}
	\IFVL(\xv)
	& \geq &
	(1 - \delta ) \, \DVLb^{2}(\xv) 
	\quad
	\text{with}
	\quad
	\delta
	=
	\CONSTi_{2} \, \bigl( \smlc_{0} + 2 \rrin \, \CONSTi_{n} \, \CONSTi_{0} \, /\sqrt{n} \bigr) 
	\leq 
	3/4 \, .
\label{we9furw3efy8ewsifhswi8dL}
\end{EQA}
\end{lemma}

\begin{proof}
By \eqref{hjwetwee3ee5555t566dsL}
\begin{EQA}
	\sumi \left| \bigl\{ \regrf_{i}(\xv) - z_{i} \bigr\} 
		\bigl\langle \nabla^{2} \regrf_{i}(\xv), \uv^{\otimes 2} \bigr\rangle 
	\right|
	& \leq &
	\CONSTi_{2} \| \regrfv(\xv) - \zv \|_{\infty} \,\, \| \DVb(\xv) \uv \|^{2} .
\label{0cndt2sghsdhtswghwsL}
\end{EQA}
Further, let \( \xv_{*} \) and \( \smlc_{0} \) be such that 
\( \| \regrfv(\xv) - \regrfv(\xv_{*}) \|_{\infty} \leq \smlc_{0} \).
Then 
\begin{EQA}
	\| \regrfv(\xv) - \zv \|_{\infty}
	& \leq &
	\| \regrfv(\xv) - \regrfv(\xv_{*}) \|_{\infty}
	+ \| \zv - \regrfv(\xv_{*}) \|_{\infty}
	\leq 
	\| \regrfv(\xv) - \regrfv(\xv_{*}) \|_{\infty} + \smlc_{0} \, .
\label{y8gwefuyfwqeguy8efwhiueL}
\end{EQA}
Given \( \xv \in \XX_{0} \), denote \( \uv = \xv - \xv_{*} \).
The definition of \( \XX_{0} \) implies \( \| \DVLbin \uv \| \leq 2 \rrin \).
The first order Taylor expansion yields
\begin{EQA}
	\regrf_{i}(\xv) - \regrf_{i}(\xv_{*})
	&=&
	\bigl\langle \nabla \regrf_{i}(\xvd), \uv \bigr\rangle 
\label{fqe9f8gefyeuyefuygydeL}
\end{EQA}
with \( \xvd = \xvs + t \uv \) for \( t \in [0,1] \).
Hence, by \eqref{wfegy7r5qrw35edfhgdyfysL} and \eqref{fchghiu87686574e5rtyyuL}
\begin{EQA}
	&& \nquad
	| \regrf_{i}(\xv) - \regrf_{i}(\xv_{*}) |
	=
	\bigl| \langle \nabla \regrf_{i}(\xvd), \uv \rangle \bigr|
	=
	\bigl| \langle \DVLbin^{-1} \, \nabla \regrf_{i}(\xvd), \DVLbin \uv \rangle \bigr|
	\\
	& \leq & 
	\| \DVLbin^{-1} \, \nabla \regrf_{i}(\xvd) \, \nabla \regrf_{i}(\xvd)^{\T} \DVLbin^{-1} \|^{1/2} \,\, \| \DVLbin \uv \|
	\leq 
	2 \rrin \, \CONSTi_{n} \, \CONSTi_{0} \, /\sqrt{n} 
\label{ggtdg2321323213689768768L}
\end{EQA}
yielding 
\begin{EQA}
	\| \regrfv(\xv) - \zv \|_{\infty}
	& \leq &
	2 \rrin \, \CONSTi_{n} \, \CONSTi_{0} \, /\sqrt{n} + \smlc_{0}
\label{wygyw3537df7ugusdgdgywydyL}
\end{EQA} 
and \eqref{we9furw3efy8ewsifhswi8dL}
follows by \eqref{bvudufgheuygfdtcsdwedeuf} and \eqref{0cndt2sghsdhtswghwsL}. 
\end{proof}

The next step is to show that 
under \nameref{theta0ref},
the point \( \xvs_{\GP} \) is also within \( \XX_{0} \).

\begin{lemma}
\label{Lnonlintrue}
Assume \nameref{theta0ref} with \( \smlc_{0} = 0 \).
Then \( \xvs_{\GP} \in \XX_{0} \).
\end{lemma}

\begin{proof}
Let \( \xvs \in \XX_{0} \) solve \( \regrfv(\xvs) = \zv \).
By definition, \( \| \DVLbin (\xvs - \xv_{0}) \| \leq \rrin \).
As \( \xvs_{\GP} \) minimizes the criteria \( \| \GP (\xv - \xv_{0}) \|^{2} + \| \zv - \regrfv(\xv) \|^{2} \), it holds
\begin{EQA}
	\| \GP (\xvs_{\GP} - \xv_{0}) \|^{2}
	\leq 
	\| \GP (\xvs_{\GP} - \xv_{0}) \|^{2} + \| \zv - \regrfv(\xvs_{\GP}) \|^{2} 
	& \leq &
	\| \GP (\xvs - \xv_{0}) \|^{2} 
\label{fkjeue7f4e65egdewghvceh}
\end{EQA}
This implies the assertion in view of \( n^{-1} \DVLbin^{2} \leq \GP^{2} \).
\end{proof}

Now we explain how local smoothness properties of \( \elll(\xv) = - \| \regrfv(\xv) - \zv \|^{2}/2 \)
can be characterized under \nameref{LLS2ref} and \nameref{LLS3ref} or \nameref{LLS4ref}.

\begin{lemma}
\label{LsmoothcalmL}
Assume \nameref{LLS2ref} and \nameref{LLS3ref}. 
Then for any \( \xv \in \UVLb \) and any \( \uv \in \R^{\dimp} \)
\begin{EQA}
	\bigl| \langle \nabla^{3} \elll(\xv), \uv^{\otimes 3} \rangle  \bigr|
	& \leq &
	\frac{4 \CONSTi_{\GP}^{3/2} \, \CONSTi_{3} \, \CONSTi_{n}}{\sqrt{n}} \, \| \DVb \uv \|^{3} \, .
\label{voedf8efyffte4er4wgwL}
\end{EQA}
Moreover, if \nameref{LLS4ref} is also valid, then
\begin{EQA}
	\bigl| \langle \nabla^{4} \elll(\xv), \uv^{\otimes 4} \rangle \bigr|
	& \leq &
	\frac{8 \CONSTi_{\GP}^{2} \, \CONSTi_{4} \, \CONSTi_{n}}{n} \, \| \DVb \uv \|^{4}.
\label{voedf8efyffte4er4wgwjw4L}
\end{EQA}
\end{lemma}

\begin{proof}
Represent \( \elll(\xv) = - \| \regrfv(\xv) - \zv \|^{2}/2 \) as a sum \( \elll(\xv) = - \sumi \elll_{i}(\xv) \) with
\( \elll_{i}(\xv) = | \regrf_{i}(\xv) - z_{i} |^{2}/2 \).
Then for any direction \( \uv \in \R^{\dimp} \) and for each summand \( \elll_{i}(\xv) \), it holds
\begin{EQA}
    \langle \nabla^{3} \elll_{i}(\xv), \uv^{\otimes 3} \rangle 
    & = &
    - 3 \langle \nabla \regrf_{i}(\xv),\uv \rangle \,
    	\bigl\langle \nabla^{2} \regrf_{i}(\xv),\uv^{2} \bigr\rangle
    - \bigl\{ \regrf_{i}(\xv) - z_{i} \bigr\} \, \langle \nabla^{3} \regrf_{i}(\xv),\uv^{3} \rangle .
    \qquad
\label{d3utd4utd4dt4cmL}
\end{EQA}
Therefore,
\begin{EQA}
	\bigl| \bigl\langle \nabla^{3} \elll(\xv), \uv^{\otimes 3} \bigr\rangle \bigr|
	& \leq &
	3 \left| \sumi \langle \nabla \regrf_{i}(\xv),\uv \rangle \,
    	\bigl\langle \nabla^{2} \regrf_{i}(\xv),\uv^{2} \bigr\rangle 
	+ \sumi \bigl\{ \regrf_{i}(\xv) - z_{i} \bigr\} \, 
		\langle \nabla^{3} \regrf_{i}(\xv),\uv^{3} \rangle \right| .
\label{hs8d8s8sduhudhuhdsuL}
\end{EQA}
By \nameref{LLS3ref}
\begin{EQA}
	\left| \sumi \langle \nabla \regrf_{i}(\xv),\uv \rangle \,
    	\bigl\langle \nabla^{2} \regrf_{i}(\xv),\uv^{2} \bigr\rangle \right|
	& \leq &
	\CONSTi_{3} 
	\sumi \bigl| \langle \nabla \regrf_{i}(\xv),\uv \rangle \bigr|^{3} .
\label{werd5wr5werw6ed37tr37tr73L}
\end{EQA}
Similarly
\begin{EQA}
	\left| \sumi \bigl\{ \regrf_{i}(\xv) - z_{i} \bigr\} \, \langle \nabla^{3} \regrf_{i}(\xv),\uv^{3} \rangle \right|
	& \leq &
	\CONSTi_{3} \max_{i \leq n} 
	\Bigl| \{ \regrf_{i}(\xv) - z_{i} \} \Bigr| \, \,  
	\sumi \bigl| \langle \nabla \regrf_{i}(\xv),\uv \rangle \bigr|^{3} \, .
\label{079u8760u9jhkigfetr6wL}
\end{EQA}
By \eqref{fchghiu87686574e5rtyyuLL},
\begin{EQA}
	\sumi \bigl| \langle \nabla \regrf_{i}(\xv),\uv \rangle \bigr|^{3}
	& \leq &
	\max_{i \leq n} \bigl| \langle \nabla \regrf_{i}(\xv),\uv \rangle \bigr| \,\,
	\sumi \langle \nabla \regrf_{i}(\xv),\uv \rangle^{2} 
	\\
	&=&
	\max_{i \leq n} \bigl| \langle \nabla \regrf_{i}(\xv),\uv \rangle \bigr| \,\,
	\| \DVb(\xv) \uv \|^{2}
	\leq 
	\CONSTi_{\GP}
	\max_{i \leq n} \bigl| \langle \nabla \regrf_{i}(\xv),\uv \rangle \bigr| \,\,
	\| \DVb \uv \|^{2}.
\label{wgfs7ws7w7td7wre2w5vhuL}
\end{EQA}
Also for any \( \uv \) with \( \| \DVLb \uv \| \leq \rr \), by \eqref{wfegy7r5qrw35edfhgdyfysL} of \nameref{LLS2ref} 
and again by \eqref{fchghiu87686574e5rtyyuLL}
\begin{EQA}
	\bigl| \langle \nabla \regrf_{i}(\xv),\uv \rangle \bigr|
	&=&
	\bigl| \langle \DVLb^{-1}(\xv) \, \nabla \regrf_{i}(\xv), \DVLb(\xv) \uv \rangle \bigr|
	\\
	& \leq &
	\| \DVLb(\xv) \uv \| \, \| \DVLb^{-1}(\xv) \, \nabla \regrf_{i}(\xv) \| 
	\leq 
	\frac{\sqrt{\CONSTi_{\GP}} \, \CONSTi_{n}}{\sqrt{n}} \, \| \DVLb \uv \|\, .
\label{dfof7f7ye4ge4trwshydyqL}
\end{EQA}
Furthermore, \( \| \regrfv(\xv) - \zv \|_{\infty} \leq 1 \) for all \( \xv \in \XX_{0} \) 
as in the proof of Lemma~\ref{SconccalmoL}, and \eqref{voedf8efyffte4er4wgwL} 
The proof of \eqref{voedf8efyffte4er4wgwjw4L} is similar with the use of
\begin{EQA}
    \langle \nabla^{4} \elll_{i}(\xv), \uv^{\otimes 4} \rangle 
    & = &
    - 3 \bigl\langle \nabla^{2} \regrf_{i}(\xv),\uv^{2} \bigr\rangle^{2}
    - 4 \langle \nabla \regrf_{i}(\xv),\uv \rangle \, \bigl\langle \nabla^{3} \regrf_{i}(\xv),\uv^{3} \bigr\rangle
    \\
    &&
    - \, \bigl\{ \regrf_{i}(\xv) - z_{i} \bigr\}  
    \bigl\langle \nabla^{4} \regrf_{i}(\xv),\uv^{4} \bigr\rangle 
\label{ncdgtdsrsdds3s4e3sewserw}	
\end{EQA}
and of \nameref{LLS4ref}.
\end{proof}

Bound \eqref{voedf8efyffte4er4wgwL} means that \eqref{gtcdsftdfvtwdsefhfdvfrvsewse3} is fulfilled with 
\( \hmax_{3} = 4 \CONSTi_{\GP}^{3/2} \, \CONSTi_{3} \, \CONSTi_{n} \).
This implies \eqref{cbc5dfedrdwewwgerg3NL} similarly to Theorem~\ref{TLaplaceTV34}.
Moreover, \eqref{voedf8efyffte4er4wgwjw4L} implies \nameref{LLcS4ref} with
\( \hmax_{4} = 8 \CONSTi_{\GP}^{2} \, \CONSTi_{4} \, \CONSTi_{n} \),
and the second statement of the theorem follows as well.

\Chapter{Some results for Gaussian quadratic forms}
\def\HVB{\mathcal{V}}
\Section{Moments of a Gaussian quadratic form}
\label{Sdevboundgen}
Let \( \gaussv \) be standard normal in \( \R^{\dimp} \) for \( \dimp \leq \infty \).
Given a self-adjoint trace operator \( \BB \), consider a quadratic form 
\( \bigl\langle \BB \gaussv, \gaussv \bigr\rangle \).

\begin{lemma}
\label{Gaussmoments}
It holds
\begin{EQA}
	\E \bigl\langle \BB \gaussv, \gaussv \bigr\rangle 
	&=& 
	\tr \BB ,
	\\ 
	\Var \bigl\langle \BB \gaussv, \gaussv \bigr\rangle 
	&=& 
	2 \tr \BB^{2} .
\label{EAarAtrA2trA2}
\end{EQA}
Moreover, 
\begin{EQA}
	\E \bigl( \bigl\langle \BB \gaussv, \gaussv \bigr\rangle - \tr \BB \bigr)^{2}
	&=&
	2 \tr \BB^{2}  ,
	\\
	\E \bigl( \bigl\langle \BB \gaussv, \gaussv \bigr\rangle - \tr \BB \bigr)^{3}
	&=&
	8 \tr \BB^{3} ,
	\\
	\E \bigl( \bigl\langle \BB \gaussv, \gaussv \bigr\rangle - \tr \BB \bigr)^{4}
	&=&
	48 \tr \BB^{4} + 12 (\tr \BB^{2})^{2} ,
\label{2pG2trD2DGm22m2}
\end{EQA}
and
\begin{EQA}
	\E \bigl\langle \BB \gaussv, \gaussv \bigr\rangle^{2}
	&=&
	(\tr \BB)^{2} + 2 \tr \BB^{2},
	\\
	\E \bigl\langle \BB \gaussv, \gaussv \bigr\rangle^{3}
	& = &
	(\tr \BB)^{3} + 6 \tr \BB \,\, \tr \BB^{2} + 8 \tr \BB^{3} ,
	\\
	\E \bigl\langle \BB \gaussv, \gaussv \bigr\rangle^{4}
	& = &
	(\tr \BB)^{4} + 12 (\tr \BB)^{2} \tr \BB^{2}
	+ 32 (\tr \BB) \tr \BB^{3}
	+ 48 \tr \BB^{4} + 12 (\tr \BB^{2})^{2} ,
\label{2pG2trD2DGm22m2}
	\\
	\Var \bigl\langle \BB \gaussv, \gaussv \bigr\rangle^{2}
	& = &
	8 (\tr \BB)^{2} \tr \BB^{2}
	+ 32 (\tr \BB) \tr \BB^{3}
	+ 48 \tr \BB^{4} + 8 (\tr \BB^{2})^{2} .
\label{2pG2trD2DGm22m4}
\end{EQA}
Moreover, if \( \BB \leq \Id_{\dimp} \) and \( \dimA = \tr \BB \), then \( \tr \BB^{m} \leq \dimA \) for 
\( m \geq 1 \) and
\begin{EQA}
	\E \bigl\langle \BB \gaussv, \gaussv \bigr\rangle^{2}
	& \leq &
	\dimA^{2} + 2 \dimA,
	\\
	\E \bigl\langle \BB \gaussv, \gaussv \bigr\rangle^{3}
	& \leq &
	\dimA^{3} + 6 \dimA^{2} + 8 \dimA ,
	\\
	\E \bigl\langle \BB \gaussv, \gaussv \bigr\rangle^{4}
	& = &
	\dimA^{4} + 12 \dimA^{3} 
	+ 44 \dimA^{2}
	+ 48 \dimA  ,
\label{2pG2trD2DGm22m2}
	\\
	\Var \bigl\langle \BB \gaussv, \gaussv \bigr\rangle^{2}
	& = &
	8 \dimA^{3} + 40 \dimA^{2} + 48 \dimA .
\label{2pG2trD2DGm22m4}
\end{EQA}
\end{lemma}

\begin{proof}
Let \( \chi = \gauss^{2} - 1 \) for \( \gauss \) standard normal.
Then \( \E \chi = 0 \), \( \E \chi^{2} = 2 \), \( \E \chi^{3} = 8 \), \( \E \chi^{4} = 60 \).
Without loss of generality assume \( \BB \) diagonal: \( \BB = \diag(\lambda_{1},\lambda_{2},\ldots,\lambda_{\dimp}) \).
Then 
\begin{EQA}
	\xi
	\eqdef
	\bigl\langle \BB \gaussv, \gaussv \bigr\rangle - \tr \BB
	&=&
	\sum_{j=1}^{\dimp} \lambda_{j} (\gauss_{j}^{2} - 1) ,
\label{j1ljgj2m1}
\end{EQA}
where \( \gauss_{j} \) are i.i.d. standard normal. 
This easily yields
\begin{EQA}
	\E \xi^{2}
	&=&
	\sum_{j=1}^{\dimp} \lambda_{j}^{2} \E (\gauss_{j}^{2} - 1)^{2}
	=
	\E \chi^{2} \, \tr \BB^{2} 
	=
	2 \tr \BB^{2}  ,
	\\
	\E \xi^{3}
	&=&
	\sum_{j=1}^{\dimp} \lambda_{j}^{3} \E (\gauss_{j}^{2} - 1)^{3}
	=
	\E \chi^{3} \, \tr \BB^{3} 
	=
	8 \tr \BB^{3} ,
	\\
	\E \xi^{4}
	&=&
	\sum_{j=1}^{\dimp} \lambda_{j}^{4} (\gauss_{j}^{2} - 1)^{4}
	+ \sum_{i\neq j} \lambda_{i}^{2} \lambda_{j}^{2} \E (\gauss_{i}^{2} - 1)^{2} \E (\gauss_{j}^{2} - 1)^{2}
	\\
	&=&
	\bigl( \E \chi^{4} - 3 (\E \chi^{2})^{2} \bigr) \tr \BB^{4} + 3 (\E \chi^{2} \, \tr \BB^{2})^{2}
	=
	48 \tr \BB^{4} + 12 (\tr \BB^{2})^{2} ,
\label{2pG2trD2DGm22m2}
\end{EQA}
ensuring
\begin{EQA}
	\E \bigl\langle \BB \gaussv, \gaussv \bigr\rangle^{2}
	&=&
	\bigl( \E \bigl\langle \BB \gaussv, \gaussv \bigr\rangle \bigr)^{2} 
	+ \E \xi^{2}
	= 
	(\tr \BB)^{2} + 2 \tr \BB^{2},
	\\
	\E \bigl\langle \BB \gaussv, \gaussv \bigr\rangle^{3}
	& = &
	\E \bigl( \xi + \tr \BB \bigr)^{3}
	=
	(\tr \BB)^{3} + \E \xi^{3}
	+ 3 \tr \BB \,\, \E \xi^{2}
	\\
	&=&
	(\tr \BB)^{3} + 6 \tr \BB \,\, \tr \BB^{2} + 8 \tr \BB^{3} ,
\label{2pG2trD2DGm22m2}
\end{EQA}
and 
\begin{EQA}
	\Var \bigl\langle \BB \gaussv, \gaussv \bigr\rangle^{2}
	& = &
	\E \bigl( \xi + \tr \BB \bigr)^{4}
	- \bigl( \E \bigl\langle \BB \gaussv, \gaussv \bigr\rangle \bigr)^{2}
	\\
	&=&
	\bigl( \tr \BB \bigr)^{4} + 6 (\tr \BB)^{2} \E \xi^{2} + 4 \tr \BB \, \E \xi^{3} + \E \xi^{4}
	- \bigl( (\tr \BB)^{2} + 2 \tr \BB^{2} \bigr)^{2}
	\\
	&=& 
	8 (\tr \BB)^{2} \tr \BB^{2}
	+ 32 (\tr \BB) \tr \BB^{3}
	+ 48 \tr \BB^{4} + 8 (\tr \BB^{2})^{2} 
\label{2pG2trD2DGm22m4}
\end{EQA}
This implies the results of the lemma.
\end{proof}

Now we compute the exponential moments of centered and non-centered quadratic forms.

\begin{lemma}
\label{Lqfexpmom}
Let \( \| \BB \|_{\oper} \leq 1 \).
Then for any \( \mu \in (0,1) \), 
\begin{EQA}
	\E \exp \Bigl\{ \frac{\mu}{2} \bigl( \langle \BB \gaussv, \gaussv \rangle - \dimA \bigr) \Bigr\}
	&=&
	\det(\Id - \mu \BB)^{-1/2} \, .
\label{m2v241m41m}
\end{EQA}
Moreover, with \( \dimA = \tr \BB \) and \( \vA^{2} = \tr \BB^{2} \)
\begin{EQA}
	\log \E \exp \Bigl\{ \frac{\mu}{2} \bigl( \langle \BB \gaussv, \gaussv \rangle - \dimA \bigr) \Bigr\}
	& \leq &
	\frac{\mu^{2} \vA^{2}}{4 (1 - \mu)} \, .
\label{m2v241m41mb}
\end{EQA}
If \( \BB \) is positive semidefinite, \( \lambda_{j} \geq 0 \), then 
\begin{EQA}
	\log \E \exp \Bigl\{ - \frac{\mu}{2} \bigl( \langle \BB \gaussv, \gaussv \rangle - \dimA \bigr) \Bigr\}
	& \leq &
	\frac{\mu^{2} \vA^{2}}{4} \, .
\label{m2v241m41mbn}
\end{EQA}
\end{lemma}

\begin{proof}
Let \( \lambda_{j} \) be the eigenvalues of \( \BB \), 
\( |\lambda_{j}| \leq 1 \).
By an orthogonal transform, one can reduce the statement to the case of a diagonal matrix 
\( \BB = \diag\bigl( \lambda_{j} \bigr) \). 
Then \( \langle \BB \gaussv, \gaussv \rangle = \sum_{j=1}^{\dimp} \lambda_{j} \eps_{j}^{2} \) and 
by independence of the \( \eps_{j} \)'s
\begin{EQA}
	&& \nquad
	\E \Bigl\{ \frac{\mu}{2} \langle \BB \gaussv, \gaussv \rangle  \Bigr\}
	=
	\prod_{j=1}^{\dimp} \E \exp \Bigl( \frac{\mu}{2} \lambda_{j} \eps_{j}^{2} \Bigr)
	=
	\prod_{j=1}^{\dimp} \frac{1}{\sqrt{1 - \mu \lambda_{j}}} 
	=
	\det \bigl( \Id - \mu \BB \bigr)^{-1/2} .
\label{dOImuBm12EB}
\end{EQA}
Below we use the simple bound: 
\begin{EQ}[rcl]
\label{lo1uusk2iukkp}
	- \log(1 - u) - u
	&=&
	\sum_{k=2}^{\infty} \frac{u^{k}}{k}
	\leq 
	\frac{u^{2}}{2} \sum_{k=0}^{\infty} u^{k} 
	=
	\frac{u^{2}}{2 (1 - u)} \, ,
	\qquad 
	u \in (0,1),
	\\
	- \log(1 - u) + u
	&=&
	\sum_{k=2}^{\infty} \frac{u^{k}}{k}
	\leq 
	\frac{u^{2}}{2} \, ,
	\qquad \qquad
	u \in (-1,0).
\label{lo1uusk2iukk}
\end{EQ}
Now it holds 
\begin{EQA}
	&& \nquad
	\log \E \Bigl\{ \frac{\mu}{2} \bigl( \langle \BB \gaussv, \gaussv \rangle - \dimA \bigr) \Bigr\}
	=
	\log \det(\Id - \mu \BB)^{-1/2} - \frac{\mu \, \dimA}{2}
	\\
	&=&
	- \frac{1}{2} \sum_{j=1}^{\dimp} \bigl\{ \log(1 - \mu \lambda_{j}) + \mu \lambda_{j} \bigr\}
	\leq 
	\sum_{j=1}^{\dimp} \frac{\mu^{2} \lambda_{j}^{2}}{4 (1 - \mu)} 
	=
	\frac{\mu^{2} \vA^{2}}{4 (1 - \mu)} \, .
\label{m2v241m4mj1pd}
\end{EQA}
The last statement can be proved similarly.
\end{proof}

Now we consider the case of a non-centered quadratic form
\( \langle \BB \gaussv,\gaussv \rangle/2 + \langle \Av,\gaussv \rangle \) for a fixed vector \( \Av \).

\begin{lemma}
\label{Lexpmomnoncen}
Let \( \lambda_{\max}(\BB) < 1 \). 
Then for any \( \Av \)
\begin{EQA}
	\E \exp\Bigl\{ \frac{1}{2}\langle \BB \gaussv,\gaussv \rangle + \langle \Av,\gaussv \rangle \Bigr\}
	&=&
	\exp\Bigl\{ \frac{\| (\Id - \BB)^{-1/2} \Av \|^{2}}{2} \Bigr\} \, \det(\Id - \BB)^{-1/2} .
\label{EeBf12BggA}
\end{EQA}
Moreover, for any \( \mu \in (0,1) \)
\begin{EQA}
	&& \nquad
	\log \E \exp\Bigl\{ 
		\frac{\mu}{2} \bigl( \langle \BB \gaussv,\gaussv \rangle - \dimA \bigr) + \langle \Av,\gaussv \rangle 
	\Bigr\}
	\\
	&=&
	\frac{\| (\Id - \mu \BB)^{-1/2} \Av \|^{2}}{2} + \log \det(\Id - \mu \BB)^{-1/2} - \mu \dimA 
	\\
	& \leq &
	\frac{\| (\Id - \mu \BB)^{-1/2} \Av \|^{2}}{2} + \frac{\mu^{2} \vA^{2}}{4 (1 - \mu)} \, .
\label{EeBf12BggAmu}
\end{EQA}
\end{lemma}

\begin{proof}
Denote \( \av = (\Id - \BB)^{-1/2} \Av \). 
It holds by change of variables \( (\Id - \BB)^{1/2} \xv = \uv \) for \( \CONSTi_{\dimp} = (2\pi)^{-\dimp/2} \)
\begin{EQA}
	&& \nquad
	\E \exp\Bigl\{ \frac{1}{2}\langle \BB \gaussv,\gaussv \rangle + \langle \Av,\gaussv \rangle \Bigr\}
	=
	\CONSTi_{\dimp}
	\int \exp\Bigl\{ - \frac{1}{2}\langle (\Id - \BB) \xv,\xv \rangle + \langle \Av,\xv \rangle \Bigr\} d\xv
	\\
	&=&
	\CONSTi_{\dimp}
	\det(\Id - \BB)^{-1/2}
	\int \exp\Bigl\{ - \frac{1}{2} \| \uv \|^{2} + \langle \av,\uv \rangle \Bigr\} d\uv
	=
	\det(\Id - \BB)^{-1/2} \, 	\ex^{\| \av \|^{2}/2}  	.
\label{EeBf12BggAp}
\end{EQA}
The last inequality \eqref{EeBf12BggAmu} follows by \eqref{m2v241m41mb}.
\end{proof}

\Section{Deviation bounds for Gaussian quadratic forms}
\label{SdevboundGauss}
The next result explains the concentration effect of 
\( \langle \BB \xiv, \xiv \rangle \)
for a centered Gaussian vector \( \xiv \sim \ND(0,\HVB^{2}) \) and a symmetric trace operator \( \BB \) in \( \R^{\dimp} \),
\( \dimp \leq \infty \).
We use a version from \cite{laurentmassart2000}.
For completeness, we present a simple proof of the upper bound.

\begin{theorem}
\label{TexpbLGA}
\label{Lxiv2LD}
\label{Cuvepsuv0}
Let \( \xiv \sim \ND(0,\HVB^{2}) \) be a Gaussian element in \( \R^{\dimp} \) and \( \BB \) be symmetric 
such that \( \BBH = \HVB \BB \HVB \) is a trace operator in \( \R^{\dimp} \).
Then with \( \dimH = \tr(\BBH) \), \( \vH^{2} = \tr(\BBH^{2}) \), and 
\( \supA = \| \BBH \| \), it holds for each \( \xx \geq 0 \)
\begin{EQA}
\label{Pxiv2dimAvp12}
	\P\Bigl( \langle \BB \xiv, \xiv \rangle - \dimH > 2 \vH \, \sqrt{\xx} + 2 \supA \xx \Bigr)
	& \leq &
	\ex^{-\xx} .
\end{EQA}
It also implies 
\begin{EQA}
	\P\bigl( \bigl| \langle \BB \xiv, \xiv \rangle - \dimH \bigr| > \zq_{2}(\BBH,\xx) \bigr)
	& \leq &
	2 \ex^{-\xx} ,
\label{PxivTBBdimA2vp}
\end{EQA}
with
\begin{EQA}
	\zq_{2}(\BBH,\xx)
	& \eqdef &
	2 \vH \, \sqrt{\xx} + 2 \supA \xx \,\, .
\label{zqdefGQF}
\end{EQA}
%
\end{theorem}

\begin{proof}
W.l.o.g. assume that \( \supA = \| \BBH \| = 1 \).
We use the identity \( \langle \BB \xiv, \xiv \rangle = \langle \BBH \gaussv, \gaussv \rangle \) with
 \( \gaussv \sim \ND(0,\Id_{\dimp}) \).
We apply the exponential Chebyshev inequality: with \( \mu > 0 \)
\begin{EQA}
	\P\Bigl( \langle \BBH \gaussv, \gaussv \rangle > \zq^{2} \Bigr)
	& \leq &
	\E \exp \Bigl( \mu \langle \BBH \gaussv, \gaussv \rangle / 2 - \mu \zq^{2} / 2 \Bigr) \, .
\label{PBggiz2E2mz2}
\end{EQA}
Given \( \xx > 0 \), fix \( \mu < 1 \) by the equation
\begin{EQA}
	\frac{\mu}{1 - \mu} 
	&=&
	\frac{2 \sqrt{\xx}}{\vH} \, 
	\quad \text{ or } \quad
	\mu^{-1} 
	=
	1 + \frac{\vH}{2 \sqrt{\xx}} \, .
\label{1v2sxm12m1m}
\end{EQA}
Let \( \lambda_{j} \) be the eigenvalues of \( \BBH \), 
\( |\lambda_{j}| \leq 1 \).
It holds with \( \dimH = \tr \BBH \) in view of \eqref{m2v241m41mb}
\begin{EQA}
	&& \nquad
	\log \E \Bigl\{ \frac{\mu}{2} \bigl( \langle \BBH \gaussv, \gaussv \rangle - \dimH \bigr) \Bigr\}
	\leq 
	\frac{\mu^{2} \vH^{2}}{4 (1 - \mu)} \, .
\label{m2v241m4mj1p}
\end{EQA}
It remains to check that the choice \( \mu \) by \eqref{1v2sxm12m1m} and 
\( \zq = \zq(\BBH,\xx) \) yields
\begin{EQA}
	\frac{\mu^{2} \vH^{2}}{4 (1 - \mu)} - \frac{\mu (\zq^{2} - \dimH)}{2}
	& = &
	\frac{\mu^{2} \vH^{2}}{4 (1 - \mu)} - \mu \bigl( \vH \sqrt{\xx} + \xx \bigr)
	=
	\mu \Bigl( \frac{\vH \sqrt{\xx}}{2} - \vH \sqrt{\xx} - \xx \Bigr)
	=
	- \xx
	\qquad
	\qquad
\label{m2vA241muz2}
\end{EQA}
as required.
The last statement \eqref{PxivTBBdimA2vp} is obtained by applying this inequality twice to \( \BBH \) and \( - \BBH \).
\end{proof}

\begin{corollary}
\label{CTexpbLGAd}
Assume the conditions of Theorem~\ref{TexpbLGA}.
Then for \( \zq > \vH \)
\begin{EQA}
	\P\bigl( \bigl| \langle \BB \xiv, \xiv \rangle - \dimH \bigr| \ge \zq \bigr)
	& \leq &
	2 \exp\biggl\{ - \frac{\zq^{2}}{\bigl( \vH + \sqrt{\vH^{2} + 2 \supA \zq} \bigr)^{2}} \biggr\}
	\\
	& \leq & 
	2 \exp\biggl( - \frac{\zq^{2}}{4\vH^{2} + 4 \supA \zq} \biggr) .
	\qquad
\label{3z2spsp2z3z2}
\end{EQA}
\end{corollary}

\begin{proof}
Given \( \zq \), define \( \xx \) by 
\( 2 \vH \sqrt{\xx} + 2 \supA \xx = \zq \) or 
\( 2 \supA \sqrt{\xx} = \sqrt{\vH^{2} + 2 \supA \zq} - \vH \).
Then
\begin{EQA}
	\P\bigl( \langle \BB \xiv, \xiv \rangle - \dimH \ge \zq \bigr)
	& \leq &
	\ex^{-\xx} 
	=
	\exp\biggl\{ - \frac{\bigl( \sqrt{\vH^{2} + 2 \supA \zq} - \vH \bigr)^{2}}{4 \supA^{2}} \biggr\}
	=
	\exp\biggl\{ - \frac{\zq^{2}}{\bigl( \vH + \sqrt{\vH^{2} + 2 \supA \zq} \bigr)^{2}} \biggr\}.
\label{3emzmsp22z2c}
\end{EQA}
This yields \eqref{3z2spsp2z3z2} by direct calculus.
\end{proof}

Of course, bound \eqref{3z2spsp2z3z2} is sensible only if \( \zq \gg \vH \).

\begin{corollary}
\label{RsochpHsA}
Assume the conditions of Theorem~\ref{TexpbLGA}.
If also \( \BB \geq 0 \), then 
\begin{EQA}
\label{Pxiv2dimAxx12}
	\P\Bigl( \langle \BB \xiv, \xiv \rangle \geq \zq^{2}(\BB,\xx) \Bigr)
	& = &
	\P\bigl( \| \BB^{1/2} \xiv \| \geq \zq(\BB,\xx) \bigr)
	\leq 
	\ex^{-\xx} 
\end{EQA}
with 
\begin{EQA}
	\zq(\BB,\xx)
	& \eqdef &
	\sqrt{\dimH + 2 \vH \, \sqrt{\xx} + 2 \supA \xx} 
	\leq 
	\sqrt{\dimH} + \sqrt{2 \supA \xx} \, .
\label{zzxxppdBlroBB}
\end{EQA}
Also
\begin{EQA}
	\P\Bigl( \langle \BB \xiv, \xiv \rangle - \dimH < - 2 \vH \, \sqrt{\xx} \Bigr)
	& \leq &
	\ex^{-\xx} .
\label{Pxiv2dimAvp12d}
\end{EQA}
\end{corollary}

\begin{proof}
The definition implies \( \vH^{2} \leq \dimH \supA \).
One can use a sub-optimal choice of the value 
\( \mu(\xx) = \bigl\{ 1 + 2 \sqrt{\supA \dimH/\xx} \bigr\}^{-1} \) yielding the statement of the corollary.
\end{proof}

As a special case, we present a bound for the chi-squared distribution 
corresponding to \( \BB = \HVB^{2} = \Id_{\dimp} \), \( \dimp < \infty \).
Then \( \tr (\BBH) = \dimp \), \( \tr(\BBH^{2}) = \dimp \) and \( \supA(\BBH) = 1 \).

\begin{corollary}
\label{Cchi2p}
Let \( \gaussv \) be a standard normal vector in \( \R^{\dimp} \).
Then for any \( \xx > 0 \)
\begin{EQA}[ccl]
\label{Pxi2pm2px}
	\P\bigl( \| \gaussv \|^{2} \geq \dimp + 2 \sqrt{\dimp \, \xx} + 2 \xx \bigr)
	& \leq &
	\ex^{-\xx},
	\\
	\P\bigl( \| \gaussv \| \,\,  \geq \sqrt{\dimp} + \sqrt{2 \xx} \bigr)
	& \leq &
	\ex^{-\xx} ,
\label{Pxi2pm2px12}
	\\
	\P\bigl( \| \gaussv \|^{2} \leq \dimp - 2 \sqrt{\dimp \, \xx} \bigr)
	& \leq &
	\ex^{-\xx}	.
\label{Pxi2pm2px22}
\end{EQA}
\end{corollary}

The bound of Theorem~\ref{TexpbLGA} 
can be represented as a usual deviation bound.

\begin{theorem}
\label{CTexpbLGA}
Assume the conditions of Theorem~\ref{TexpbLGA} with \( \BB \geq 0 \).
Then for \( \zq > \sqrt{\dimH} + 1 \)
\begin{EQA}
	\P\bigl( \langle \BB \xiv, \xiv \rangle \ge \zq^{2} \bigr)
	& \leq &
	\exp\Bigl\{ - \frac{(\zq - \sqrt{\dimH})^{2}}{2 \supA} \Bigr\} ,
\label{3emzmsp22z2}
	\\
	\E \bigl\{ \langle \BB \xiv, \xiv \rangle^{1/2} \Ind\bigl( \langle \BB \xiv, \xiv \rangle \ge \zq^{2} \bigr) \bigr\}
	& \leq &
	\exp\Bigl\{ - \frac{(\zq - \sqrt{\dimH})^{2}}{2 \supA} \Bigr\} ,
\label{3emzmsp22z21}
	\\
	\E \bigl\{ \langle \BB \xiv, \xiv \rangle \Ind\bigl( \langle \BB \xiv, \xiv \rangle \ge \zq^{2} \bigr) \bigr\}
	& \leq &
	\frac{2 \zq}{\zq - \sqrt{\dimH}} \exp\Bigl\{ - \frac{(\zq - \sqrt{\dimH})^{2}}{2 \supA} \Bigr\} .
\label{3emzmsp22z2e}
\end{EQA}
\end{theorem}

\begin{proof}
Bound \eqref{3emzmsp22z2} follows from 
\eqref{Pxiv2dimAxx12}.
It obviously suffices to check the bound for the excess risk for \( \supA = 1 \).
It follows with \( \eta = \| \BB^{1/2} \xiv \| \) for \( \zq \geq \sqrt{\dimH} + 1 \)
\begin{EQA}
	\E \bigl\{ \eta \Ind(\eta > \zq) \bigr\}
	&=&
	\int_{\zq}^{\infty} \P(\eta \geq \zq) \, d\zq
	\leq 
	\int_{\zq}^{\infty} \exp\bigl\{ - \frac{(x - \sqrt{\dimH})^{2}}{2} \bigr\} \, dx
	\leq 
	\exp\bigl\{ - \frac{(\zq - \sqrt{\dimH})^{2}}{2} \bigr\}.
\label{zEe2Iezz2c2H23}
\end{EQA} 
Similarly
\begin{EQA}
	\E \bigl\{ \eta^{2} \Ind(\eta^{2} > \zq^{2}) \bigr\}
	&=&
	\int_{\zq^{2}}^{\infty} \P(\eta^{2} \geq \zz) \, d\zz
	\leq 
	\int_{\zq^{2}}^{\infty} \exp\bigl\{ - \frac{(\sqrt{\zz} - \sqrt{\dimH})^{2}}{2} \bigr\} \, d\zz .
\label{zEe2Iezz2c2H23d}
\end{EQA} 
By change of variables \( \sqrt{\zz} - \sqrt{\dimH} = u \) for \( \zq > \sqrt{\dimH} + 1 \)
\begin{EQA}
	&& \nquad
	\int_{\zq^{2}}^{\infty} \exp\bigl\{ - \frac{(\sqrt{\zz} - \sqrt{\dimH})^{2}}{2} \bigr\} \, d\zz
	\leq 
	2 \int_{\zq - \sqrt{\dimH}}^{\infty} (u + \sqrt{\dimH}) \, \exp\{ - u^{2}/2 \} \, du
	\\
	& \leq &
	2 \left( 1 + \frac{\sqrt{\dimH}}{\zq - \sqrt{\dimH}} \right) 
	\exp\bigl\{ - (\zq - \sqrt{\dimH})^{2}/2 \bigr\} 
	=
	\frac{2 \zq}{\zq - \sqrt{\dimH}} \exp\bigl\{ - (\zq - \sqrt{\dimH})^{2}/2 \bigr\}\, .
\label{dz21fsHzsHe22}
\end{EQA}
\end{proof}

\bibliography{exp_ts,listpubm-with-url}

\end{document}